\documentclass{amsart}

\usepackage{amssymb}
\usepackage{graphicx}

\newtheorem{theorem}{Theorem}[section]
\newtheorem{lemma}[theorem]{Lemma}
\newtheorem{proposition}[theorem]{Proposition}

\newtheorem{theorem_definition}[theorem]{Theorem-Definition}
\newtheorem{corollary}[theorem]{Corollary}

\newtheorem{mytheorem}{Theorem}
\newtheorem{mycorollary}[mytheorem]{Corollary}

\newtheorem{_definition}[theorem]{Definition}
\newenvironment{definition}{\begin{_definition}\rm}{\end{_definition}}

\newtheorem{_example}[theorem]{Example}
\newenvironment{example}{\begin{_example}\rm}{\end{_example}}

\newtheorem{_remark}[theorem]{\it Remark}
\newenvironment{remark}{\begin{_remark}\rm}{\end{_remark}}

\newtheorem{_claim}{Claim}

\numberwithin{equation}{subsection}
\numberwithin{table}{subsection}
\numberwithin{figure}{subsection}
\numberwithin{theorem}{subsection}

\newcommand{\A}{\mathord{\mathbb A}}
\newcommand{\C}{\mathord{\mathbb C}}
\newcommand{\E}{\mathord{\mathbb E}}
\newcommand{\F}{\mathord{\mathbb F}}

\renewcommand{\P}{\mathord{\mathbb  P}}
\newcommand{\Q}{\mathord{\mathbb  Q}}
\newcommand{\R}{\mathord{\mathbb R}}
\newcommand{\Z}{\mathord{\mathbb Z}}

\newcommand{\AAA}{\mathord{\mathcal A}}
\newcommand{\CCC}{\mathord{\mathcal C}}
\newcommand{\DDD}{\mathord{\mathcal D}}
\newcommand{\EEE}{\mathord{\mathcal E}}
\newcommand{\FFF}{\mathord{\mathcal F}}
\newcommand{\GGG}{\mathord{\mathcal G}}
\newcommand{\HHH}{\mathord{\mathcal H}}
\newcommand{\III}{\mathord{\mathcal I}}

\newcommand{\LLL}{\mathord{\mathcal L}}
\newcommand{\MMM}{\mathord{\mathcal M}}

\newcommand{\OOO}{\mathord{\mathcal O}}
\newcommand{\PPP}{\mathord{\mathcal P}}
\newcommand{\QQQ}{\mathord{\mathcal Q}}
\newcommand{\RRR}{\mathord{\mathcal R}}
\newcommand{\SSS}{\mathord{\mathcal S}}
\newcommand{\TTT}{\mathord{\mathcal T}}

\newcommand{\WWW}{\mathord{\mathcal W}}
\newcommand{\XXX}{\mathord{\mathcal X}}
\newcommand{\YYY}{\mathord{\mathcal Y}}
\newcommand{\ZZZ}{\mathord{\mathcal Z}}

\font\mathgot=eufm10

\newcommand{\mmm}{\mathord{\hbox{\mathgot m}}}
\newcommand{\ppp}{\mathord{\hbox{\mathgot p}}}

\font\smallmathgot=eufm8
\newcommand{\smallppp}{\mathord{\hbox{\smallmathgot p}}}
\newcommand{\smallmmm}{\mathord{\hbox{\smallmathgot m}}}

\newcommand{\maprightsp}[1]{\; \smash{\mathop{\; \longrightarrow \; }\limits\sp{#1}}\; }

\newcommand{\mapleftsp}[1]{\; \smash{\mathop{\; \longleftarrow \; }\limits\sp{#1}}\; }

\newcommand{\inj}{\hookrightarrow}
\newcommand{\leftinj}{\hookleftarrow}

\newcommand{\isom}{\mathbin{\,\raise -.6pt\rlap{$\to$}\raise 3.5pt%
\hbox{\hskip .3pt$\mathord{\sim}$}\,}}

\newcommand{\set}[2]{\{\; {#1} \; \mid \; {#2} \;  \}}
\newcommand{\shortset}[2]{\{ {#1} \,|\, {#2}   \}}

\newcommand{\map}[3]{ #1 \, : \, #2 \, \to \, #3}

\newcommand{\mapisomet}[3]{#1\;:\; #2\;\isomet\;#3}

\newcommand{\shortmap}[3]{ #1  :  #2 \to #3}

\newcommand{\sprime}{\sp\prime}
\newcommand{\spprime}{\sp{\prime\prime}}
\newcommand{\sptimes}{\sp{\times}}
\newcommand{\sperp}{\sp{\perp}}

\newcommand{\dual}{\sp{\vee}}

\newcommand{\inv}{\sp{-1}}

\newcommand{\NS}{\mathord{\rm NS}}

\newcommand{\Ker}{\operatorname{\rm Ker}\nolimits}
\newcommand{\Coker}{\operatorname{\rm Coker}\nolimits}
\renewcommand{\Im}{\operatorname{\rm Im}\nolimits}

\newcommand{\Gal}{\operatorname{\rm Gal}\nolimits}
\newcommand{\End}{\operatorname{\rm End}\nolimits}
\newcommand{\pr}{\operatorname{\rm pr}\nolimits}
\newcommand{\Hom}{\operatorname{\rm Hom}\nolimits}
\newcommand{\Spec}{\operatorname{\rm Spec}\nolimits}
\newcommand{\rank}{\operatorname{\rm rank}\nolimits}
\newcommand{\disc}{\operatorname{\rm disc}\nolimits}
\newcommand{\ch}{\operatorname{\rm char}\nolimits}
\newcommand{\GL}{\operatorname{\it GL}\nolimits}
\newcommand{\SL}{\operatorname{\it SL}\nolimits}
\newcommand{\Km}{\operatorname{\rm Km}\nolimits} 
\newcommand{\Emb}{\operatorname{\rm Emb}\nolimits} 
\newcommand{\Pic}{\operatorname{\rm Pic}\nolimits}
\newcommand{\Lie}{\operatorname{\rm Lie}\nolimits}
\newcommand{\Tr}{\operatorname{\rm Tr}\nolimits}
\newcommand{\Leff}{\operatorname{\rm Leff}\nolimits}
\newcommand{\Frob}{\operatorname{\rm Frob}\nolimits}

\newcommand{\wt}[1]{\widetilde{#1}}

\newcommand{\ang}[1]{\langle #1\rangle}

\newcommand{\rmand}{\textrm{and}}
\newcommand{\quand}{\quad\rmand\quad}

\newcommand{\weight}{\operatorname{\rm wt}\nolimits}

\newcommand{\cohomH}{\mathord{\rm H}}
\newcommand{\cohom}{\mathord{\rm H}^2}
\newcommand{\Cl}{\mathord{\it Cl}}

\newcommand{\sbsqppp}{_{[\smallppp]}}
\newcommand{\sbppp}{\sb{\smallppp}}

\newcommand{\isomet}{\inj}
\newcommand{\SSps}{\SSS}

\newcommand{\lats}{\LLL}
\newcommand{\speven}{\sp{\textrm{even}}}
\newcommand{\sppos}{\sp{\textrm{pos}}}

\newcommand{\Disc}[1]{D_{#1}}
\newcommand{\Discform}[1]{q_{#1}}

\newcommand{\RS}{\rm{RS}}

\newcommand{\id}{\mathord{\rm{id}}}

\newcommand{\notdiv}{{\not|}\,}

\newcommand{\primitive}{\sp{*}}

\newcommand{\ori}[1]{\widetilde{#1}}

\newcommand{\lat}{\Lambda}
\newcommand{\RSlat}{\Lambda}

\newcommand{\nr}{\operatorname{\rm nr}\nolimits}

\newcommand{\Bl}[1]{\widetilde{#1}} 

\newcommand{\DRS}[1]{(\Disc{#1}^{\RS}, \Discform{#1}^{\RS})}

\newcommand{\tiliota}{\tilde\iota} 
\newcommand{\phiA}{\varphi_{A}} 
\newcommand{\phiY}{\varphi_{Y}}

\newcommand{\KKmm}{\mathord{\mathcal{K}}}
\newcommand{\SSII}{\mathord{\mathcal{SI}}}
\newcommand{\SSIIKK}{\mathord{\mathcal{SIK}}}
\newcommand{\ELL}{\mathord{\mathcal{ELL}}}

\newcommand{\ZK}{\Z_K}

\newenvironment{rmenumerate}
{\begin{enumerate}
\renewcommand{\labelenumi}{(\theenumi)}
\renewcommand{\theenumi}{\roman{enumi}}}
{\end{enumerate}}

\newcommand{\an}{\mathord{\rm{an}}}
\newcommand{\ana}{{}^{\an}\hskip -1pt}

\newcommand{\cloQ}{\overline{\Q}}
\newcommand{\cloF}{\overline{\F}}

\newcommand{\wcong}{\;\cong\;}
\newcommand{\Ea}{E_{\alpha}}
\newcommand{\YYYgamma}{\YYY^{\,\gamma}}

\newcommand{\conj}[1]{\overline{#1}}

\newcommand{\smallmath}[1]{\raisebox{-2pt}{\hbox{${}^{#1}$}} }

\newcommand{\sqD}{\sqrt{D}}
\newcommand{\setR}{\mathord{\hbox{\bf R}}}
\newcommand{\quaternion}{B}
\newcommand{\NB}{N_B}

\newcommand{\typeI}{\mathord{\rm I}}
\newcommand{\typeIV}{\mathord{\rm IV}}

\newcommand{\textsum}{\textstyle{\sum}}

\newcommand{\fMor}{\mathord{\hbox{\mathgot Mor}}}
\newcommand{\Mor}{\mathord{\hbox{\rm Mor}}}
\newcommand{\Hilb}{\mathord{\hbox{\rm Hilb}}}
\newcommand{\fHom}{\mathord{\mathcal Hom}}
\newcommand{\Sch}{\mathord{\mathcal Sch}}

\begin{document}

\newcommand{\upinj}{\raisebox{-3pt}{\rotatebox{90}{$\inj$}}}
\newcommand{\downinj}{\raisebox{8pt}{\rotatebox{-90}{$\inj$}}}
\newcommand{\vertcong}{\raisebox{6pt}{\rotatebox{-90}{$\cong$}}} 
\newcommand{\NEinj}{\raisebox{1pt}{\rotatebox{45}{$\inj$}}}
\newcommand{\NWinj}{\raisebox{1pt}{\rotatebox{135}{$\inj$}}}
\newcommand{\SEinj}{\raisebox{10pt}{\rotatebox{-45}{$\inj$}}}
\newcommand{\SWinj}{\raisebox{10pt}{\rotatebox{-135}{$\inj$}}}
\newcommand{\NEisom}{\raisebox{-2pt}{\rotatebox{45}{$\isom$}}}
\newcommand{\NWisom}{\raisebox{-2pt}{\rotatebox{135}{$\isom$}}}
\newcommand{\SEisom}{\raisebox{6pt}{\rotatebox{-45}{$\isom$}}}
\newcommand{\SWisom}{\raisebox{6pt}{\rotatebox{-135}{$\siom$}}}

\title[Lattices of a singular $K3$ surface]{Transcendental lattices and supersingular reduction lattices
of a singular $K3$ surface}

\author{Ichiro Shimada}
\address{
Department of Mathematics,
Faculty of Science,
Hokkaido University,
Sapporo 060-0810,
JAPAN
}
\email{shimada@math.sci.hokudai.ac.jp
}


\begin{abstract}
A $K3$ surface $X$ defined over a field $k$ 
of characteristic $0$ is called \emph{singular} if the N\'eron-Severi lattice
$\NS (X)$ of $X\otimes \overline{k}$ is of rank $20$. Let  $X$ 
be a singular $K3$ surface defined over a number field $F$.
For each embedding $\sigma: F\inj \C$, we denote by $T(X\sp\sigma)$  the 
transcendental lattice of the complex $K3$ surface $X\sp\sigma$ obtained from $X$ 
by $\sigma$. For each prime  $\smallppp$ of $F$ 
at which $X$ has a supersingular reduction $X\sbppp$, 
we define $L(X, \smallppp)$ to be the orthogonal complement of  $\NS(X)$ in $\NS(X\sbppp)$.
We investigate the relation between these lattices $T(X\sp\sigma)$ and $L(X, \smallppp)$.
As an application,
we give a lower bound for the degree of a number field over which
a singular $K3$ surface with a given transcendental lattice can be defined.
\end{abstract}

\maketitle

\section{Introduction}\label{sec:Introduction}
For a smooth projective surface $X$ defined over a field $k$,
we denote by $\Pic (X)$ the Picard group of $X$,
and by $\NS (X)$ the  N\'eron-Severi lattice of $X\otimes \bar k$,
where $\bar k$ is the algebraic closure of $k$.
When  $X$ is a $K3$ surface,
we have a natural isomorphism $\Pic (X\otimes \bar k)\cong \NS(X)$.
We say that a $K3$ surface $X$ in characteristic $0$ is \emph{singular}
if $\NS (X)$ is of rank $20$,
while 
a $K3$ surface $X$ in characteristic $p>0$ is \emph{supersingular}
if $\NS (X)$ is of rank $22$.
It is known (\cite{MR578868}, \cite{MR0441982}, \cite{Shioda_murre}) that 
every complex singular $K3$ surface
is defined over a number field.
\par
\smallskip
For a number field $F$, we denote  
by $\Emb (F)$ the set of embeddings of $F$ into $\C$,
by $\Z_F$ the integer ring of $F$, and 
by $\pi_F: \Spec \Z_F\to\Spec \Z$ the natural projection.
%
%
%
%
Let $X$ be a singular $K3$ surface defined over a number field $F$,
and let $\XXX\to U$ be a smooth proper family of $K3$ surfaces
over   a non-empty open subset $U$ of $\Spec \Z_F$
such that the generic fiber is isomorphic to $X$.
We put
$$
d(X)\;:=\;\disc (\NS (X)).
$$
Remark that 
we have $d(X)<0$ by Hodge index theorem.
For $\sigma\in \Emb (F)$, we denote by $X\sp{\sigma}$ 
the complex analytic $K3$ surface obtained from $X$ by $\sigma$.
The \emph{transcendental lattice} $T(X\sp{\sigma})$ of $X\sp{\sigma}$
is defined to be the orthogonal complement of $\NS(X)\cong \NS(X\sp{\sigma})$
in the second Betti cohomology group $\cohom (X\sp{\sigma}, \Z)$,
which we regard  as a lattice by the cup-product.
Then  $T(X\sp{\sigma})$ is an even positive-definite lattice of rank $2$
with discriminant $-d(X)$.
For a closed point $\ppp$ of $U$,
we denote by $X\sbppp$ the reduction of $\XXX$ at $\ppp$.
Then $X\sbppp$  is a $K3$ surface defined over 
the  finite field $\kappa\sbppp:=\Z_F/\ppp$.
For a  prime integer $p$, we put
$$
\SSps_p (\XXX)\;:=\;\set{\ppp\in U}{\textrm{$\pi_F(\ppp)=p$ and $X\sbppp$ is supersingular}}.
$$
For each $\ppp\in \SSps_p (\XXX) $, we have the specialization homomorphism
$$
\map{\rho\sbppp }{ \NS(X)}{\NS(X\sbppp)},
$$
which preserves the intersection pairing
(see~\cite[Exp.~X]{MR0354655}, \cite[\S4]{MR0404257} or \cite[\S20.3]{MR1644323}),
and hence is injective.
We denote by $L(\XXX,\ppp)$ the orthogonal complement of $\NS(X)$ in $\NS(X\sbppp)$,
and call $L(\XXX,\ppp)$
the \emph{supersingular reduction lattice} of $\XXX$ at $\ppp$.
Then  $L(\XXX,\ppp)$ is an even negative-definite lattice of rank $2$.
We will see that, 
if $p\notdiv 2d(X)$, 
then the discriminant of $L(\XXX,\ppp)$ is $-p^2 d(X)$.
\par
\smallskip
For an odd prime integer $p$ not dividing  $x\in \Z$, we denote by
$$
\chi_p (x)\;:=\;\left(\frac{x}{\;p\;}\right)\;\;\in\;\;\{1, -1\}
$$
the Legendre character.
In~\cite[Proposition 5.5]{normalK3}, we have proved the following.
(See Theorem-Definition~\ref{thm:geomRSlat} for the definition of the Artin invariant.)

\begin{proposition}\label{prop:normalK3}
Suppose that $p\notdiv 2 d(X)$.
\par
{\rm (1)}
If $\chi_p(d(X))=1$, then $\SSps_p (\XXX)$ is empty.
\par
{\rm (2)}
If $\ppp \in \SSps_p (\XXX)$, then the Artin invariant of $X\sbppp$ is $1$.
\end{proposition}
The first main result of this paper,
which will be proved in \S\ref{subsec:thm:almostallp},  is as follows:
\begin{mytheorem}\label{thm:almostallp}
There exists a finite set $N$ of  prime integers 
containing the prime divisors of $2 d(X)$
such that the following holds:
\begin{equation}\label{eq:pnotinN1}
p\notin N\;\;\Rightarrow\;\;
\SSps_p(\XXX)=\begin{cases}
\emptyset   & \textrm{if $\chi_p(d(X))=1$,}\\
\pi_F\inv (p) & \textrm{if $\chi_p(d(X))=-1$.}
\end{cases}
\end{equation}
\end{mytheorem}
We put $\Z_\infty:=\R$.
Let $R$ be $\Z$ or $\Z_l$,
where $l$ is a prime integer or $\infty$.
An \emph{$R$-lattice} is a free $R$-module $\lat$ 
of finite rank with a non-degenerate symmetric bilinear form
$$
\map{(\phantom{b}, \phantom{b})}{ \lat\times \lat}{ R}.
$$
The \emph{discriminant} $\disc (\lat) \in R/(R\sptimes )^2$
of an $R$-lattice $\lat$
is the determinant modulo $(R\sptimes )^2$
of a symmetric matrix expressing $(\phantom{i}, \phantom{i})$.
\par
\smallskip
A $\Z$-lattice is simply called a \emph{lattice}.
For a lattice $\lat$ and a non-zero integer $n$,
we denote by $\lat [n]$ the lattice 
obtained from $\lat$ by multiplying 
the symmetric bilinear form $(\phantom{i}, \phantom{i})$ 
by $n$.
A lattice $\lat$ is said to be  \emph{even} if $(v,v)\in 2\Z$ holds for any $v\in \lat$.
Let $\lat$ and $\lat\sprime$ be lattices.
We denote by $\lat\perp\lat\sprime$ the orthogonal direct-sum of $\lat$ and $\lat\sprime$.
A homomorphism $\lat\to\lat\sprime$
preserving the symmetric bilinear form
is called an \emph{isometry}.
Note that an isometry is injective
because of  the non-degeneracy  of the symmetric bilinear forms.
An isometry $\lat\isomet \lat\sprime$
(or a \emph{sublattice} $\lat$ of  $\lat\sprime$)
is said to be \emph{primitive}
if the cokernel $\lat\sprime/\lat$ is torsion-free.
The \emph{primitive closure}
of a sublattice $\lat\isomet \lat\sprime$
is the intersection of $\lat\otimes \Q$ and $\lat\sprime$ in $\lat\sprime\otimes\Q$.
For an isometry $\lat\isomet\lat\sprime$,
we put
$$
(\lat\isomet\lat\sprime)\sperp\;:=\;
\set{x\in \lat\sprime}{(x, y)=0\quad\textrm{for all}\quad y\in \lat}.
$$
Note that $(\lat\isomet\lat\sprime)\sperp$ is  primitive in $\lat\sprime$.
Let $r$ be a positive integer, and $d$ a non-zero integer.
We denote by $\lats (r, d)$ the set of isomorphism classes of
lattices  of rank $r$ with discriminant $d$,
and by $[\lat]\in \lats (r, d)$ 
the isomorphism class of a lattice $\lat$.
If $[\lat]\in \lats (r, d)$, then we have 
$[\lat[n]]\in \lats (r, n^r d)$,
and the map $\lats (r, d)\to  \lats (r, n^r d)$
given by $[\lat]\mapsto[\lat[n]]$ is injective.
We denote by $\lats\speven (r,  d)$ (resp.~$\lats\sppos (r,  d)$)
the set of isomorphism classes in $\lats (r, d)$
of even lattices (resp.~of positive-definite lattices).
\par
\smallskip
We recall the notion of \emph{genera} of lattices.
See~\cite{MR522835}, for example, for details.
Two lattices $\lat$ and $\lat\sprime$ are said to be 
\emph{in the same genus}
if $\lat\otimes\Z_l$ and $\lat\sprime \otimes\Z_l$
are isomorphic as $\Z_l$-lattices for any $l$ (including $\infty$).
If  $\lat$ and $\lat\sprime$ are in the same genus,
then we have $\rank (\lat)=\rank (\lat\sprime)$
and $\disc (\lat)=\disc(\lat\sprime)$.
Therefore the set $\lats (r, d)$ is decomposed into the  disjoint union of
genera. 
For each non-zero integer $n$,
$\lat$ and $\lat\sprime$ are in the same genus
if and only if $\lat[n]$ and $\lat\sprime[n]$ are in the same genus.
Moreover,
if $\lat\spprime$ is in the same  genus as $\lat[n]$,
then there exists $\lat\sprime$
in the same genus as $\lat$
such that $[\lat\spprime]=[\lat\sprime[n]]$ holds.
Therefore, for each genus $\GGG\subset \lats (r, d)$,
we can define the genus $\GGG[n]\subset \lats (r, n^r d)$
by 
$$
\GGG[n]\;:=\;\set{[\lat[n]]}{[\lat]\in \GGG}.
$$
The map
from the set of genera in $\lats (r, d)$
to the set of genera in $\lats (r, n^r d)$
given by $\GGG\mapsto \GGG[n]$ is injective.
Suppose that $\lat$ and $\lat\sprime$ are in the same genus.
If $\lat$ is even (resp. positive-definite),
then so is $\lat\sprime$.
Hence   $\lats\speven (r,  d)$ and 
$\lats\sppos (r,  d)$  are 
also disjoint unions of
genera.
We  say that a genus $\GGG\subset \lats(r, d)$
is \emph{even} (resp.~\emph{positive-definite})
if $\GGG\subset \lats\speven (r,  d)$
(resp.~$\GGG\subset \lats\sppos(r,  d)$) holds.
\par
\smallskip
We review the theory of discriminant forms due to Nikulin~\cite{MR525944}.
Let $\lat$ be an even lattice.
We put $\lat\dual \;:=\;\Hom (\lat, \Z)$.
Then $\lat$ is embedded into $\lat\dual$ naturally as a submodule of finite index,
and there exists a unique $\Q$-valued symmetric bilinear form on $\lat\dual$
that extends the $\Z$-valued symmetric bilinear form on $\lat$.
We put
$$
\Disc{\lat}\;:=\;\lat\dual/\lat,
$$
which is a finite abelian group of order  $|\disc (\lat)|$,
and define a quadratic form
$$
\map{\Discform{\lat}}{\Disc{\lat}}{\Q/2\Z}
$$
by $\Discform{\lat} (x+\lat):=(x,x)+ 2\Z$ for $x\in \lat\dual$.
The finite quadratic form $(\Disc{\lat}, \Discform{\lat})$
is called the \emph{discriminant form} of $\lat$.
\begin{theorem_definition}[Corollary 1.9.4 in~\cite{MR525944}]
Let $\lat$ and $\lat\sprime$ be even lattices.
Then $\lat$ and $\lat\sprime$ are in the same genus
if and only if
the following hold{\rm :}
\begin{enumerate}
\renewcommand{\labelenumi}{(\theenumi)}
\renewcommand{\theenumi}{\roman{enumi}}
\item $\lat\otimes\Z_\infty$ and $\lat\sprime\otimes \Z_\infty$ are isomorphic as $\Z_\infty$-lattices,
and
\item the finite quadratic forms $(\Disc{\lat},\Discform{\lat})$ and $(\Disc{\lat\sprime},\Discform{\lat\sprime})$
are isomorphic.
\end{enumerate}
Therefore, for an even genus $\GGG$, 
we can define the \emph{discriminant form $(\Disc{\GGG}, \Discform{\GGG})$
of $\GGG$}.
\end{theorem_definition}
Next, we define Rudakov-Shafarevich lattices.
\begin{theorem_definition}[Section 1 of \cite{MR633161}]\label{thm:RSlat}
For each odd prime $p$ and a positive integer $\sigma\le 10$,
there exists, uniquely up to isomorphism,  an even lattice $\RSlat_{p, \sigma}$ of rank $22$
with signature $(1, 21)$
such that  the discriminant group is 
isomorphic to $(\Z/p\Z)^{\oplus 2\sigma}$.
We  call $\RSlat_{p, \sigma}$ a \emph{Rudakov-Shafarevich lattice}.
\end{theorem_definition}
%
%
%
%
\begin{theorem_definition}[Artin~\cite{MR0371899}, Rudakov-Shafarevich~\cite{MR633161}]\label{thm:geomRSlat}
For a supersingular $K3$ surface $Y$ in odd characteristic $p$,
there exists a positive integer $\sigma\le 10$,
which is called the \emph{Artin invariant of $Y$},
such that $\NS(Y)$ is isomorphic to the Rudakov-Shafarevich lattice $\RSlat_{p, \sigma}$.
\end{theorem_definition}
We denote by $\DRS{p, \sigma}$
the discriminant form  of the Rudakov-Shafarevich lattice $\RSlat_{p, \sigma}$.
The finite quadratic form $\DRS{p, \sigma}$
has been calculated explicitly in our previous paper~\cite[Proof of Proposition 4.2]{normalK3}.
\par
\smallskip
Our second main result, which will be 
proved in \S\ref{sec:thm:genus},  is as follows:
\begin{mytheorem}\label{thm:genus}
Let $X$ be a singular $K3$ surface defined over a number field $F$,
and let $\XXX\to U$ be a smooth proper family of $K3$ surfaces
over a non-empty open subset $U$ of $\Spec \Z_F$
such that the generic fiber is isomorphic to $X$.
We put $d(X):=\disc(\NS(X))$.
\par
\smallskip
{\rm (T)}
There exists a unique genus $\GGG_{\C}(X)\subset \lats(2, -d(X))$
such that $[\,T(X\sp{\sigma})\,]$ is contained in $\GGG_{\C}(X)$ for any $\sigma\in \Emb(F)$.
This  genus $\GGG_{\C}(X)$ is determined  by the properties that
it is even, positive-definite,  and that  the discriminant form
is isomorphic to $(\Disc{\NS(X)}, -\Discform{\NS(X)})$.
\par
\smallskip
{\rm (L)}
Let $p$ be a prime integer not dividing  $2d(X)$.
Suppose that $\SSps_p(\XXX)\ne \emptyset$.
Then there exists  a unique genus $\GGG_p (\XXX)\subset \lats(2, -d(X))$
such that $[L(\XXX, \ppp)]$ is contained in $\GGG_{p}(\XXX)[-p]$ for any $\ppp\in \SSps_p(\XXX)$.
This  genus $\GGG_p(\XXX)$ is determined  by the properties that
it is even, positive-definite,  and that
the discriminant form of $\GGG_p(\XXX)[-p]$ is isomorphic to 
$\DRS{p, 1}\oplus
(\Disc{\NS(X)}, -\Discform{\NS(X)})$.
\end{mytheorem}
To ease  notation, we put
$$
M[a,b,c]\;:=\;\left[\begin{matrix} 2a & b \\ b & 2c\end{matrix}\right].
$$
Let $D$ be a negative integer.
We then put
\begin{eqnarray}\label{eq:QQQ}
\phantom{aaaaaaa}\QQQ_D&:=&%
\set{M[a,b,c]}{%
a,b,c\in\Z,\;\; a>0,\;\; c>0,\;\; b^2-4ac=D},  \\
\phantom{aaaaaaa}\QQQ\primitive_D&:=&%
\set{M[a,b,c]\in \QQQ_D}{%
\gcd (a, b, c)=1}. \label{eq:QQQprim}
\end{eqnarray}
The group $\GL_2(\Z)$ acts on $\QQQ_D$ from right by
$(M, g)\mapsto g^T M g$ for $M\in \QQQ_D$ and $g\in \GL_2(\Z)$,
and the subset  $\QQQ\primitive_D$ of $\QQQ_D$ is stable by this action.
We put
$$
\renewcommand{\arraystretch}{1.2}
\begin{array}{ll}
\lats_D\;:=\;\QQQ_D/\GL_2(\Z),
\quad&
\lats\primitive_D\;:=\;\QQQ\primitive_D/\GL_2(\Z),
\quad\\
\ori{\lats}_D\;:=\;\QQQ_D/\SL_2(\Z),
\quad&
\ori{\lats}\primitive_D\;:=\;\QQQ\primitive_D/\SL_2(\Z).
\end{array}
$$
Then we have a natural identification
$$
\lats_D\;=\;\lats\sppos(2, -D)\cap \lats\speven(2, -D),
$$
and  $\ori{\lats}_D$
is regarded as the set of isomorphism classes
of even positive-definite \emph{oriented} lattices of rank $2$ with discriminant $-D$.
\par
\smallskip
Let $S$ be a complex $K3$ surface or a complex abelian surface.
Suppose that 
the transcendental lattice 
$T(S):=(\NS(S)\isomet \cohom(S,\Z))\sperp$ of $S$
 is of rank $2$.
Then  $T(S)$ is even, positive-definite and of discriminant $-d(S)$,
where $d(S):=\disc (\NS(S))$.
By the Hodge structure  
$$
T(S)\otimes \C=\cohomH^{2,0}(S)\oplus  \cohomH^{0,2}(S)
$$
of $T(S)$, 
we can  define a \emph{canonical orientation} on $T(S)$
as follows.
An ordered  basis $(e_1, e_2)$ of $T(S)$ is said to be \emph{positive}
if the imaginary part of
$(e_1, \omega_S)/(e_2, \omega_S)\in \C$
is positive, where $\omega_S$ is a basis of $\cohomH^{2,0}(S)$.
We denote by $\ori{T}(S)$ the \emph{oriented transcendental lattice} of $S$,
and by $[\,\ori{T}(S)\,]\in \ori{\lats}_{d(S)}$ the isomorphism class of $\ori{T}(S)$.
We have the following important theorem due to Shioda and Inose~\cite{MR0441982}:
\begin{theorem}[\cite{MR0441982}]\label{thm:SI}
The map $S\mapsto [\,\ori{T}(S)\,]$
gives rise to a bijection from the set of isomorphism classes 
of complex singular $K3$ surfaces $S$ 
to the set of isomorphism classes of 
even positive-definite oriented lattices of rank $2$.
\end{theorem}
If a genus 
$\GGG\subset\lats_D$ satisfies $\GGG\cap\lats\primitive_D\ne\emptyset$,
then $\GGG\subset \lats\primitive_D$ holds.
Therefore $\lats\primitive_D$ is a disjoint union of genera.
For a genus $\GGG\subset \LLL_D$, 
we denote by  $\ori{\GGG}$ the pull-back of $\GGG$ by the natural projection
$\ori{\lats}_{D}\to {\lats}_{D}$,
and call $\ori{\GGG}\subset \ori{\lats}_{D}$ a \emph{lifted genus}.
\par
\smallskip
A negative integer $D$ is called a \emph{fundamental discriminant}
if it is the discriminant of an imaginary quadratic field.
\par
\smallskip
Our third main result,
which  will be proved in \S\ref{subsec:thm:toFT} and  \S\ref{subsec:thm:toFL},
is as follows:
\par
\smallskip
\begin{mytheorem}\label{thm:toF}
Let $S$ be a complex singular $K3$ surface.
Suppose that
$D:=\disc (\NS(S))$
is a fundamental discriminant,
and that $[T(S)]$ is contained in $\lats\primitive _{D}$.
\par
\smallskip
{\rm (T)}
There exists a singular $K3$ surface $X$ defined over a number field $F$
such that 
$\shortset{[\,\ori{T} (X\sp{\sigma})\,]}{\sigma\in \Emb(F)}$
is equal to   the lifted genus in $\ori{\lats}\primitive _{D}$
that contains  $[\ori{T} (S)]$.
In particular, there exists $\sigma_0 \in \Emb(F)$ such that
$X^{\sigma_0}$ is isomorphic to $S$ over $\C$.
\par
\smallskip
{\rm (L)}
Suppose further that $D$ is odd.
Then there exists 
a smooth proper family 
$\XXX\to U$ of $K3$ surfaces over a non-empty open subset $U$ of $\Spec \Z_F$,
where $F$ is a number field, such that the following hold:
\begin{enumerate}
\renewcommand{\theenumi}{\roman{enumi}}
\item\label{item:L1} 
the generic fiber $X$ of $\XXX\to U$ satisfies the property in~{\rm{(T)}}  above,
\item\label{item:L2}
if $p\in \pi_F(U)$, then $p\notdiv 2D$, and 
\item\label{item:L3}
if $p\in \pi_F(U)$ and $\chi_p(D)=-1$, then
$\SSps_p(\XXX)=\pi_F\inv (p)$ holds, and  the set
$\shortset{[\,L(\XXX,\ppp)\,]}{\ppp\in \SSps_p(\XXX)}$
coincides with a genus  in $\lats(2, -p^2 D)$.
\end{enumerate}
\end{mytheorem}
Suppose that $D$ is a negative fundamental discriminant.
The set $\ori{\lats}\primitive_D$ and its decomposition into lifted genera are very well understood
by the work of Gauss.
We review the theory briefly.
We put $K:=\Q(\sqrt{D})$,
and denote by  $\III_D$ the multiplicative group  of non-zero fractional ideals of $K$,
by $\PPP_D\subset \III_D$ the subgroup of non-zero principal fractional ideals,
and by $\Cl_D:=\III_D/\PPP_D$ the ideal class group of  $K$.
Let $I$ be an element of $\III_D$. 
We denote by $[I]\in \Cl_D$ the ideal class of $I$.
We put
\begin{equation*}\label{eq:NI}
N(I)\;:=\;[\Z_K:nI]/n^2,
\end{equation*}
%
%
where $n$ is an integer $\ne 0$  such that $nI\subset\Z_K$, and 
define a bilinear form on $I$ by 
\begin{equation}\label{eq:bilK}
(x, y)\;:=\;(x\bar y + y \bar x)/ N(I)\;=\;\Tr_{K/\Q} (x\bar y)/ N(I).
\end{equation}
We say that 
an ordered basis $(\omega_1, \omega_2)$ of $I$ as a $\Z$-module
is positive if
\begin{equation}\label{eq:oriinK}
(\omega_1\bar \omega_2-\omega_2\bar \omega_1)/\sqrt{D}>0.
\end{equation}
By the bilinear form~\eqref{eq:bilK} and the orientation~\eqref{eq:oriinK}, 
the $\Z$-module $I$  of rank $2$ obtains a structure of 
even positive-definite  oriented lattice  with discriminant $-D$.
 The isomorphism class of this oriented lattice is 
 denoted by 
$\psi(I)\in \ori{\lats}_D$.
For the following,
see~\cite[Theorems 5.2.8 and 5.2.9]{MR1228206} and \cite[Theorem 3.15]{MR1028322}, for example.
\begin{proposition}\label{prop:ClD}
{\rm (1)}
The map $\psi: \III_D\to \ori{\lats}_D$ defined above induces a bijection 
$\Psi: \Cl_D\cong \ori{\lats}\primitive_D $
with the inverse     given by the following.
Let  $[\lat]\in \ori{\lats}\primitive_D$
be  represented by $M[a,b,c]\in \QQQ\primitive_D$,
and let  $I\in \III_D$ be 
the fractional ideal
generated by $\omega_1=(-b+\sqrt{D})/2$ and $\omega_2=a$.
Then $\Psi([I])=[\psi(I)]$
is equal to $[\lat]$.
\par
{\rm (2)}
Let $[I]$ and $[J]$ be elements of $\Cl_D$.
Then $\Psi([I])$ 
and $\Psi([J])$ are in the same lifted genus if and only if
$[I][J]\inv$ is contained in $\Cl_D^2:=\shortset{x^2}{x\in\Cl_D}$.
In particular,
every lifted genus in $\ori{\lats}\primitive_D$
consists of the same number of isomorphism classes,
and the cardinality is equal to $|\Cl_D^2|$.
\end{proposition}
%
%
%
%
%
%
%
%
Using Theorems~\ref{thm:SI} and~\ref{thm:toF}(T),
we obtain the following:
\begin{mycorollary}\label{cor:lowerbound}
Let $S$ be a complex singular  $K3$ surface such that 
$D:=\disc(\NS(S))$ is a fundamental discriminant,
and that 
$[T(S)]$ is contained in $\lats\primitive _{D}$.
Let $Y$ be  a $K3$ surface   defined over a number field $L$
such that $Y\sp{\tau_0}$ is isomorphic to $S$ over $\C$ for some $\tau_0\in \Emb (L)$.
Then we have $[L:\Q]\ge |\Cl_D^2|$.
\end{mycorollary}
\begin{proof}
Let $X$ be the $K3$ surface defined over a number field $F$ given
in Theorem~\ref{thm:toF}(T).
Then
the complex $K3$ surfaces
$X\sp{\sigma_0}$ and $Y\sp{\tau_0}$ are  isomorphic over $\C$,
and hence 
 there exists a number field $M\subset \C$ containing both of
$\sigma_0(F)$ and $\tau_0(L)$ such that $X\otimes M$ and $Y\otimes M$ are isomorphic over $M$.
Therefore, for each $\sigma\in \Emb (F)$,
there exists $\tau\in \Emb (L)$ such that $X\sp{\sigma}$ is  isomorphic to $Y\sp{\tau}$  over $\C$.
Since there exist exactly $|\Cl_D^2|$ 
isomorphism classes of complex $K3$ surfaces
among  $X\sp{\sigma}$ $({\sigma\in \Emb (F)})$,
we have  $|\Emb (L)|\ge |\Cl_D^2|$.
\end{proof}
The proof of Theorem~\ref{thm:genus}
is in fact an easy application of Nikulin's theory
of discriminant forms, and is given 
 in \S\ref{sec:thm:genus}.
The main tool of the proof of Theorems~\ref{thm:almostallp}  and~\ref{thm:toF}
is the Shioda-Inose-Kummer construction~\cite{MR0441982}.
This construction makes a singular $K3$ surface $Y$
from a pair of elliptic curves $E\sprime$ and $E$.
Shioda and Inose~\cite{MR0441982} proved
that, over $\C$,  the transcendental lattices of $Y$ 
and $E\sprime\times E$ are isomorphic.
We present their construction in our setting,
and show that,
over a number field,
the supersingular reduction lattices of
$Y$ and $E\sprime\times E$
are also isomorphic under certain assumptions.
The supersingular reduction lattice of $E\sprime\times E$
is calculated by 
the specialization homomorphism
$\Hom (E\sprime, E)\to \Hom(E\sprime\sbppp, E\sbppp)$.
In \S\ref{sec:A},
we investigate 
the $\Hom$-lattices of elliptic curves.
After examining
the Kummer construction in~\S\ref{sec:Kummer}
and the Shioda-Inose construction in~\S\ref{sec:ShiodaInose},
we prove Theorems~\ref{thm:almostallp} and~\ref{thm:toF}
in~\S\ref{sec:proofs}.
For Theorem~\ref{thm:toF}(T),
we use 
the  Shioda-Mitani theory~\cite{MR0382289}.
For Theorem~\ref{thm:toF}(L),
we need a description of 
embeddings of $\Z_K$ into maximal orders of
a quaternion algebra 
over $\Q$.
We use Dorman's description~\cite{MR1024555},
which we expound  in~\S\ref{sec:Endss}. 
\par
\smallskip
In~\cite{MR1465448},
Shafarevich studied, by means of  the Shioda-Inose-Kummer construction,  
number fields over which a singular $K3$ surface 
with a prescribed N\'eron-Severi  lattice 
can be  defined, 
and proved a certain finiteness theorem.
\par
\smallskip
The supersingular reduction lattices and their  relation to 
the transcendental lattice 
were first studied by Shioda~\cite{MR2023754}
for certain  $K3$ surfaces.
Thanks are due to  Professor Tetsuji Shioda
for stimulating conversations
and many  comments.
\par
\smallskip
After the first version of this paper appeared on the e-print archive,
Sch\"utt~\cite{MS} has succeeded in removing the assumptions in Theorem~\ref{thm:toF}(T)
and Corollary~\ref{cor:lowerbound}
that $D=\disc(\NS(S))$ be a fundamental discriminant,
and that $[T(S)]$ be in $\lats\primitive_D$.
Interesting examples 
of singular $K3$ surfaces defined over number fields
are also given in~\cite[\S7]{MS}.
\par
\smallskip
Applications to  topology of 
Theorem~\ref{thm:toF}(T) and its  generalization by Sch\"utt~\cite{MS}
are given  in~\cite{AZP} and~\cite{nonhomeo}.
\par
\smallskip
The author expresses gratitude  to the referee for many comments and suggestions improving the exposition.
\par
\smallskip
%
%
Let $W$ be a Dedekind domain.
For $P\in \Spec W$, we put
\begin{equation}\label{eq:kappa}
\kappa_P\;:=\;\begin{cases}
\textrm{the quotient field of $W$} &\textrm{if $P$ is the generic point,}\\
W/\ppp &\textrm{if $P$ is a closed point $\ppp$.}
\end{cases}
\end{equation}
\section{Proof of Theorem~\ref{thm:genus}}\label{sec:thm:genus}
\subsection{The discriminant form of an orthogonal complement}
The following can be  derived from~\cite[Proposition 1.5.1]{MR525944}.
We give a simple and direct proof.
\begin{proposition}\label{prop:nikulin}
Let $L$ be an even lattice,  and  $M\subset L$   a primitive sublattice.
We put
$N:=(M\isomet L)\sperp$.
Suppose that $\disc  (M)$ and $\disc (L)$ are prime to each other.
Then there exists an isomorphism 
$$
(\Disc{N}, \Discform{N})\wcong (\Disc{L}, \Discform{L})\oplus (\Disc{M}, -\Discform{M})
$$
of finite quadratic forms.
In particular, we have
$\disc (N)=\disc (L)\disc (M)$.
\end{proposition}
\begin{proof}
We put 
$d_L:=|\disc (L)|=|D_L|$.
The multiplication by $d_L$ 
induces  an automorphism $\delta_L: \Disc{M}\isom \Disc{M}$
of $\Disc{M}$ by the assumption.
We regard  $L$, $M$, $N$ and $L\dual$, $M\dual$, $N\dual$ as submodules of
$L\otimes \Q=(M\otimes \Q)\oplus (N\otimes \Q)$.
First we show that
\begin{equation}\label{eq:cap}
L\dual\cap M\dual =M.
\end{equation}
The inclusion $\supseteq$ is obvious.
Suppose that $x\in L\dual \cap M\dual$.
Then we have $d_L x\in L$.
Since $M$ is primitive in $L$,
we have $L\cap M\dual =M$, and  hence $\delta_L (x+M)=0$ holds in $\Disc{M}$.
Because $\delta_L$ is an automorphism of $\Disc{M}$, we have $x\in M$.
Next we show that the composite of natural homomorphisms 
\begin{equation}\label{eq:composite}
L\;\;\inj\;\; L\dual \;\to\; M\dual \;\to\; \Disc{M}
\end{equation}
is surjective.
Let $\xi \in \Disc{M}$ be given.
There exists $\eta\in \Disc{M}$ such that $\delta_L(\eta)=\xi$.
Since $L\dual\to M\dual$ is surjective by the primitivity of $M\isomet L$,
there exists $y\in L\dual $ that is mapped  to $\eta$.
Then $x:=d_L y$ is in $L$ and is mapped to $\xi$.
\par
We define a homomorphism
$\tau : \Disc{N}\to \Disc{L}\oplus \Disc{M}$ as follows.
Let $x\in N\dual$ be given.
Since $L\dual\to N\dual$ is surjective by the primitivity of $N\isomet L$,
there exists $z\in L\dual $ that is mapped to $x$.
Let $y\in M\dual$ be the image of $z$ by $L\dual\to M\dual$.
We put
$$
\tau (x+N)\;:=\;(z+L, y+M).
$$
The well-definedness of $\tau$ follows from the formula~\eqref{eq:cap}.
Since  $z=(y, x)$ in $L\dual\subset M\dual\oplus N\dual$,
we have $q_N(x+N)=q_L(z+L)-q_M(y+M)$.
The infectivity of $\tau$ follows from $L\cap N\dual =N$.
Since the homomorphism~\eqref{eq:composite} is surjective,
the homomorphism $\tau$ is also surjective.
\end{proof}
\subsection{The cokernel of the  specialization isometry}\label{subsec:coker}
Let $W$ be a Dedekind domain with the quotient field $F$ being a number field,
and let $\XXX\to U:=\Spec W$ be a smooth proper family of $K3$ surfaces.
We put $X:=\XXX\otimes F$.
In this subsection, we do \emph{not} assume that  $\rank (\NS (X))=20$.
Let $\ppp$ be a closed point of $U$ such that $X_0:=\XXX\otimes\kappa\sbppp$ is 
 supersingular.
We consider the specialization isometry
$$
\rho\;\;:\;\; \NS(X)=\Pic (X\otimes \overline{F})\;\; \isomet \;\; \NS(X_0)=\Pic (X_0\otimes \bar{\kappa}\sbppp),
$$
whose definition is given in \cite[Exp.~X]{MR0354655} or \cite[\S4]{MR0404257}.
We put $p:=\ch \kappa\sbppp$.
\begin{proposition}
Every torsion element of $\Coker(\rho)$ has order  a power of $p$.
\end{proposition}
\begin{proof}
We denote by $\hat{F}$ the completion of $F$ at $\ppp$,
and by $\hat{A}$ the valuation ring of $\hat{F}$
with the maximal ideal $\hat{\ppp}$.
Let $\hat{L}$ be a finite extension of $\hat{F}$ with the valuation ring $\hat{B}$,
the maximal ideal $\hat{\mmm}$, and the residue field   $\kappa_{\hat{\smallmmm}}$ such that
there exist  natural isomorphisms
$\Pic (X\otimes \hat{L})\cong \NS(X)$ and 
$\Pic (X_0\otimes \kappa_{\hat{\smallmmm}})\cong \NS(X_0)$.
Then $\rho$ is obtained from the restriction isomorphism
$$
\Pic (\XXX\otimes \hat{B}) \;\isom\;  \Pic (X\otimes \hat{L})
$$
to the generic fiber, 
whose inverse is given by taking the closure of divisors,
and the restriction homomorphism 
\begin{equation}\label{eq:resPic}
\Pic (\XXX\otimes \hat{B}) \;\to\;  \Pic (X_0\otimes \kappa_{\hat{\smallmmm}})
\end{equation}
to the central fiber.
Therefore it is enough to show that the order of  any   torsion element
of the cokernel of the homomorphism~\eqref{eq:resPic} is a power of $p$.
We put
$$
\YYY:=\XXX\otimes \hat{B}\quand Y_n:=\YYY\otimes (\hat{B}/\hat{\mmm}^{n+1}).
$$
Let $\hat{Y}$ be  the formal scheme obtained by completing $\YYY$
along $Y_0=X_0\otimes \kappa_{\hat{\smallmmm}}$.
Note that  $(\YYY, Y_0)$ satisfies the effective Lefschetz condition
$\Leff (\YYY, Y_0)$ in~\cite[Exp.~X]{MR0476737}.
(See~\cite[Theorem 9.7 in Chap.~II]{MR0463157}.)
Hence, by~\cite[Proposition 2.1 in Exp.~XI]{MR0476737},
we have
$\Pic (\YYY)\cong \Pic (\hat{Y})$.
On the other hand, we have
$\Pic (\hat{Y})= \projlim_n \Pic (Y_n)$
by~\cite[Exercise 9.6 in Chap.~II]{MR0463157}.
Let $\OOO_n$ denote the structure sheaf of $Y_n$.
From the natural exact sequence 
$0 \to  \OOO_0 \to  \OOO_{n+1}\sptimes \to  \OOO_n\sptimes \to  1$
(see~\cite[Exp.~XI]{MR0476737}),
we obtain an exact sequence
$$
0\;\to\;  \Pic (Y_{n+1})\;\to\;  \Pic (Y_{n})\;\to\;  \cohom(Y_0, \OOO_0).
$$
In particular, the projective limit of $\Pic (Y_n)$ is equal to $\cap_n \Pic (Y_{n})$.
Since every  non-zero element of $\cohom(Y_0, \OOO_0)$ is of order $p$,
every torsion element of $\Pic (Y_0)/\cap_n \Pic (Y_{n})$
is of order a power of $p$.
\end{proof}
Let $\overline{\NS(X)}$ be the primitive closure of $\NS(X)$ in $\NS(X_0)$.
Then the index of $\NS(X)$ in $\overline{\NS(X)}$ is a divisor of $\disc (\NS(X))$.
Therefore we obtain the following:
\begin{corollary}\label{cor:sp}
If $p$  does not divide $\disc(\NS(X))$, then 
the specialization isometry
$\rho: \NS(X)\isomet \NS(X_0)$ is primitive.
\end{corollary}
\begin{remark}
Artin~\cite[\S1]{MR0371899} showed a similar result over an  equal characteristic base.
Note  that the definition of supersingularity in~\cite[Definition (0.3)]{MR0371899} differs from ours.
\end{remark}
\subsection{Proof of Theorem~\ref{thm:genus}}\label{subsec:proofgenus}
Let $X\to\Spec F$ and $\XXX\to U$ be 
as in the statement of Theorem~\ref{thm:genus}.
Note that $\NS(X)$ is of signature $(1, 19)$, while 
the lattice $\cohom(X\sp{\sigma}, \Z)$ is 
even, unimodular and of signature $(3, 19)$
for any  $\sigma\in \Emb(F)$.
Hence 
$T(X\sp{\sigma})$ is even, positive-definite of rank $2$,
and its discriminant form is isomorphic to  $(\Disc{\NS(X)}, -\Discform{\NS(X)})$
by Proposition~\ref{prop:nikulin}.
Therefore $[T(X\sp{\sigma})]$ is contained in the genus 
$\GGG\subset \lats_{d(X)}$
characterized by 
$(\Disc{\GGG}, \Discform{\GGG})\cong (\Disc{\NS(X)}, -\Discform{\NS(X)})$.
\par
%
%
Let $\ppp$ be a point of $ \SSps_p(\XXX)$ with $p\notdiv 2 d(X)$.
Since  the Artin invariant of
$X\sbppp$ is $1$ by Proposition~\ref{prop:normalK3},
we have $\NS(X\sbppp)\cong \RSlat_{p, 1}$ by Theorem~\ref{thm:geomRSlat}.
Therefore  $L(\XXX,\ppp)$ is 
even, negative-definite of rank $2$.
On the other hand,  Corollary~\ref{cor:sp} implies that
the specialization isometry  $\rho$ is primitive, 
and hence the  discriminant form of $L(\XXX,\ppp)$ is isomorphic to
$\DRS{p,1}\oplus(\Disc{\NS(X)}, -\Discform{\NS(X)})$
by Proposition~\ref{prop:nikulin}.
It remains to show that
there exists $[M]\in \lats_{d(X)}$ such that $L(\XXX,\ppp)\cong M[-p]$,
or equivalently,
we  have $(x, y)\in p\Z$  for any $x, y\in L(\XXX,\ppp)$.
This follows from the following  lemma,
whose proof was given in~\cite{Milnor20}.
\begin{lemma}
Let $p$ be an odd prime integer, and 
$L$  an even  lattice of rank $2$.
If the $p$-part of $\Disc{L}$
is isomorphic to $(\Z/p\Z)\sp{\oplus 2}$,
then $(x, y)\in p\Z$ holds for any $x, y\in L$.
\end{lemma}
\section{$\Hom$-lattice}\label{sec:A}
\subsection{Preliminaries}
Let $E\sprime$ and $E$ be elliptic curves defined over a field $k$.
We denote by $\Hom_k (E\sprime, E)$
the $\Z$-module of homomorphisms from 
$E\sprime$ to $E$ defined over $k$,
and put 
\begin{eqnarray*}
&&\Hom(E\sprime, E)\;:=\;\Hom_{\bar k}(E\sprime\otimes \bar k, E\otimes \bar k), \\
&& \End_k (E):=\Hom_k (E, E)\quand \End (E):=\Hom (E, E)=\End_{\bar k}(E\otimes\bar{k}).
\end{eqnarray*}
The  Zariski  tangent space $T_O(E)$ of $E$ at the origin $O$
is a one-dimensional $k$-vector space,
and hence $\End_k (T_O(E))$ is canonically isomorphic to $k$.
By  the action of $\End_k (E)$ on $T_O(E)$, we have a representation
$$
\map{\Lie}{\End_k(E)}{\End_k (T_O(E))=k}.
$$
According to~\cite[\S6 in Chap.~III]{MR817210},
we define a lattice structure  on $\Hom (E\sprime, E)$
by
\begin{equation*}\label{eq:fg}
(f, g)\;:=\;\deg (f+g)-\deg(f)-\deg (g).
\end{equation*}
We consider the product abelian surface
$$
A\;:=\; E\sprime\times E.
$$
Let $O\sprime\in E\sprime$  and $O\in E$  be the origins.
We put
\begin{equation}\label{eq:xieta}
\xi:=[E\sprime \times \{O\}]\;\in \;\NS(A),
\qquad 
\eta:=[\{O\sprime\}\times E]\;\in\; \NS(A),
\end{equation}
and denote by $U(A)$ the sublattice of $\NS(A)$ spanned by $\xi$ and $\eta$,
which is even, unimodular and of signature $(1, 1)$.
The following is classical.
See~\cite{MR0029522}, for example.
\begin{proposition}\label{prop:NSA}
The lattice $\NS(A)$ is isomorphic to
$U(A)\perp \Hom (E\sprime, E)[-1]$.
In particular,
the lattice $\Hom (E\sprime, E)$ is  even and positive-definite.
\end{proposition}
One can easily prove the following propositions
by means of, for example, 
the results in~\cite[\S9 in Chap.~III]{MR817210} and~\cite[\S3]{MR0206004}.
\begin{proposition}\label{prop:Hom_char0}
Suppose that $\ch k=0$.
Then the following  are equivalent:
\begin{rmenumerate}
\item $\rank( \Hom (E\sprime, E))=2$.
\item $E\sprime$ and $E$ are isogenous over $\bar k$, and $\rank (\End (E\sprime))=2$. 
\item There exists an imaginary quadratic field $K$
such that both of $\End (E\sprime)\otimes\Q$ and $\End (E)\otimes\Q$ are isomorphic to $K$.
\end{rmenumerate}
\end{proposition}
\begin{proposition}\label{prop:Hom_charp}
Suppose that $\ch k>0$.
Then the following are equivalent:
\begin{rmenumerate}
\item $\rank (\Hom (E\sprime, E))=4$.
\item $E\sprime$ and $E$ are isogenous over $\bar k$, and $\rank (\End (E\sprime))=4$. 
\item Both of $E\sprime$ and $E$ are supersingular.
\end{rmenumerate}
\end{proposition}
\subsection{The elliptic curve $E^J$}\label{subsec:EJ}
To the end of \S\ref{subsec:Hom_p}, 
we work over an algebraically closed  field $k$.
For an elliptic curve $E$, 
we denote by $k(E)$ the function field of $E$.
\begin{definition}\label{def:isomisog}
Two non-zero isogenies $\phi_1: E\to E_1$
and $\phi_2: E\to E_2$ are \emph{isomorphic} if 
there exists an isomorphism $\psi: E_1\isom E_2$ such that
$\psi\circ\phi_1=\phi_2$ holds,
or equivalently, if the subfields $\phi_1\sp * k(E_1)$ and $\phi_2\sp * k(E_2)$ of $k(E)$ are equal.
\end{definition}
For a non-zero endomorphism $a\in \End(E)$, we denote by $E^a$ the image of $a$, 
that is, $E^a$ is an elliptic curve isomorphic to $E$
with an isogeny $a:E\to E^a$.
The function field $k(E^a)$ is canonically identified with the  subfield 
$a^* k(E)=\shortset{a^* f}{f\in k(E)}$
of $k(E)$, and we have
$[k(E):k(E^a)]=\deg a$.
\begin{definition}\label{def:EJ}
Let $J\subset \End (E)$ be a non-zero left-ideal of $\End (E)$.
We denote by $k(E^J)\subset k(E)$ the composite of the subfields
$k(E^a)$ for all non-zero $a\in J$.
Then $k(E^J)$ is a function field of an elliptic curve $E^J$.
We denote by 
$$
\map{\phi^J}{ E}{ E^J}
$$
the isogeny corresponding to $k(E^J)\inj k(E)$.
\end{definition}
\begin{remark}\label{rem:ab}
Let  $a, b\in \End (E)$ be non-zero.
Since $ba(x)=b(a(x))$, we have canonical inclusions 
$k(E^{ba})\subset k(E^a)\subset k(E)$.
Hence,
if the left-ideal $J$ is generated by non-zero elements $a_1, \dots, a_t$,
then $k(E^J)$ is the composite of 
$k(E^{a_1}), \dots, k(E^{a_t})$.
\end{remark}
\begin{remark}\label{rem:characterization}
The isogeny $\phi^J : E\to E^J$ is 
characterized by the following properties:
(i) every $a\in J$ factors through $\phi^J$, and
(ii) if  every $a\in J$ factors through an isogeny $\psi: E\to E\sprime$,
then $\phi^J$ factors through $\psi$.
\end{remark}
\subsection{The $\Hom$-lattice in characteristic $0$}\label{subsec:Hom_0}
In this subsection,
we assume that $k=\bar k$ is of characteristic $0$,
and that the conditions in Proposition~\ref{prop:Hom_char0}
are satisfied.
We denote by $D$ the discriminant of the imaginary quadratic field  $K$ in the condition (iii)
of Proposition~\ref{prop:Hom_char0}.
Note that $\End (E)$ is isomorphic to a $\Z$-subalgebra of $\ZK$
with $\Z$-rank $2$,
and that there exist two embeddings 
of $\End (E)$ into  $\ZK$ as a $\Z$-subalgebra 
that are conjugate over $\Q$.
Each embedding $\End (E)\inj \Z_K$ is  an isometry of lattices,
where $\Z_K$ is considered as a lattice by the formula~\eqref{eq:bilK},
because the dual endomorphism corresponds to the conjugate element  over $\Q$.
\begin{proposition}\label{prop:mnsq}
There exist non-zero integers $m$ and $n$ such that
\begin{equation}\label{eq:mnsq}
m^2\disc(\Hom(E\sprime, E))=-n^2 D.
\end{equation}
\end{proposition}
\begin{proof}
There exists a non-zero isogeny $\alpha: E\to E\sprime$.
Then the map $g\mapsto g\circ\alpha$ induces an isometry $\Phi_\alpha$
from $\Hom (E\sprime, E)[\deg \alpha]$ to $\End(E)$.
Putting $m:=\deg \alpha$ and 
$n:=[\Z_K: \End(E)]\cdot | \Coker \Phi_\alpha|$, 
we obtain the equality~\eqref{eq:mnsq}.
\end{proof}
\begin{definition}
Since $\bar k=k$, we have $\Lie :\End(E)\to k$.
Suppose that an embedding $i:K\inj k$ is fixed.
Then 
an embedding $\iota: \End (E)\inj \ZK$ as a $\Z$-subalgebra 
is called   \emph{$\Lie$-normalized} if 
$\Lie :\End(E)\to k$
coincides with the composite of
$\iota: \End (E)\inj \ZK$,
the inclusion 
$\ZK\inj K$ and $i: K\inj k$.
\end{definition}
\begin{definition}\label{def:Jstar}
Suppose that $k=\C$,  and that $\End (E)\cong \Z_K$.
We fix an embedding $K\inj \C$.
Let $\Lambda\subset \C$ be a $\Z$-submodule of rank $2$ such that
$E\cong \C/\Lambda$ as  a Riemann surface.
For an ideal class $[I]$  of $\Z_K$
represented by a  fractional ideal $I\subset K\subset\C$,
we denote by $[I]* E$ the complex elliptic curve
$\C/I\inv\Lambda$,
where $I\inv\Lambda$ is the $\Z$-submodule of $\C$
generated by $x\lambda$ $(x\in I\inv, \lambda\in \Lambda)$.
When $I\subseteq \Z_K$, we have 
$I\inv \Lambda\supset\Lambda$, and
the identity map $\id_{\C}$  of $\C$ induces an isogeny
$$
\map{{}^{\an}\phi^I\;\;}{\;\;E=\C/\Lambda\;}{\;[I]*E=\C/I\inv \Lambda}.
$$
\end{definition}
\begin{proposition}\label{prop:anaEJ}
Suppose that $k=\C$, and that $\End(E)\cong \ZK$.
For   an ideal $J\subset \End (E)$,
the isogeny  $\phi^J: E\to E^J$ is isomorphic to  ${}^{\an}\phi^{J}: E\to [J]*E$,
where $J$ is regarded as an ideal of $\ZK$ by the $\Lie$-normalized isomorphism
$\End(E)\cong \ZK$. 
\end{proposition}
\begin{proof}
Suppose that $E=\C/\Lambda$.
We choose $\Lambda\sprime\subset \C$ such that
$E^J=\C/\Lambda\sprime$, and that  
$\phi^J: E=\C/\Lambda\to E^J=\C/\Lambda\sprime$  
is given by  $\id_{\C}$.
For a non-zero $a\in J$,
we have $(1/a)\Lambda\supset\Lambda$ and 
there exists 
a canonical isomorphism $E^a=\C/(1/a)\Lambda$
such that $a: E\to E^a$ is given by  $\id_{\C}$.
Therefore $\Lambda\sprime$ is the largest $\Z$-submodule of $\C$
that is contained in  $(1/a)\Lambda$
for any non-zero $a\in J$.
Hence we have $\Lambda\sprime=J\inv \Lambda$.
\end{proof}
From this analytic description of $\phi^J: E\to E^J$, 
we obtain the following,
which holds in  any field of characteristic $0$.
\begin{proposition}\label{prop:PhiJchar0}
Suppose that $\ch k=0$, and that 
$\End (E)\cong \Z_K$.
Let $J$ be an ideal of $\End(E)$.
Then 
$\End (E^J)$ is also isomorphic to $\Z_K$.
Moreover, 
$\deg\phi^J$ is equal to $|\End(E)/J|$,
and the image of the map
$$
\map{\Phi^J}{ \Hom (E^J, E)}{ \End(E)}
$$ 
given by $g\mapsto g\circ \phi^J$ coincides with $J$.
\end{proposition}
\subsection{The $\Hom$-lattice  of supersingular elliptic curves}\label{subsec:Hom_p}
In this subsection,
we assume that $k=\bar k$ is of characteristic $p>0$,
and that the conditions in Proposition~\ref{prop:Hom_charp}
are satisfied.
In particular, $E$ is a supersingular elliptic curve.
\par
\smallskip
We denote by $\quaternion$ the quaternion algebra over $\Q$ that ramifies exactly 
at $p$ and $\infty$.
It is well-known that $\quaternion$ is unique up to isomorphism.
We denote by $x\mapsto x\sp*$  the canonical involution of 
$\quaternion$.
Then $\quaternion$ is equipped with a positive-definite 
$\Q$-valued symmetric bilinear form
defined by
\begin{equation}\label{eq:bilB}
(x, y)\;:=\;xy^*+yx^*.
\end{equation}
A subalgebra of $\quaternion$ is called an \emph{order} if its $\Z$-rank is $4$.
An order is said to be \emph{maximal} if it is maximal 
among  orders with  respect to the inclusion.
If $R$ is an order of $\quaternion$,
then the bilinear form~\eqref{eq:bilB} takes values in $\Z$ on $R$,
and $R$ becomes an even lattice.
It is known that  $R$ is maximal if and only if the discriminant of $R$ is $p^2$.
The following are the classical results due to Deuring~\cite{MR0005125}.
(See also~\cite[Chapter~13, Theorem~9]{MR890960}):
\begin{proposition}\label{prop:Endss1}
There exists a maximal order $R$ of $\quaternion$
such that 
$\End(E)$ is isomorphic to $ R$ as a $\Z$-algebra.
The canonical involution of $R$ corresponds to
the involution $\phi\mapsto\phi^*$ of $\End(E)$,
where $\phi^*$ is the dual endomorphism. Hence
the  lattice $\End(E)$ is isomorphic to the lattice $R$,
and we have  $\disc (\End(E))=p^2$.
\end{proposition}
Conversely, we have the following:
\begin{proposition}\label{prop:Endss2}
Let $R$ be a maximal order of $B$.
Then there exists a supersingular elliptic curve $E_R$
such that
$\End(E_R)$ is isomorphic to $R$  as a  $\Z$-algebra.
\end{proposition}
We fix an isomorphism $\End(E)\otimes\Q\cong \quaternion$
such that $\End (E)$ is mapped to a maximal order $R$ of $\quaternion$.
Let $J$ be a non-zero left-ideal of $\End (E)$.
Consider the left- and right-orders 
$$
O_l(J):=\shortset{x\in \quaternion}{xJ\subset J},
\quad
O_r(J):=\shortset{x\in \quaternion}{Jx\subset J}
$$
of $J$.
Since $O_l(J)$ contains  $R$ and $R$ is maximal,
$O_l(J)$ is maximal,
and hence $O_r(J)$ is also maximal by~\cite[Theorem (21.2)]{MR1972204}.
In other words, $J$ is a \emph{normal ideal} of $\quaternion$.
We denote by $\nr (J)$ the greatest common divisor of the integers
$$
\nr (\phi)\;:=\;\phi\phi^*=\deg \phi \qquad(\phi\in J).
$$
(See~\cite[Corollary (24.12)]{MR1972204}.)
Then, by~\cite[Theorem (24.11)]{MR1972204},  we have
\begin{equation}\label{eq:nrsq}
\nr (J)^2=|R/J|.
\end{equation}
On the other hand, Deuring~\cite[(2.3)]{MR0005125} proved the following:
\begin{equation}\label{eq:degDeuring}
\deg \phi^J=\nr (J).
\end{equation}
\begin{proposition}\label{prop:PhiJcharp}
The image of the map $\Phi^J: \Hom (E^J, E)\to \End (E)$
given by $g\mapsto g\circ \phi^J$ is equal to $J$.
\end{proposition}
\begin{proof}
By Remark~\ref{rem:characterization},
we have $J\subseteq \Im \Phi^J$.
Suppose that there exists $a \in \Im \Phi^J$ such that 
$a\notin J$.
Let $J\sprime$ be the left-ideal of $\End (E)$ 
generated by $J$ and $a$.
Then we have $\nr(J\sprime)<\nr(J)$ by the formula~\eqref{eq:nrsq}.
On the other hand,
since $a$ factors through $\phi^J$,
we have $k(E\sp a)\subset k(E^J)$ and hence $k(E^{J\sprime})= k(E^J)$.
This contradicts Deuring's formula~\eqref{eq:degDeuring}.
\end{proof}
%
%
%
%
\begin{proposition}\label{prop:ImPsi}
Let $\psi: E\to E\spprime$ be a non-zero isogeny, and let
$$
\map{\Psi}{ \Hom (E\spprime, E)}{ \End (E)}
$$ 
be the homomorphism of $\Z$-modules given by $g\mapsto g\circ \psi$.
We denote  by $J_\psi$ the image of $\Psi$,
which is a left-ideal of $\End (E)$.
Then $\psi$ is equal to 
$\phi^{J_\psi}$. 
\end{proposition}
\begin{proof}
Since $k(E^{g\circ \psi})\subset k(E\spprime)$ as subfields of $k(E)$
for any non-zero $g\in \Hom(E\spprime, E)$, we have
$k(E^{J_\psi})\subset k(E\spprime)$, and hence 
$\phi^{J_\psi}: E\to E^{J_\psi}$ factors through
$\psi: E\to E\spprime$.
The greatest common divisor of the degrees 
of $g \in \Hom (E\spprime, E) $ is $1$
by Proposition~\ref{prop:gcddeg}
in the next subsection.
Hence we have $\nr (J_\psi)=\deg \psi$
by the definition of $\nr$.
Since $\deg \phi^{J_\psi}=\nr (J_\psi)$ by Deuring's formula~\eqref{eq:degDeuring},
we have $\psi=\phi^{J_\psi}$ and  $E\spprime= E^{J_\psi}$.
\end{proof}
\begin{corollary}
The map $J\mapsto \phi^J$
establishes 
a one-to-one correspondence between
the set of non-zero left-ideals of $\End (E)$ and the set of
isomorphism classes of  non-zero isogenies from $E$.
\end{corollary}
\begin{proposition}\label{prop:psq}
Let $E\sprime$ and $E$ be supersingular. 
Then the discriminant of the lattice $\Hom (E\sprime, E)$ is equal to $p^2$.
\end{proposition}
\begin{proof}
Since $E\sprime$ and $ E$ are isogenous,
there exists a  non-zero left-ideal $J$ of $\End (E)$ such that $E\sprime\cong E^J$.
Then we have an isomorphism
$\Hom (E\sprime, E) \cong J$
of $\Z$-modules given by $g\mapsto g\circ \phi^J$,
and hence we have $\Hom (E\sprime, E)[\deg \phi^J]\cong J$
as a lattice, 
from which we obtain
$$
\disc (\Hom (E\sprime, E))=\frac{\disc (J)}{(\deg \phi^J)^4}
=\frac{\disc (\End (E)) \cdot [\End(E):J]^2}{(\deg \phi^J)^4}=\disc (\End (E))
$$
by the formulae~\eqref{eq:nrsq} and~\eqref{eq:degDeuring}.
Thus we have $\disc (\Hom (E\sprime, E))=p^2$ by Proposition~\ref{prop:Endss1}.
\end{proof}
\subsection{The specialization isometry of $\Hom$-lattices}\label{subsec:spHom}
Let $E$ be an elliptic curve defined over  a finite extension  $L\subset \cloQ _p$ of $\Q_p$ such that
the $j$-invariant $j(E)\in L$ is integral over $\Z_p$.
This condition is satisfied, for example, 
if 
$\rank (\End(E))=2$.
Then $E$ has potentially good reduction, that is,
there exist a finite extension $M\subset \cloQ_p$ of $L$ and 
a smooth proper morphism  $\EEE_M  \to \Spec\Z_M$
over the valuation ring $\Z_M$ of  $M$ such that 
$\EEE_M\otimes  M$ is isomorphic to $E\otimes M$.
Let $E_0$ be the central fiber of $\EEE_M$.
Then we have a specialization isometry
$$
\rho\;:\; \End (E)\;\isomet \; \End (E_0),
$$
which is obtained from  the specialization isometry $\NS(E\times E)\isomet \NS(E_0\times E_0)$ 
and Proposition~\ref{prop:NSA}.
The following  follows, for example,
from the existence and the uniqueness of the N\'eron model~\cite[Chap.~IV]{MR1312368}.
\begin{proposition}\label{prop:indep}
The isomorphism class of $E_0$  over $\overline{\F}_p$ and the specialization isometry $\rho$
do not depend on the choice of $M$ and $\EEE_M$.
\end{proposition}
Replacing $L$ by a finite extension  if necessary, we assume that
$$
\End (E)\;=\; \End_L(E), 
$$
so that 
$\shortmap{\Lie }{ \End (E)}{ L}$
is defined.
Let $E\sprime$ be another elliptic curve defined over  a finite extension  $L\sprime\subset \cloQ _p$ of $\Q_p$
such that $j(E\sprime)\in L\sprime$ is integral over $\Z_p$.
Then we have a specialization isometry
$\shortmap{\rho\sprime}{ \End (E\sprime)}{\End (E\sprime_0)}$,
where $E\sprime_0$ is the central fiber of a N\'eron model of $E\sprime$.
Replacing $L\sprime$ by a finite extension,
we assume that
$\End (E\sprime)=\End_{L\sprime}(E\sprime)$.
The following is easy to prove.
\begin{proposition}\label{prop:indep2}
Suppose that there exists $g\in \Gal (\cloQ _p/\Q_p)$ such that
$j(E\sprime)=j(E)\sp g$.
Then there exist isomorphisms
$\End (E)\cong\End (E\sprime)$ and $\End (E_0)\cong\End (E_0\sprime)$ 
induced from  $g$ such that the following diagram is commutative:
$$
\renewcommand{\arraystretch}{1.4}
\begin{array}{ccccccc}
\cloQ _p & \leftinj & L&\mapleftsp{\Lie} &\End(E) & \maprightsp{\rho} & \End(E_0) \\
\smallmath{g} \big\downarrow\wr &  & &&\big\downarrow\wr &                 &\big\downarrow\wr \\                       
\cloQ _p & \leftinj & L\sprime &\mapleftsp{\Lie} &\End(E\sprime) & \maprightsp{\rho\sprime} & \End(E\sprime_0).
\end{array}
$$
\end{proposition}
Suppose that $\End(E)\otimes \Q$ is isomorphic to
an imaginary quadratic field $K$.
The following result is again due to Deuring~\cite{MR0005125}.
(See also~\cite[Chapter~13, Theorem~12]{MR890960}).
\begin{proposition}\label{prop:redss}
The elliptic curve $E_0$ is supersingular if and only if 
$p$ is inert or ramifies in $K$.
\end{proposition}
We now work over $\cloQ_p$, and 
assume that $\End (E)$ is isomorphic to $\Z_K$.
Suppose that $E_0$ is supersingular.
We put
$R:=\End(E_0)$.
Let $J$ be an ideal of $\End(E)$, and 
consider the elliptic curve $E^J$.
Since $\End (E^J)$ is also isomorphic to $\Z_K$ by Proposition~\ref{prop:PhiJchar0},
the reduction $(E^J)_0$ of $E^J$ is supersingular by Proposition~\ref{prop:redss},
and we have a reduction
$$
\map{\rho(\phi^J)}{ E_0}{(E^J)_0}
$$
of the isogeny $\phi^J: E\to E^J$.
On the other hand, 
we have the left-ideal $R\cdot \rho(J)$
of $R$ generated by $\rho(J)\subset R$,
and  the associated  isogeny
$$
\map{\phi^{RJ}}{ E_0}{(E_0)^{R \rho(J)}}.
$$
\begin{proposition}\label{prop:isogs}
The isogenies $\rho(\phi^J)$ and $\phi^{RJ}$ are isomorphic.
\end{proposition}
\begin{proof}
We choose  $a_1, \dots, a_t\in J$ such that
 $J$ is generated by $a_1, \dots, a_t$, and that
 $[\End(E):J]$ is  equal to the greatest common divisor  of $\deg a_1, \dots, \deg a_t$.
By Proposition~\ref{prop:PhiJchar0},
we have  $\deg \rho(\phi^J)=\deg \phi^J=[\End(E):J]$.
By Deuring's formula~\eqref{eq:degDeuring}, we see that $\deg \phi^{RJ}$ is a common divisor of
$\deg \rho(a_i)=\deg a_i$ for $i=1, \dots, t$,
and hence $\deg \phi^{RJ}$ divides $\deg \rho(\phi^J)$.
On the other hand,
the left-ideal $R\cdot \rho(J)$ is generated by 
$\rho(a_1), \dots, \rho(a_t)$, and hence, 
by Remarks~\ref{rem:ab} and~\ref{rem:characterization},
we see that $\phi^{RJ}$ factors through $ \rho(\phi^J)$.
Therefore we obtain $\rho(\phi^J)=\phi^{RJ}$.
\end{proof}
By Proposition~\ref{prop:isogs}, the following diagram is   commutative:
$$
\renewcommand{\arraystretch}{1.2}
\begin{array}{ccc}
\Hom(E^J, E) & \isomet & \Hom((E^J)_0, E_0)\\
\llap{\smallmath{\Phi^J}}\downinj & & \downinj\rlap{\;\;\smallmath{\Phi^{RJ}}}  \\
\End(E) & \isomet &\End(E_0),
\end{array}
$$
where 
the horizontal arrows are the specialization isometries.
By  Propositions~\ref{prop:PhiJchar0} and~\ref{prop:PhiJcharp}, we obtain the following:
%
%
\begin{proposition}\label{prop:JRJsperp}
We put $d_J:=\deg \phi^J=\deg \rho(\phi^J)=\deg \phi^{RJ}$.
Then we have an isomorphism of lattices
$$
(\Hom(E^J, E)  \isomet  \Hom((E^J)_0, E_0))\sperp [d_J]\;\wcong\; (J\isomet R\cdot\rho(J))\sperp,
$$
where, in the right-hand side,
$J$ and $R\cdot\rho(J)$ are regarded as sublattices
of the lattices $\End(E)\cong\ZK$ and $\End(E_0)=R$,
respectively,
and $J\isomet R\cdot\rho(J)$ is given by the specialization isometry $\rho:\End(E)\isomet R$.
\end{proposition}
Finally, we state 
the lifting theorem of Deuring~\cite{MR0005125}.
See also~\cite[Chapter~13, Theorem~14]{MR890960} and~\cite[Proposition 2.7]{MR772491}.
\begin{proposition}\label{prop:lift}
Let $E_0$ be a supersingular elliptic curve defined over 
a field $\kappa_0$ of characteristic $p$,
and $\alpha_0$ an endomorphism of $E_0$.
Then there exist a smooth proper family of elliptic curves $\EEE\to \Spec \Z_L$
over the valuation ring $\Z_L$ of a finite extension $L$ of $\Q_p$ and 
an endomorphism $\alpha$ of $\EEE$  over $\Z_L$
such that $(\EEE,\alpha)\otimes \bar{\kappa}\sbppp$ is
isomorphic to $(E_0, \alpha_0)\otimes \bar{\kappa}_0$,
where $\ppp$ is the closed point of $\Spec \Z_L$.
\end{proposition}
\subsection{Application of Tate's theorem~\cite{MR0206004}}\label{subsec:Tate}
In this subsection,
we prove the following result, which was used in the proof of Proposition~\ref{prop:ImPsi}.
\begin{proposition}\label{prop:gcddeg}
Let $E\sprime$ and $E$ be supersingular elliptic curves.
Then 
the greatest common divisor of the degrees 
of $g \in \Hom (E\sprime, E) $ is $1$.
\end{proposition}
\begin{proof}
Without loss of generality,
we can assume that $E\sprime$ and $E$ are defined over 
a finite field $\F_q$
of characteristic $p$.
Replacing $\F_q$ by a finite extension,
we can assume that 
$\End(E\sprime)=\End_{\F_q}(E\sprime)$,
$\End(E)=\End_{\F_q}(E)$ and 
$\Hom (E\sprime, E)=\Hom _{\F_q}(E\sprime, E)$
hold.
Let $l$ be a prime integer $\ne p$,
and consider the $l$-adic Tate module $T_l(E\sprime)$ of $E\sprime$.
By the famous theorem of Tate~\cite{MR0206004},
we see that
$$
\End_{\Gal (\overline{\F}_q/\F_q)} (T_l(E\sprime))
\wcong 
\End_{\F_q}(E\sprime)\otimes {\Z_l}
\,=\,
\End(E\sprime)\otimes {\Z_l}
$$ 
is of rank $4$,
and hence we can assume that
the $q$-th power Frobenius morphism $\Frob_{E\sprime}$  acts on $T_l(E\sprime)$
as a scalar multiplication by $\sqrt{q}$.
In the same way,
we can  assume that
$\Frob_E$ acts on $T_l(E)$
as a scalar multiplication by $\sqrt{q}$.
Then, by the theorem of Tate~\cite{MR0206004} again, we have a natural  isomorphism 
$$
\Hom (E\sprime, E) \otimes \Z_l \wcong\Hom(T_l(E\sprime), T_l(E)) \wcong \End_{\Z_l}(\Z_l\sp{\oplus 2}).
$$
Hence there exists $g\in \Hom (E\sprime, E)$
such that $\deg g$ is not divisible by $l$.
Therefore the greatest common divisor of the degrees 
of $g \in \Hom (E\sprime, E) $ is a power of $p$.
Let $\shortmap{F}{ E\sprime}{ E\sp{\prime (p)}}$ 
be the $p$-th power Frobenius morphism of $E\sprime$.
If the degree of $g: E\sprime\to E$ is divisible by $p$, then $g$ factors 
as $g\sprime\circ F$ with $\deg g\sprime=\deg g/p$.
Therefore it is enough to show the following:
\par\smallskip\noindent
{\bf Claim.} 
For any supersingular elliptic curve $E$ in characteristic $p$,
there exists $g\in \Hom (E, E\sp{(p)})$
such that $\deg g$ is prime to $p$.
\par\smallskip
Note that $j(E)\in \F_{p^2}$
and $j(E^{(p)})=j(E)^p$.
By~Proposition~\ref{prop:lift},
there exists an elliptic curve $E\sp{\sharp}$ defined over a finite extension  $L$ of $\Q_p$
such that $\End (E\sp{\sharp})$ is of rank $2$,
and that $E\sp{\sharp}$ has a  reduction 
isomorphic to $E$
at the  closed point  $\ppp$ of $\Z_L$.
We  assume that $L$ is Galois over $\Q_p$,
and fix an embedding $L\inj \C$.
Then $\End (E\sp{\sharp})$ is an order $\OOO$ of an imaginary quadratic field,
and $E\sp{\sharp}\otimes\C$ is isomorphic to $\C/I$
as a Riemann surface
for some invertible $\OOO$-ideal $I$
(\cite[Corollary 10.20]{MR1028322}).
Note that $j(E\sp{\sharp})$ is a root of the Hilbert class polynomial of the order $\OOO$
(\cite[Proposition 13.2]{MR1028322}).
There exists an element $\gamma\in \Gal(L/\Q_p)$
such that 
$$
j(E\sp{\sharp})\sp{\gamma}\equiv j(E\sp{\sharp})\sp p\;\bmod\,\ppp.
$$
We put  $E\sp{\flat}:=(E\sp{\sharp})\sp{\gamma}$.
Then $E\sp{\flat}$ has a  reduction 
isomorphic to $E\sp{(p)}$
at  $\ppp$, and 
we have $E\sp{\flat}\otimes\C\cong \C/J$
as a Riemann surface
for some invertible $\OOO$-ideal $J$.
The degree of homomorphisms
in $\Hom (E\sp{\sharp}, E\sp{\flat})=\Hom(\C/I, \C/J)$ is 
given by a primitive binary form 
corresponding to the ideal class of the 
proper $\OOO$-ideal $I\inv J$
by~\cite[Theorem 7.7]{MR1028322}. By~\cite[Lemma 2.25]{MR1028322}, 
we see that 
$\Hom (E\sp{\sharp}, E\sp{\flat})$
has an element whose degree is prime to $p$.
Since the specialization homomorphism
$\Hom (E\sp{\sharp}, E\sp{\flat})\to \Hom (E, E\sp{(p)})$
preserves the degree, we obtain the proof.
\end{proof}
\section{Kummer construction}\label{sec:Kummer}
We denote by $k$  an algebraically closed field of characteristic $\ne 2$.
\subsection{Double coverings}\label{subsec:double}
We work over   $k$.
Let $W$ and $Z$ be smooth projective surfaces,
and $\phi:W\to Z$ a finite double covering.
Let $\iota: W\isom W$ be the deck-transformation of $W$ over $Z$.
Then we have homomorphisms
$$
\phi_*:\NS(W)\to\NS(Z)\quand \phi^*:\NS(Z)\to\NS(W).
$$
Let $\NS(W)_{\Q}^+\subset \NS(W)\otimes \Q$ be 
the eigenspace of $\iota_*$ with the eigenvalue $1$.
We put
$$
\NS(W)^+\;:=\;\NS(W)\cap \NS(W)_{\Q}^+.
$$
When the base field $k$ is $\C$,
we assume that $\cohom (W, \Z)$ and $\cohom(Z, \Z)$ are torsion-free,
so that they can be regarded as lattices.
We have homomorphisms
$$
\phi_*:\cohom(W, \Z)\to \cohom(Z, \Z)\quand \phi^*:\cohom(Z, \Z)\to \cohom(W,\Z).
$$
Note that  $\phi^*$ preserves the Hodge structure.
We define 
$\cohom (W, \Z)^+:=\cohom(W,\Z)\cap \cohom(W,\Q)^+$
in the same way as $\NS(W)^+$.
\begin{lemma}\label{lem:double}
The homomorphism $\phi_*$ induces an isometry
$$
\mapisomet{\phi_*^+}{\NS(W)^+[2]}{\NS(Z)}
$$
with a  finite $2$-elementary cokernel.
When $k=\C$, 
$\phi_*$ induces an isometry
$$
\mapisomet{\phi_*^+}{\cohom(W,\Z)^+[2]}{ \cohom(Z,\Z)}
$$
with a  finite $2$-elementary cokernel
that preserves the Hodge structure.
\end{lemma}
\begin{proof}
The proof follows immediately from the following:
\begin{equation*}\label{eq:double}
\renewcommand{\arraystretch}{1.4}
\begin{array}{ccl}
&&\phi^*\circ \phi_*(w) =w+\iota_*(w),\quad
\phi_*\circ \phi^*(z) =2z,\quad
\iota_*\circ \phi^* (z)=\phi^*(z),\\
&&(\phi^*(z_1), \phi^*(z_2)) = 2 (z_1, z_2), \quad
(\iota_*(w_1), \iota_*(w_2)) =  (w_1,w_2).
\end{array}
\end{equation*}
The inverse of the isomorphism $\phi_*^+\otimes\Q$ is given by $(1/2) \phi^*\otimes\Q$.
\end{proof}
\subsection{Disjoint $(-2)$-curves}\label{subsec:disjoint}
We continue to work over $k$.
Let $C_1$, \dots, $C_m$ be $(-2)$-curves on a $K3$ surface $X$  that 
are disjoint to  each other,
$\Delta\subset \NS(X)$ the sublattice
generated by $[C_1]$, \dots, $[C_m]$,
and $\overline{\Delta}\subset \NS(X)$ the primitive closure of $\Delta$.
The discriminant group $\Disc{\Delta}$ of $\Delta$ is isomorphic to $\F_2^{\oplus m}$
with basis
$$
\gamma_i\;:=\;-[C_i]/2+\Delta\quad (i=1, \dots, m).
$$
For $x=x_1\gamma_1+\dots+x_m\gamma_m\in \Disc{\Delta}$,
we denote by $\weight (x)$ the \emph{Hamming weight} of $x$, 
that is, the number of $x_i\in \F_2$ with $x_i\ne 0$.
Then $\Discform{\Delta}:\Disc{\Delta}\to\Q/2\Z$ is 
given by
$$
\Discform{\Delta}(x)=(-\weight(x)/2) +2\Z\;\;\in\;\;\Q/2\Z.
$$
\begin{lemma}\label{lem:nodal}
We put $H_\Delta:=\overline{\Delta}/\Delta\subset \Disc{\Delta}$. Then, 
for every $x\in H_{\Delta}$, we have $\weight(x)\equiv 0\bmod 4$
and $\weight(x)\ne 4$.
\end{lemma}
\begin{proof}
Since $H_{\Delta}$ is totally isotropic
with respect to  $\Discform{\Delta}$,  we have $\weight(x)\equiv 0\bmod 4$ for any $x\in H_{\Delta}$.
Let $\gamma: X\to Y$ be the contraction of  $C_1$, \dots, $C_m$,
and  $\LLL_Y$  a very ample line bundle on the normal $K3$ surface $Y$.
Then $\{[C_1], \dots, [C_m]\}$ is a fundamental system of roots
in the the root system
\begin{equation}\label{eq:rootsystem}
\set{r\in \NS(X)}{(r, [\gamma\sp *\LLL_Y])=0, \;\; (r,r)=-2}
\end{equation}
of type $m\A_1$.
(See~\cite[Proposition 2.4]{normalK3}.)
If there were $x\in H_{\Delta}$ with $\weight (x)=4$,
then there would exist a vector $r$ in the set~\eqref{eq:rootsystem}
such that $r\ne\pm  [C_i]$ for any $i$, 
which is a contradiction.
\end{proof}
\subsection{Double Kummer pencil}
Let $E\sprime$ and $E$ be elliptic curves defined over $k$.
We put $A:=E\sprime\times E$,
and denote by $\Km (A)$ the Kummer surface associated with
$A$,
that is, $\Km (A)$ is the minimal resolution of the quotient surface
$A/\ang{\iota_A}$,
where $\iota_A: A\isom  A$ is the inversion automorphism $x\mapsto -x$.
Let $u\sprime_i $ and $u_j$ $ (1\le i, j\le 4)$ be the points of order $\le 2$
in $E\sprime$ and $E$, respectively, and 
let $\beta_A: \Bl{A}\to A$  be the blowing-up of $A$
at the fixed points $(u\sprime_i, u_j)$ of $\iota_A$.
Let  $\phiA: \Bl{A}\to \Km(A)$ denote the natural finite double covering.
The involution $\iota_A$ lifts to an involution $\tiliota_A$
of $\Bl{A}$,
and $\phiA$ is the quotient morphism  $\Bl{A}\to \Bl{A}/\ang{\tiliota_A}=\Km(A)$.
\begin{definition}
The diagram
$$
\Km(A)\mapleftsp{\phiA} \Bl{A}\maprightsp{\beta_A} A=E\sprime\times E
$$
is called the \emph{Kummer diagram} of $E\sprime$ and $E$.
We denote by $E_{ij}\subset \Km(A)$ the image by $\phiA$ of the exceptional curve 
of $\beta_A$ over the point $(u\sprime_i, u_j)\in A$,
and by $F_j$ and $G_i$ the image by $\phiA$
of the strict transforms of $E\sprime\times \{u_j\}$ and $\{u\sprime_i\}\times E$,
respectively.
These $(-2)$-curves $E_{ij}$, $F_j$ and  $G_i$ on $\Km(A)$ form 
the configuration depicted  in Figure~\ref{fig:Km},
which is called  the \emph{double Kummer pencil} (see~\cite{MR0441982}).
\end{definition}
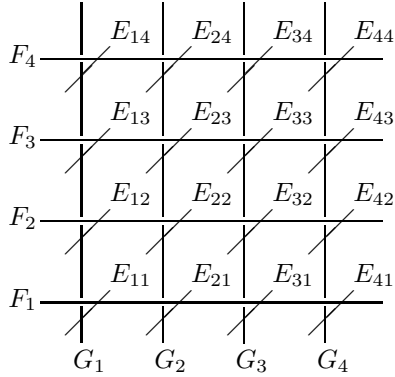
\begin{figure}
\begin{center}
\setlength{\unitlength}{1.08mm}
\begin{picture}(60, 47)
%
%
\put(5, 10){\line(1, 0){42}}
\put(1,9.5){$F_1$}
\put(5, 20){\line(1, 0){42}}
\put(1,19.5){$F_2$}
\put(5, 30){\line(1, 0){42}}
\put(1,29.5){$F_3$}
\put(5, 40){\line(1, 0){42}}
\put(1,39.5){$F_4$}
%
%
\put(10, 5.1){\line(0, 1){4.4}}
\put(10, 10.6){\line(0, 1){8.9}}
\put(10, 20.6){\line(0, 1){8.9}}
\put(10, 30.6){\line(0, 1){8.9}}
\put(10, 40.6){\line(0, 1){6.4}}
\put(9,2){$G_1$}
\put(20, 5.1){\line(0, 1){4.4}}
\put(20, 10.6){\line(0, 1){8.9}}
\put(20, 20.6){\line(0, 1){8.9}}
\put(20, 30.6){\line(0, 1){8.9}}
\put(20, 40.6){\line(0, 1){6.4}}
\put(19,2){$G_2$}
\put(30, 5){\line(0, 1){4.4}}
\put(30, 10.6){\line(0, 1){8.9}}
\put(30, 20.6){\line(0, 1){8.9}}
\put(30, 30.6){\line(0, 1){8.9}}
\put(30, 40.6){\line(0, 1){6.4}}
\put(29,2){$G_3$}
\put(40, 5.1){\line(0, 1){4.4}}
\put(40, 10.6){\line(0, 1){8.9}}
\put(40, 20.6){\line(0, 1){8.9}}
\put(40, 30.6){\line(0, 1){8.9}}
\put(40, 40.6){\line(0, 1){6.4}}
\put(39,2){$G_4$}
%
%
\put(8, 6){\line(1,1){5.5}}
\put(13.5,12.5){$E_{11}$}
\put(8, 16){\line(1,1){5.5}}
\put(13.5,22.5){$E_{12}$}
\put(8, 26){\line(1,1){5.5}}
\put(13.5,32.5){$E_{13}$}
\put(8, 36){\line(1,1){5.5}}
\put(13.5,42.5){$E_{14}$}
\put(18, 6){\line(1,1){5.5}}
\put(23.5,12.5){$E_{21}$}
\put(18, 16){\line(1,1){5.5}}
\put(23.5,22.5){$E_{22}$}
\put(18, 26){\line(1,1){5.5}}
\put(23.5,32.5){$E_{23}$}
\put(18, 36){\line(1,1){5.5}}
\put(23.5,42.5){$E_{24}$} 
\put(28, 6){\line(1,1){5.5}}
\put(33.5,12.5){$E_{31}$}
\put(28, 16){\line(1,1){5.5}}
\put(33.5,22.5){$E_{32}$}
\put(28, 26){\line(1,1){5.5}}
\put(33.5,32.5){$E_{33}$}
\put(28, 36){\line(1,1){5.5}}
\put(33.5,42.5){$E_{34}$} 
\put(38, 6){\line(1,1){5.5}}
\put(43.5,12.5){$E_{41}$}
\put(38, 16){\line(1,1){5.5}}
\put(43.5,22.5){$E_{42}$}
\put(38, 26){\line(1,1){5.5}}
\put(43.5,32.5){$E_{43}$}
\put(38, 36){\line(1,1){5.5}}
\put(43.5,42.5){$E_{44}$} 
\end{picture}
\end{center}
\vskip -.3cm
\caption{Double Kummer pencil}\label{fig:Km}
\end{figure}
\par
\smallskip
Let $B_{16}\subset \NS(\Bl{A})$ be the sublattice generated by the classes of
the sixteen $(-1)$-curves contracted by $\beta_A$.
Then we have 
$$
\NS(\Bl{A})=\NS(A)\perp B_{16}=U(A)\perp\Hom (E\sprime, E)[-1]\perp B_{16}.
$$
Since $\tiliota_A$ acts on $\NS(\Bl{A})$ trivially,
we see that $\phiA$ induces an isometry
\begin{equation*}\label{eq:phiA_NS}
\mapisomet{(\phiA)_*^+}{U(A)[2]\perp\Hom (E\sprime, E)[-2]\perp B_{16}[2]\;\;}{\;\;\NS(\Km(A))}
\end{equation*}
with a finite $2$-elementary cokernel. Hence we obtain the following:
\begin{proposition}\label{prop:NSKmHom}
We have $\rank (\NS(\Km(A)))=18 + \rank (\Hom (E\sprime, E))$.
\end{proposition}
Recall the basis $\xi$ and $\eta$ of $U(A)$ defined by~\eqref{eq:xieta}.
We  denote by 
$\tilde{\xi}\in \NS(\Km(A))$ and $ \tilde{\eta}\in \NS(\Km(A))$ the image of 
$\xi$ and $\eta$ by $(\varphi_{A})_*^+\circ (\beta_A)^*$,
where $(\beta_A)^*$ denotes the total transformation of divisors.
Then we have, for any $i$ and $j$,
\begin{equation}\label{eq:xietatil}
\tilde\xi\,=\,2[F_j]+\textsum_{\mu=1}^4  [E_{\mu j}] \quand
\tilde\eta\,=\,2[G_i]+\textsum_{\nu=1}^4 [E_{i\nu}], 
\end{equation}
and they are orthogonal to $[E_{ij}]$.
We have $\tilde\xi^2=\tilde\eta^2=0$ and $\tilde\xi\tilde\eta=2$.
We then put
$$
N(\Km(A))\;:=\;\ang{\tilde\xi,\tilde\eta}\perp\ang{\,[E_{ij}]\mid 1\le i, j\le 4\,}\;\subset\;  \NS(\Km(A)),
$$
which is the image of 
$U(A)[2]\perp B_{16}[2]$ by the isometry $(\phiA)_*^+$, and 
we denote  by $\overline{N}(\Km(A))\subset \NS(\Km(A))$ 
the primitive closure of $N(\Km(A))$.
\begin{proposition}\label{prop:N}
The lattice $\overline{ N}(\Km(A))$ is generated by the classes of $(-2)$-curves in the double Kummer pencil, and 
we have 
$[\overline{ N}(\Km(A)): N(\Km (A))]=2^7$.
In particular, we have $\disc (\overline{N}(\Km(A)))=-2^4$.
\end{proposition}
\begin{proof}
For simplicity, we put $N:=N(\Km (A))$ and $\overline{ N}:=\overline{N}(\Km (A))$.
Let $N\sprime\subset \NS(\Km(A))$  be the sublattice  
generated by $[E_{ij}]$, $[F_j]$ and  $[G_i]$.
It is obvious from the equalities~\eqref{eq:xietatil} 
that 
$N\sprime$ is contained in $\overline{N}$,
and 
it is easy to calculate that  $[N\sprime:N]=2^7$.
We will show that $N\sprime=\overline{N}$.
Let $\tilde\xi\dual$, $\tilde\eta\dual$ and $[E_{ij}]\dual (1\le i, j\le 4)$
be the basis of $N\dual$ dual to  the basis 
$\tilde\xi$, $\tilde\eta$ and $ [E_{ij}] \,(1\le i, j\le 4)$ of $N$.
The discriminant group $\Disc{N}=N\dual/N$
is isomorphic to $\F_2\sp{\oplus 18}$ with basis
$\tilde\xi\dual+N$, $\tilde\eta\dual+N$ and $[E_{ij}]\dual+N$.
With respect to this basis, we write an element of $\Disc{N}$ by
$[\,x, y\,|\, z_{11}, \dots, z_{44}]$ with $x, y, z_{ij}\in \F_2$.
Then $\Discform{N}:\Disc{N}\to\Q/2\Z$ is given by
$$
\Discform{N}([\,x, y\,|\,  z_{11}, \dots, z_{44}])=(\,xy-\weight([z_{11}, \dots, z_{44}])/2\,)\, +\, 2\Z,
$$
where $\weight([z_{11}, \dots, z_{44}])$ is the Hamming weight of
$[z_{11}, \dots, z_{44}]\in \F_2\sp{\oplus 16}$.
We put
$H\sprime:=N\sprime/N$ and  $\overline{H}:=\overline{N}/N$.
Then we have $H\sprime\subseteq \overline{H}$.
By Lemma~\ref{lem:nodal},
if a code word $[\,0,0\,|\,  z_{11}, \dots, z_{44}]$ is in $\overline{H}$,
then $\weight([z_{11}, \dots, z_{44}])\ne 4$ holds.
We can confirm  by computer that
every element $v$ of the finite abelian group $D_N$ 
of order $2^{18}$ satisfies the following:
if $v\notin H\sprime$, then
the linear code $\ang{H\sprime, v}\subset D_N$ 
spanned by $H\sprime$ and $v$
is either
not totally isotropic with respect to $\Discform{N}$,
or containing $[\,0,0\,|\,  z_{11}, \dots, z_{44}]$ with
$\weight([z_{11}, \dots, z_{44}])=4$.
Therefore $\overline{H}=H\sprime$ holds.
\end{proof}
\subsection{The transcendental lattice of $\Km(A)$}\label{subsec:TKm}
In this subsection,
we work over $\C$.
Then we have
$\cohom(\Bl{A}, \Z)\;=\;\cohom(A, \Z)\perp B_{16}$.
Since $\tiliota_A$ acts on $\cohom(\Bl{A}, \Z)$ trivially, we have an isometry
\begin{equation*}\label{eq:phiA_H}
\mapisomet{(\phiA)_*^+}{\cohom(A, \Z)[2]\perp B_{16}[2]}{ \cohom(\Km(A), \Z)}
\end{equation*}
with a finite $2$-elementary cokernel.
We put
\begin{eqnarray*}
P(A)&:=&(U(A)\perp B_{16}\isomet \cohom(\Bl{A}, \Z))\sperp\;=\;(U(A)\isomet \cohom(A, \Z))\sperp\quand\\
Q(\Km(A))&:=& (\overline{N}(\Km(A))\isomet \cohom(\Km(A),\Z))\sperp.
\end{eqnarray*}
\begin{proposition}\label{prop:PQ}
The isometry $(\phiA)_*^+$ induces the following commutative diagram,
in which the horizontal isomorphisms of lattices preserve the Hodge structure:
\begin{equation}\label{eq:PQTT}
\renewcommand{\arraystretch}{1.2}
\begin{array}{ccc}
T(A)[2] & \wcong & T(\Km(A))\\
\downinj & & \downinj \\
P(A)[2] & \wcong & Q(\Km(A)).
\end{array}
\end{equation}
\end{proposition}
\begin{proof}
First we prove that $(\phiA)_*^+$ induces
$P(A)[2]  \cong  Q(\Km(A))$.
By the definition of $N(\Km(A))$,
the isometry $(\phiA)_*^+$ maps
$(U(A)\perp B_{16})[2]$ to $N(\Km(A))$
isomorphically,
and hence
$(\phiA)_*^+$ induces an isometry
from $P(A)[2]$ to $ Q(\Km(A))$
with a finite $2$-elementary cokernel.
Since $U(A)\perp B_{16}$ and $\cohom(\Bl{A}, \Z)$  are unimodular,
we have $\disc( P(A)[2])=2^4$.
Since $\cohom(\Km(A),\Z)$ is unimodular and 
$\disc (\overline{N}(\Km(A)))=-2^4$ by Proposition~\ref{prop:N},
we have $\disc( Q(\Km(A)))=2^4$.
Therefore the isometry $P(A)[2]\isomet Q(\Km(A))$
is in fact an isomorphism.
\par
By definition, we have 
 $T(\Bl{A})\subset P(A)$ and 
$T(\Km(A))\subset Q(\Km(A))$.
Since $(\phiA)_*^+$ preserves the Hodge structure,
the isomorphism $P(A)[2]  \cong  Q(\Km(A))$
induces $T(A)[2] \cong  T(\Km(A))$.
\end{proof}
\begin{remark}
The isomorphism  $T(A)[2]  \cong  T(\Km(A))$ was proved in~\cite[\S4]{MR0284440}
using the sublattice $\ang{[E_{ij}]\mid 1\le i, j\le 4}$
instead of $N(\Km(A))$.
We need the diagram~\eqref{eq:PQTT}
for the proof of Proposition~\ref{prop:LKm}.
\end{remark}
\subsection{The supersingular reduction lattice of $\Km(A)$}\label{subsec:LKm}
Let $W$ be either a number field, or a Dedekind domain with the quotient field $F$ being
a number field.
We assume that $2$ is invertible in $W$.
Let $\EEE\sprime$ and $\EEE$ be smooth proper families of
elliptic curves over $U:=\Spec W$.
We put
$\AAA:=\EEE\sprime\times_U\EEE$.
\begin{definition}
A diagram
$$
(\KKmm)\;\;:\;\;  \Km(\AAA)\mapleftsp{} \Bl{\AAA}\maprightsp{} \AAA=\EEE\sprime\times_U\EEE
$$
of schemes and morphisms over $U$ is called 
the \emph{Kummer diagram} over $U$ of $\EEE\sprime$ and $\EEE$
if the following hold:
\begin{rmenumerate}
\item
$\Km(\AAA)$ and $\Bl{\AAA}$ are smooth and proper over $U$,
\item
$\Bl{\AAA}\to \AAA$ is the blowing-up
along the fixed locus (with the reduced structure)
of the inversion automorphism $\iota_{\AAA}:\AAA\isom\AAA$
over $U$, and 
\item 
$\Km(\AAA)\leftarrow \Bl{\AAA}$ is the quotient morphism by a lift $\tiliota_{\AAA}$
of $\iota_{\AAA}$.
\end{rmenumerate}
%
%
\end{definition}
In this subsection, we consider the case where $W$ is a Dedekind domain.
\par
\smallskip
Suppose that the Kummer diagram $(\KKmm)$
over $U$ of $\EEE\sprime$ and $\EEE$ is given.
Then, at every point $P$ of $U$ (closed or generic, see the definition~\eqref{eq:kappa}), 
the diagram $(\KKmm)\otimes\bar\kappa_P$ is the Kummer diagram
of the elliptic curves $\EEE\sprime\otimes \bar\kappa_P$ and $\EEE\otimes \bar\kappa_P$.
\par
\smallskip
Let $\ppp$ be a closed point of  $U$
with $\kappa:=\kappa\sbppp$ 
being  of characteristic $p$.
Note that $p \ne 2$ by the assumption $1/2\in W$.
We put
\begin{equation}\label{eq:EEAEEA}
\renewcommand{\arraystretch}{1.2}
\begin{array}{lll}
E\sprime:=\EEE\sprime\otimes \overline {F},\quad&
E:=\EEE\otimes \overline{F},\quad &
A:=E\sprime\times E=\AAA\otimes\overline{F}\quad \quand 
\\
E\sprime_0:=\EEE\sprime\otimes \bar \kappa,\quad& 
E_0:=\EEE\otimes \bar \kappa,\quad&
A_0:=E_0\sprime\times E_0=\AAA\otimes\bar \kappa.
\end{array}
\end{equation}
Then we have 
$\Km(\AAA) \otimes \overline {F}=\Km(A)$ and 
$\Km(\AAA) \otimes \bar \kappa=\Km(A_0)$.
We  assume that
\begin{equation*}\label{eq:rankHom}
\rank (\Hom (E\sprime, E))=2\quand
\rank (\Hom(E\sprime_0, E_0))=4.
\end{equation*}
Then, by Proposition~\ref{prop:NSA},  we have
$\rank (\NS(A))=4$ and 
$\rank (\NS(A_0))=6$.
By Proposition~\ref{prop:NSKmHom}, 
we see that
$\Km (A)$ is singular
and  $\Km(A_0)$ is supersingular.
We consider the supersingular reduction lattices 
\begin{eqnarray*}
L(\AAA, \ppp) &:=& (\NS(A)\isomet \NS(A_0))\sperp \quand\\
L(\Km(\AAA), \ppp) &:=& (\NS(\Km(A))\isomet \NS(\Km(A_0)))\sperp.
\end{eqnarray*}
Note that, by Proposition~\ref{prop:NSA}, we have
\begin{equation}\label{eq:LHom}
L(\AAA, \ppp)= (\Hom (E\sprime, E)\isomet \Hom(E_0\sprime, E_0))\sperp[-1].
\end{equation}
\begin{proposition}\label{prop:LKm}
Suppose that $p$ is prime to $\disc(\NS(\Km(A)))$.
Then the Kummer diagram $(\KKmm)$ induces an isomorphism
$L(\AAA,\ppp)[2]\cong L(\Km(\AAA),\ppp)$.
\end{proposition}
We use the following lemma,
whose proof is quite elementary and is omitted:
\begin{lemma}\label{lem:perp}
Let 
 $\rho: \lat_1\isomet\lat_2$  be an isometry of lattices.
Suppose that $\rho$ maps a sublattice $N_1$ 
of $\lat_1$ to a sublattice $N_2$ of $\lat_2$.
For $i=1,2$, we put $M_i:=(N_i\isomet \lat_i)\sperp$.
If $\rank (N_1)=\rank (N_2)$, then $\rho$ maps $M_1$ into $M_2$
and we have
$$
(\lat_1\isomet\lat_2)\sperp\;=\; (M_1\isomet M_2)\sperp.
$$
\end{lemma}
\begin{proof}[Proof of Proposition~\ref{prop:LKm}]
Note that the double Kummer pencil is defined on $\Km(\AAA)$ and is flat over $U$.
Hence, by Proposition~\ref{prop:N},
the specialization isometry 
$$
\mapisomet{\rho_{\Km(A)}}{ \NS(\Km(A))}{\NS(\Km(A_0))}
$$
maps $\overline{N} (\Km(A))$ to $\overline{N} (\Km(A_0))$ isomorphically.
We put
\begin{eqnarray*}
\Sigma (\Km(A)) &:=& (\overline{N} (\Km(A))\isomet \NS(\Km(A)))\sperp \quand\\
\Sigma (\Km(A_0)) &:=& (\overline{N} (\Km(A_0))\isomet \NS(\Km(A_0)))\sperp.
\end{eqnarray*}
Then  $\rho_{\Km (A)}$ maps 
$\Sigma (\Km(A)) $ to $\Sigma (\Km(A_0)) $, and we have
\begin{equation}\label{eq:LSigmaSigma}
L(\Km(\AAA), \ppp) \;=\; (\Sigma (\Km(A))\isomet \Sigma (\Km(A_0)))\sperp
\end{equation}
by Lemma~\ref{lem:perp}.
The isometry $(\phiA)_*\sp+$
maps $(U(A)\perp B_{16})[2]\subset \NS(A)[2]$ to $N(\Km(A))$ isomorphically.
Hence $(\phiA)_*\sp+$ 
induces an isometry
\begin{equation}\label{eq:HomSigma}
\Hom(E\sprime, E)[-2]\;\isomet\; \Sigma (\Km(A))
\end{equation}
with a finite $2$-elementary cokernel.
In the same way, we see that 
$(\varphi_{A_0})_*\sp+$ induces an isometry
\begin{equation}\label{eq:HomSigmazero}
\Hom(E_0\sprime, E_0)[-2]\;\isomet\; \Sigma (\Km(A_0))
\end{equation}
with a finite $2$-elementary cokernel.
By the equalities~\eqref{eq:LHom} and~\eqref{eq:LSigmaSigma},
it is enough to show that the isometries~\eqref{eq:HomSigma} and~\eqref{eq:HomSigmazero}
are isomorphisms.
\par
To show that~\eqref{eq:HomSigma} is an  isomorphism,
we choose an embedding $\sigma$ of $\overline{F}$ into $\C$,
and consider the analytic manifolds
$E\sp{\prime\sigma}$, $E\sp{\sigma}$, ${\Bl{A}}\sp{\sigma}$, 
$A\sp{\sigma}$ and $\Km(A\sp{\sigma})=\Km(A)\sp{\sigma}$.
By Proposition~\ref{prop:NSA}, we have
$$
\Hom(E\sp{\prime\sigma}, E\sp{\sigma})[-1]\;=\;
P(A\sp{\sigma})\cap \NS(A\sp{\sigma}) \;=\;
(T(A\sp{\sigma})\isomet P(A\sp{\sigma}))\sperp,
$$
where $P(A\sp{\sigma})$ is the lattice defined in the previous subsection.
By the definition of $\Sigma(\Km(A)\sp{\sigma})$, we have
$$
\Sigma(\Km(A)\sp{\sigma})\;=\;
Q(\Km(A)\sp{\sigma})\cap \NS(\Km(A)\sp{\sigma}) \;=\;
(T(\Km(A)\sp{\sigma})\isomet Q(\Km(A)\sp{\sigma}))\sperp,
$$
where $Q(\Km(A)\sp{\sigma})$ is the lattice defined in the previous subsection.
Therefore, from~\eqref{eq:PQTT}, we see that
the analytic Kummer diagram $(\KKmm)\sp{\sigma}$
induces 
$$
\Hom(E\sp{\prime\sigma}, E\sp{\sigma})[-2]\wcong \Sigma (\Km(A)\sp{\sigma}).
$$
Hence the isometry~\eqref{eq:HomSigma} is an  isomorphism.
\par
Since $p\notdiv 2\disc(\NS(\Km(A)))$,
the Artin invariant of $\Km(A_0)$ is $1$ by Proposition~\ref{prop:normalK3},
and hence $\disc(\NS(\Km(A_0)))$ is equal to $-p^2$.
By Proposition~\ref{prop:N},
we have
$\disc(\overline{N} (\Km(A_0)))=-2^4$,
and hence
we obtain
$\disc(\Sigma(\Km(A_0)))=2^4 p^2$
by Proposition~\ref{prop:nikulin}.
On the other hand, 
we have
$\disc(\Hom(E\sprime_0, E_0)[-2])=2^4 p^2$
by~Proposition~\ref{prop:psq}.
Comparing the discriminants, we conclude  that the isometry~\eqref{eq:HomSigmazero} is an  isomorphism.
\end{proof}
\section{Shioda-Inose construction}\label{sec:ShiodaInose}
We continue to denote by $k$  an algebraically closed field  of characteristic $\ne 2$.
\subsection{Shioda-Inose configuration}\label{subsec:SIconfig}
Let $Z$ be a $K3$ surface defined over $k$.
\begin{definition}
We say that a pair $(C, \Theta)$ of reduced effective divisors  on $Z$
is a \emph{Shioda-Inose configuration}
if the following hold:
\begin{rmenumerate}
\item $C$ and $\Theta$ are disjoint, 
\item
$C=C_1+\dots+C_8$ is  an $ADE$-configuration of $(-2)$-curves of type $\E_8$,
\item $\Theta=\Theta_1+\dots+\Theta_8$ is  an $ADE$-configuration of $(-2)$-curves of type $8 \A_1$, 
 \item\label{condSI:double}
there exists a class $[\LLL]\in \NS(Z)$ such that
$2[\LLL]=[\Theta]$.
\end{rmenumerate}
\end{definition}
Let  $(C, \Theta)$ be a Shioda-Inose configuration on $Z$.
Then there exists a finite double covering $\varphi_{Y}: \Bl{Y}\to Z$
that branches exactly along $\Theta$
by the condition~(\ref{condSI:double}).
Let $T_i\subset \Bl{Y}$ be the reduced part of the pull-back of $\Theta_i$ by $\phiY$,
which is a $(-1)$-curve on $\Bl{Y}$, 
and let
$\beta_Y:\Bl{Y}\to  Y$ be the contraction of $T_1, \dots, T_8$.
Then $Y$ is a $K3$ surface.
Let $\tiliota_Y$ be the deck-transformation of $\Bl{Y}$ over $Z$.
Then $\tiliota_Y$ is the  lift of an  involution $\iota_Y: Y\isom Y$ of $Y$,
which has  eight fixed points.
\begin{definition}
The diagram
$$
(\SSII)\;\;:\;\;Y\mapleftsp{\beta_Y}\Bl{Y}\maprightsp{\phiY} Z
$$
is called the \emph{Shioda-Inose diagram}
associated with the Shioda-Inose configuration
$(C, \Theta)$ on $Z$.
We call $Y$ a \emph{Shioda-Inose surface} 
of  $Z$.
\end{definition}
We denote by $\Gamma\subset \NS(Z)$
the sublattice generated by $[C_1], \dots, [C_8]$.
Then $\Gamma$ is a negative-definite root lattice of type $\E_8$.
In particular, $\Gamma$ is unimodular.
We then denote by
$M(Z)\subset \NS(Z)$ the sublattice generated by
$[\Theta_1], \dots, [\Theta_8]$,
and by $\overline{M}(Z)$ the primitive closure of $M(Z)$ in $\NS(Z)$.
Note that $\overline{M}(Z)$ and $\Gamma$ are orthogonal in $\NS(Z)$.
By the condition~\eqref{condSI:double} of the Shioda-Inose configuration, $\overline{M}(Z)$
contains  the class
$[\LLL]=[\Theta]/2$.
In fact, we can see that $\overline{M}(Z)$ is generated by $M(Z)$ and $[\LLL]$
from the following equality~\cite[Lemma 3.4]{MR0441982}:
\begin{equation}\label{eq:M}
\disc( \overline{M}(Z))=2^6.
\end{equation}
Note that one can also  prove the equality~\eqref{eq:M} easily  using  Lemma~\ref{lem:nodal}
in the same way as the proof Proposition~\ref{prop:N}.
We then put
\begin{equation}\label{eq:Xi}
\Xi(Z)\;:=\;(\Gamma\perp \overline{M} (Z)\isomet \NS(Z))\sperp.
\end{equation}
Next we define several sublattices of $\NS(\Bl{Y})$ and $\NS(Y)$.
Since $\phiY$ is \'etale in a neighborhood of $C\subset Z$ and $C$ is simply-connected,
the pull-back  of $C$ by $\phiY$ consists of 
two connected components $C^{[1]}$ and $C^{[2]}$.
Let  $\Gamma^{[1]}$ and $\Gamma^{[2]}$ be the sublattices of $\NS(\Bl{Y})$
generated by the classes of the irreducible components of
$C^{[1]}$ and $C^{[2]}$, respectively.
We  put
$$
\wt{\Gamma}:=\Gamma^{[1]}\perp\Gamma^{[2]}.
$$
The sublattice  $\wt{\Gamma}$ is mapped by $\beta_{Y}$ isomorphically to a  sublattice of $\NS(Y)$,
which we will denote by the same letter $\wt{\Gamma}$.
We  denote by $B_8\subset\NS(\Bl{Y})$ the sublattice generated by
the classes of the $(-1)$-curves $[T_1], \dots, [T_8]$.
Then $B_8$ is orthogonal to $\wt{\Gamma}$, 
and we have a canonical isomorphism 
$\NS(\Bl{Y})\cong \NS(Y)\perp B_8$.
We put
\begin{equation}\label{eq:Pi}
\Pi(Y)\;:=\;(\wt{\Gamma}\perp B_8\isomet \NS(\Bl{Y}))\sperp=(\wt{\Gamma}\isomet \NS(Y))\sperp.
\end{equation}
Since $\wt{\Gamma}$ and $B_8$ are unimodular, we have 
\begin{equation}\label{eq:NSY}
\NS(\Bl{Y})=\wt{\Gamma}\perp B_8\perp \Pi(Y)
\quand
\NS(Y)=\wt{\Gamma}\perp \Pi(Y).
\end{equation}
The action of $\tiliota_Y$ on $\NS(\Bl{Y})$
and the action of $\iota_Y$ on $\NS(Y)$
preserve the orthogonal direct-sum decompositions~\eqref{eq:NSY},
and the action of $\tiliota_Y$ is trivial on $B_8$.
We put
$$
\wt{\Gamma}\sp+\;:=\;\wt{\Gamma}\cap \wt{\Gamma}_{ \Q}\sp+
\quand 
\Pi(Y)\sp+\;:=\;\Pi(Y)\cap \Pi(Y)_{\Q}\sp+,
$$
where $\wt{\Gamma}_{ \Q}\sp+$ (resp.~$\Pi(Y)_{\Q}\sp+$) is the eigenspace of 
$(\tiliota_Y)_*$ 
on $\wt{\Gamma}\otimes\Q$ (resp.~$\Pi(Y)\otimes\Q$) with the eigenvalue $1$.
Since  $\tiliota_Y$ acts on $\wt{\Gamma}$ by interchanging   $\Gamma^{[1]}$ and $\Gamma^{[2]}$,
we have
$\rank (\wt{\Gamma}\sp+)=8$.
By Lemma~\ref{lem:double},
we see that $\varphi_{Y}$ induces an isometry
\begin{equation*}\label{eq:toNSZ}
\mapisomet{(\varphi_{Y})_*\sp+}{\wt{\Gamma}\sp+[2]\perp B_8[2]\perp \Pi(Y)\sp+[2]\;}{\; \NS(Z)}
\end{equation*}
with a finite $2$-elementary cokernel.
Since  $(\varphi_{Y})_*\sp+$  induces an isometry
$\wt{\Gamma}\sp+[2]\isomet \Gamma$
with  a finite $2$-elementary cokernel
and an isomorphism  $ B_8[2]\cong M(Z)$,
we obtain the following:
\begin{proposition}\label{prop:NSZY}
{\rm (1)}
We have
$$
\rank (\NS(Y))=16+\rank(\Pi(Y)) \;\ge \; \rank (\NS(Z))=16 +\rank(\Pi(Y)\sp+).
$$

{\rm (2)}
If $\rank (\NS(Z))$ is equal to $\rank (\NS(Y))$, then 
$(\phiY)\sb*\sp+$ induces an isometry
$\Pi(Y)[2]\isomet\Xi(Z)$
with a finite $2$-elementary cokernel.
\end{proposition}
\subsection{The transcendental lattice of the  Shioda-Inose surface}\label{subsec:TSI}
In this subsection, we work over $\C$.
Note that we have 
$\cohom(\Bl{Y},\Z)=\cohom (Y, \Z)\perp B_8$.
We consider the isometry
$$
\mapisomet{(\phiY)_*^+}{\cohom(\Bl{Y}, \Z)\sp +[2]}{\cohom(Z, \Z)}
$$
with a finite $2$-elementary cokernel.
We put
\begin{eqnarray*}
R(Y)&:=&(\wt{\Gamma}\perp B_8 \isomet \cohom(\Bl{Y}, \Z))\sperp\;=\;
(\wt{\Gamma} \isomet \cohom(Y, \Z))\sperp \quand\\
S(Z)&:=&(\Gamma\perp \overline{M}(Z)\isomet \cohom(Z, \Z))\sperp.
\end{eqnarray*}
\begin{proposition}\label{prop:RS}
The isometry $(\phiY)_*^+$ induces the following commutative diagram,
in which the horizontal  isomorphisms of lattices preserve the Hodge structure:
\begin{equation}\label{eq:RSTT}
\renewcommand{\arraystretch}{1.2}
\begin{array}{ccc}
T(Y)[2] & \wcong & T(Z)\\
\downinj & & \downinj \\
R(Y)[2] & \wcong & S(Z).
\end{array}
\end{equation}
\end{proposition}
\begin{proof}
First we prove $R(Y)[2] \cong S(Z)$.
Since $\wt{\Gamma}$ and $B_8$ are unimodular,
we have
\begin{equation}\label{eq:decompY}
\cohom(\Bl{Y}, \Z)=\wt{\Gamma}\perp B_8\perp R(Y).
\end{equation}
The action of $\tiliota_Y$ on $\cohom(\Bl{Y}, \Z)$ preserves the decomposition~\eqref{eq:decompY},
and is trivial on $B_8$.
Since $\rank (\wt{\Gamma}\sp +)=8$ and $\rank (\cohom(\Bl{Y}, \Z)\sp +)=22$, 
we see that $\tiliota_Y$ acts on $R(Y)$ trivially.
(This fact was also proved in~\cite[Lemma 3.2]{MR0441982}.)
We thus obtain an isometry
$$
\mapisomet{(\phiY)_*\sp + }{\wt{\Gamma}\sp + [2] \perp B_8[2]\perp R(Y)[2] }{ \cohom (Z, \Z)}
$$
with a finite $2$-elementary cokernel.
Since $(\phiY)_*\sp + $ maps $ \wt{\Gamma}\sp + [2]$ to $\Gamma$ and $B_8[2]$ to $\overline{M}(Z)$
with finite $2$-elementary cokernels,
it induces an isometry from $R(Y)[2]$ to $S(Z)$ with a finite $2$-elementary cokernel.
From the decomposition~\eqref{eq:decompY},  we have $\disc (R(Y)[2])=-2^6$.
Since  $\disc (\overline{M}(Z))=2^6$ by the equality~\eqref{eq:M},
we have $\disc (S(Z))=-2^6$.
Therefore the isometry $R(Y)[2] \isomet S(Z)$ is in fact an isomorphism.
\par
The proof of the isomorphism $T(Y)[2]  \cong  T(Z)$
is completely parallel to  the second paragraph of the proof of Proposition~\ref{prop:PQ}.
\end{proof}
\begin{remark}
The isomorphism  $T(Y)[2]  \cong  T(Z)$ is due to Shioda and Inose~\cite{MR0441982}.
We need the diagram~\eqref{eq:RSTT}
for the proof of Proposition~\ref{prop:LSI}.
\end{remark}
\begin{corollary}\label{cor:SIoverC}
In characteristic $0$, 
we have
$\rank (\NS(Y))=\rank (\NS(Z))$ and 
$2^{22-r} \disc (\NS(Y))= \disc (\NS(Z))$, 
where $r:=\rank (\NS(Y))=\rank (\NS(Z))$.
\end{corollary}
\subsection{The supersingular reduction lattice of the  Shioda-Inose surface}\label{subsec:LSI}
Let $W$ be either a number field, or a Dedekind domain with   the quotient field $F$ being 
a number field
such that $1/2\in W$.
Let $\ZZZ$ be a smooth proper family of $K3$ surfaces over $U:=\Spec W$.
\begin{definition}
A diagram
$$
(\SSII)\;\;:\;\; \YYY\mapleftsp{} \Bl{\YYY}\maprightsp{} \ZZZ
$$
of schemes and morphisms over $U$ is called 
a \emph{Shioda-Inose  diagram} over $U$ 
if there exists a pair of  reduced effective divisors $(\CCC, \varTheta)$ of $\ZZZ$
such that the following hold:
\begin{rmenumerate}
\item $\CCC$ and $\varTheta$ are flat over $U$, 
\item
$\YYY$ and $\Bl{\YYY}$ are smooth and proper over $U$,
\item at each point $P$ of $U$ (closed or generic, see the definition~\eqref{eq:kappa}),
the pair of divisors $(\CCC\otimes\bar\kappa_P, \varTheta\otimes\bar\kappa_P)$
is a Shioda-Inose configuration on $\ZZZ\otimes\bar\kappa_P$,
\item
$\Bl{\YYY}\to \ZZZ$ is 
a finite double covering  that branches  exactly along $\varTheta$, and 
\item 
$\YYY\leftarrow \Bl{\YYY}$ is a contraction of the inverse image of $\varTheta$ in $\Bl{\YYY}$.
\end{rmenumerate}
%
\end{definition}
In this subsection,
we consider the case where $W$ is a Dedekind domain.
\par
\smallskip
Suppose that a Shioda-Inose  diagram $(\SSII)$
over $U$  is given.
Then, at every point $P$ of $U$,  
the diagram $(\SSII)\otimes\bar\kappa_P$ is a Shioda-Inose diagram
of $\ZZZ\otimes\bar\kappa_P$.
\par
\smallskip
Let $\ppp$ be a closed point of  $U$
with  $\kappa:=\kappa\sbppp$ being  of characteristic $p$.
Note that $p \ne 2$ by the assumption $1/2\in W$.
We put
\begin{equation}\label{eq:YYZYYZ}
\renewcommand{\arraystretch}{1.4}
\begin{array}{lll}
Y:=\YYY\otimes\overline{F}, & \Bl{Y}:=\Bl{\YYY}\otimes\overline{F}, & Z:=\ZZZ\otimes\overline{F},\\
Y_0:=\YYY\otimes\bar \kappa, & \Bl{Y}_0:=\Bl{\YYY}\otimes\bar \kappa, & Z_0:=\ZZZ\otimes\bar \kappa.
\end{array}
\end{equation}
We assume that
$Z$ is  singular  and 
$Z_0$ is  supersingular.
Then $Y$ is singular and $Y_0$ is supersingular 
by Proposition~\ref{prop:NSZY}.
We consider the supersingular reduction lattices 
$$
L(\ZZZ,\ppp):=(\NS(Z)\isomet \NS(Z_0))\sperp\;\;\quand
\quad
L(\YYY,\ppp):=(\NS(Y)\isomet \NS(Y_0))\sperp.
$$
By Corollary~\ref{cor:SIoverC},
we have
$4\disc (\NS(Y))= \disc (\NS(Z))$.
Since $p$ is odd, the following are equivalent:
(i)  $p\notdiv  \disc (\NS(Z))$, and  (ii) $p\notdiv  \disc (\NS(Y))$.
\begin{proposition}\label{prop:LSI}
Suppose that $p$  satisfies the conditions {\rm (i)} and {\rm (ii)} above.
Then the Shioda-Inose diagram $(\SSII)$
induces an isomorphism
$L(\YYY,\ppp)[2]\cong L(\ZZZ,\ppp)$.
\end{proposition}
\begin{proof}
Recall the definition~\eqref{eq:Pi} of $\Pi$.
Since the specialization isometry
$\NS(Y)\isomet \NS(Y_0)$ maps
$\wt{\Gamma}\subset \NS(Y)$ to $\wt{\Gamma}\subset \NS(Y_0)$
isomorphically,
it maps $\Pi(Y)$ to $\Pi(Y_0)$, and we have
$$
L(\YYY,\ppp)\;=\;(\Pi(Y)\isomet \Pi(Y_0))\sperp
$$
by Lemma~\ref{lem:perp}.
Recall the definition~\eqref{eq:Xi} of $\Xi$.
Since the specialization isometry
$\NS(Z)\isomet \NS(Z_0)$ maps
$\Gamma\subset \NS(Z)$ to $\Gamma\subset \NS(Z_0)$
and 
$M(Z)\subset \NS(Z)$ to $M(Z_0)\subset \NS(Z_0)$
isomorphically,
it maps $\Xi(Z)$ to $\Xi(Z_0)$, and we have
$$
L(\ZZZ,\ppp)\;=\; (\Xi(Z)\isomet \Xi(Z_0))\sperp
$$
by Lemma~\ref{lem:perp}.
By Proposition~\ref{prop:NSZY}, we have isometries
\begin{equation}\label{eq:PiXi}
\Pi(Y)[2]\;\isomet\;\Xi(Z)\quand 
\Pi(Y_0)[2]\;\isomet\;\Xi(Z_0)
\end{equation}
with  finite $2$-elementary cokernels
induced by $(\SSII)\otimes\overline{F}$ and $(\SSII)\otimes\bar{\kappa}$,
respectively.
It is therefore enough to show that
both of the isometries in~\eqref{eq:PiXi} are isomorphisms.
\par
We choose an embedding $\sigma$ of $\overline{F}$ into $\C$,
and consider the transcendental lattices
$T(Y\sp{\sigma})$ and $T(Z\sp{\sigma})$.
We have
\begin{eqnarray*}
\Pi(Y) &\cong&\Pi(Y\sp{\sigma})\;=\; R(Y\sp{\sigma})\cap \NS(Y\sp{\sigma}) \;=\;(T(Y\sp{\sigma})\isomet R(Y\sp{\sigma}))\sperp\quand\\
\Xi(Z) &\cong &\Xi(Z\sp{\sigma})\,\;=\; S(Z\sp{\sigma})\cap \NS(Z\sp{\sigma})\,\;=\; (T(Z\sp{\sigma})\isomet S(Z\sp{\sigma}))\sperp,
\end{eqnarray*}
where $R(Y\sp{\sigma})$ and $S(Z\sp{\sigma})$ are the lattices defined in the previous subsection.
Since the analytic Shioda-Inose diagram $(\SSII)\sp{\sigma}$
induces  the commutative diagram~\eqref{eq:RSTT} 
for $Y\sp\sigma$ and $Z\sp\sigma$,
we see that the first isometry of~\eqref{eq:PiXi} is 
an isomorphism.
\par
%
%
By~Proposition~\ref{prop:normalK3}, we have
$\disc (\NS(Y_0))=\disc (\NS(Z_0))=-p^2$.
Since $\wt{\Gamma}$ is unimodular, we have
$\disc (\Pi (Y_0))= -p^2$ by Proposition~\ref{prop:nikulin},
and hence 
$\disc (\Pi (Y_0)[2])$ is equal to $-2^6 p^2$.
Since $\Gamma$ is unimodular and $\disc (\overline{M}(Z_0))=2^6$ by the equality~\eqref{eq:M},
we see that 
$\disc (\Xi (Z_0))$ is equal to $-2^6 p^2$
by Proposition~\ref{prop:nikulin} again.
Therefore 
the second isometry of~\eqref{eq:PiXi} is 
also an isomorphism.
\end{proof}
\section{Proof of Theorems~\ref{thm:almostallp} and~\ref{thm:toF}}\label{sec:proofs}
\subsection{Preliminaries}\label{subsec:FGA}
In this subsection,
we quote fundamental facts in algebraic geometry from
Grothendieck's FGA~\cite[no.~221]{MR0146040}.
See also~\cite[Chap.~5]{MR2222646}.
\par
\smallskip
Let $S$ be a noetherian scheme,
and let $\WWW$ and $\ZZZ$ be schemes flat and projective over $S$.
We denote by $\fMor_S(\WWW, \ZZZ)$ the functor from the category $\Sch_S$ of 
locally noetherian schemes over $S$ 
to the category of sets such that,
for an object  $T$ of $\Sch_S$, we have 
$$
\fMor_S(\WWW, \ZZZ)(T)=\textrm{the set of $T$-morphisms from $\WWW\times_S T$ to $\ZZZ\times_S T$}.
$$
Then we have the following (\cite[no.~221, Section 4]{MR0146040}, \cite[Theorem 5.23]{MR2222646}):
\begin{theorem}
The functor $\fMor_S(\WWW, \ZZZ)$ is representable by 
an open subscheme $\Mor_S(\WWW,\ZZZ)$ of the Hilbert scheme
$\Hilb_{\WWW\times_S\ZZZ/S}$ parameterizing  
 closed subschemes of  $\WWW\times_S\ZZZ$ flat over $S$. 
\end{theorem}
Let $F$ be a number field,
and let $X$ and $Y$ be smooth projective varieties defined over $F$.
By the flattening  stratification~(\cite[Theorem 5.12]{MR2222646}, \cite[Lecture 8]{MR0209285}),
we have a non-empty open subset $U$ of $\Spec \Z_F$ and smooth projective $U$-schemes 
$\XXX$ and $\YYY$ such that
the generic fibers $\XXX\times_U F$ and $\YYY\times_U F$
are isomorphic to $X$ and $Y$, respectively.
We will consider the schemes 
$$
\Mor_{V}(\XXX_V, \YYY_V)=\Mor_{U}(\XXX, \YYY)\times_U V
$$
for non-empty open subsets $V$ of $U$,
where
$\XXX_V:=\XXX\times_U V$ and $\YYY_V:=\YYY\times_U V$.
\begin{proposition}\label{prop:morext}
Let $\varphi : X\to Y$ be an $F$-morphism.
Then there exist a non-empty open subset $V\subset U$
and a $V$-morphism 
$\wt{\varphi }_V:\XXX_V\to \YYY_V$  that extends $\varphi $.
If $\wt{\varphi }_{V\sprime}:\XXX_{V\sprime}\to \YYY_{V\sprime}$ 
is a morphism over a non-empty  open subset $V\sprime \subset  U$ that extends $\varphi $,
then $\wt{\varphi }_{V}|_{V\cap V\sprime}=\wt{\varphi }_{V\sprime}|_{V\cap V\sprime}$ holds,
where $\wt{\varphi }_{V}|_{V\cap V\sprime}$ and $\wt{\varphi }_{V\sprime}|_{V\cap V\sprime}$
denote the restrictions of  $\wt{\varphi }_{V}$ and $\wt{\varphi }_{V\sprime}$
 to $\XXX_{V\cap V\sprime}$.
\end{proposition}
\begin{proof}
We denote by $[\varphi]:\Spec F\to \Mor_U(\XXX,\YYY)$
the $U$-morphism corresponding to $\varphi : X\to Y$.
Let $\Phi$ be the Hilbert polynomial of the graph $\Gamma(\varphi)\subset X\times_F Y$
of $\varphi$
with respect to a relatively ample invertible sheaf $\OOO(1)$ of 
$\XXX\times _U \YYY\to U$, 
so that $[\varphi]$ is an $F$-rational point of
$\Mor_U(\XXX,\YYY)\cap H^\Phi$,
where $H^\Phi:=\Hilb^\Phi_{\XXX\times _U \YYY/U}$
is  the Hilbert scheme
 parameterizing  
 closed subschemes of  $\XXX\times_U\YYY$ flat over $U$
with the Hilbert polynomial of fibers with respect to $\OOO(1)$
being  equal to $\Phi$.
Since $H^\Phi$ is projective over $U$, the morphism $[\varphi]$
extends to a morphism $[\varphi]\sp\sim_U: U\to H^\Phi$.
Since $\Mor_U(\XXX,\YYY)\cap H^\Phi$ is open in $H^\Phi$,
there exists a non-empty open subset $V$ of $U$ such that 
$[\varphi]$ extends to a $U$-morphism 
$$
\map{[\varphi]\sp\sim_V }{ V}{ \Mor_U(\XXX,\YYY)}.
$$
Hence the existence of a morphism  $\wt{\varphi }_V:\XXX_V\to\YYY_V$
extending $\varphi$  over some non-empty open subset $V\subset U$  is proved.
The equality  $\wt{\varphi }_{V}|_{V\cap V\sprime}=\wt{\varphi }_{V\sprime}|_{V\cap V\sprime}$
follows from the fact that  $H^\Phi\to U$ is separated.
\end{proof}
We call $\wt{\varphi}_V$ the \emph{extension of $\varphi$ over $V$}.
By the uniqueness of the extension, we obtain the following:
\begin{corollary}\label{cor:compext}
Let $\ZZZ$ be a smooth projective $U$-scheme
with the generic fiber  $Z$.
Let $\psi: Y\to Z$
be an $F$-morphism,
and let $\wt{\psi}_{V\sprime}: \YYY_{V\sprime}\to\ZZZ_{V\sprime}$ be the extension of $\psi$ 
over a non-empty open subset  ${V\sprime}\subset U$.
Then 
$(\wt{\psi}_{V\sprime}|_{V\cap V\sprime})\circ (\wt{\varphi}_V|_{V\cap V\sprime})$
is the extension of $\psi\circ \varphi: X\to Z$ over $V\cap V\sprime$.
\end{corollary}
Applying Corollary~\ref{cor:compext} to an $F$-isomorphism 
and its inverse,
we obtain the following,
which plays a key role in the proof of Theorem~\ref{thm:almostallp}:
\begin{corollary}\label{cor:isomV}
If $X$ and $Y$ are isomorphic over $F$,
then there exists a non-empty open subset $V\subset U$
such that $\XXX_V$ and $\YYY_V$ are isomorphic over $V$.
\end{corollary}
We give three  applications that will be used in 
the proof of Proposition~\ref{prop:ext}.
\begin{example}\label{example:blowup}
Let $Q$ be an $F$-rational point of $Y$, and 
$\varphi: X\to Y$  the blowing-up of $Y$ at $Q$,
which is defined over $F$.
By shrinking $U$ if necessary,
we can assume that $Q$ is the generic fiber of 
a closed subscheme $\QQQ\subset \YYY$
that is smooth over $U$.
Let $\beta_U: \XXX\sprime\to \YYY$ be 
the blowing-up of $\YYY$ along $\QQQ$,
which is defined over $U$.
Then the restriction $\beta_\eta: X\sprime\to Y$ of $\beta_U$  to 
the generic fiber $X\sprime$ of $\XXX\sprime\to U$ is isomorphic  to $\varphi$,
that is, there exists an $F$-isomorphism $\tau: X\sprime \isom X$
such that $\beta_\eta=\phi\circ \tau$.
Hence, by Corollaries~\ref{cor:compext}~and~\ref{cor:isomV},  
there exists a non-empty open subset $V\subset U$ such that
the restriction $\beta_V:\XXX\sprime_V\to \YYY_V$
of $\beta_U$  to $\XXX\sprime_V$ coincides with the composite of  the $V$-isomorphism 
$\wt{\tau}_V: \XXX\sprime_V\isom \XXX_V$ and  the extension $\wt{\varphi}_V: \XXX_V\to \YYY_V$
of $\varphi$ over $V$.
\end{example}
\begin{example}\label{example:doublecov}
Let $D$ be a reduced smooth divisor of $Y$
such that every irreducible component $D_i$ of $D$
is defined over $F$.
By shrinking $U$ if necessary,
we can assume that each $D_i$ is the generic fiber of 
a closed subscheme $\DDD_i\subset \YYY$
that is smooth over $U$.
We can also assume that 
these $\DDD_i$ are mutually disjoint.
Then  $\DDD:=\sum \DDD_i$ is smooth over $U$.
\begin{proposition}\label{prop:doubleV}
Let 
$\varphi: X\to Y$ be an $F$-morphism that is a double covering branching exactly along $D$.
Then
there exists an  open subset $V\ne \emptyset$ of $U$
such that 
the extension of $\varphi$ over $V$  is 
a double covering of $\YYY_V$   branching exactly along $\DDD_V$.
\end{proposition}
\begin{proof}
Let $L$ be an invertible sheaf on $Y$ defined by the exact sequence 
$$
0\;\to\; \OOO_Y\;\to\;\varphi_* \OOO_X\;\to\; L\inv\;\to 0.
$$
Then $L$ is defined over $F$, and 
we have  an isomorphism 
$$
\rho\;:\; L^{\otimes 2}\;\isom\;\OOO_Y(D)
$$
on $Y$
that corresponds to the double covering $\varphi$
in the way described in~\cite[Chap.~0]{MR986969}.
There exist a non-empty open subset $V$ of $U$
and an invertible sheaf $\LLL$ on $\YYY_V$ such that
$\LLL\otimes_{\OOO_{\YYY_V}} \OOO_Y=L$.
We consider the invertible sheaves
$$
\MMM:=\fHom_{\OOO_{\YYY_V} }(\LLL^{\otimes 2}, \OOO_{\YYY_V} (\DDD_V))
\quand 
M:=\fHom_{\OOO_Y}(L^{\otimes 2}, \OOO_Y (D)).
$$
on $\YYY_V$ and $Y$, respectively.
Then we have 
$M=\MMM\otimes_{\OOO_{\YYY_V}}\OOO_Y \cong \OOO_Y$.
By~\cite[Proposition 9.3 in Chap.~III]{MR0463157}, 
the restriction homomorphisms
$$
H^0 (\YYY_V, \MMM)\to H^0(Y,  M)
\quand
H^0 (\YYY_V, \MMM\inv )\to H^0(Y,  M\inv)
$$
to the generic fiber $Y$ induce isomorphisms 
$$
H^0 (\YYY_V, \MMM)\otimes_R F \cong H^0(Y, M)
\quand
H^0 (\YYY_V, \MMM\inv)\otimes_R F \cong H^0(Y, M\inv),
$$
where $R:=\Gamma(V, \OOO_V)$.
Hence, by shrinking $V=\Spec R$,
we have elements $f\in H^0 (\YYY_V, \MMM)$ and $g\in H^0 (\YYY_V, \MMM\inv )$
that restrict to $\rho$ and $\rho\inv$, respectively.
Then the composites $f\circ g$ and $g\circ f$,
considered as elements of $H^0(\YYY_V, \MMM\otimes \MMM\inv)= R$,
are mapped to the $1\in H^0(Y, M\otimes M\inv)= F$.
Since $R\inj F$,
we see that $f$ and $g$ are isomorphisms.
Thus $\rho$ extends to an isomorphism 
$$
\wt{\rho}\;:\; \LLL^{\otimes 2}\;\isom\; \OOO_{\YYY_V} (\DDD_V).
$$
By means of  $\wt{\rho}$,  a  double covering $\delta_V: \XXX\sprime_V\to \YYY_V$  that branches exactly along $\DDD_V$
is constructed as a closed subscheme of the line bundle on $\YYY_V$
corresponding to the invertible sheaf $\LLL$.
By construction,
the restriction $\delta_\eta:X\sprime\to Y$ of  $\delta_V$ to 
the generic fiber is isomorphic to $\varphi: X\to Y$.
By Corollaries~\ref{cor:compext}~and~\ref{cor:isomV},
it follows that,  making  $V$ smaller if necessary,
we have a $V$-isomorphism $\XXX_V\sprime\cong\XXX_V$ 
under which $\delta_V$ coincides with the extension $\wt{\varphi}_V$ of $\varphi$ over $V$.
\end{proof}
\end{example}
\begin{example}\label{example:invol}
In this example, we assume $1/2\in R:=\Gamma(U, \OOO_U)$.
Let $\iota: X\isom X$ be an involution defined over $F$,
and $\varphi: X\to Y$  the quotient morphism by the group $\ang{\iota}$.
Suppose that the extension $\wt{\iota}_U:\XXX\isom \XXX$
of $\iota$ over $U$ exists.
Then $\wt{\iota}_U$ is an involution over $U$ by Corollary~\ref{cor:compext}.
Let $q_U: \XXX\to\YYY\sprime$ be the quotient morphism
by the group $\ang{\wt{\iota}_U}$,
which is defined over $U$ by $1/2\in R$.
Then, by Corollaries~\ref{cor:compext},~\ref{cor:isomV} and Lemma~\ref{lem:Ai} below,
we have a non-empty open subset $V\subset U$ and a $V$-isomorphism $\YYY_V\cong \YYY\sprime_V$
 under which 
the extension $\wt{\varphi}_V$ of $\varphi$ over $V$ 
coincides with the restriction $q_V: \XXX_V\to\YYY\sprime_V$
of $q_U$ to $\XXX_V$.
\begin{lemma}\label{lem:Ai}
Let $A$ be an $R$-algebra on which an involution $i$ acts.
Then we have
$A\sp{\ang{i}}\otimes_R F=(A\otimes_R F)\sp{\ang{i}}$,
where $A\sp{\ang{i}}:=\shortset{a\in A}{i(a)=a}$.
\end{lemma}
\begin{proof}
Since $1/2\in R$,
we see that the $R$-module $A$ is the direct-sum of
$A\sp{\ang{i}}=\shortset{(a+i(a))/2}{a\in A}$
and $\shortset{(a-i(a))/2}{a\in A}$.
\end{proof}
\end{example}
\subsection{Shioda-Inose configuration on $\Km (E\sprime\times E)$}
The following result is due to Shioda and Inose~\cite{MR0441982}.
We briefly recall the proof.
\begin{proposition}\label{prop:SIexists}
Let $E\sprime$ and $E$ be elliptic curves defined over 
an algebraically closed field $k$ of characteristic $0$.
Then there exists a Shioda-Inose configuration
$(C, \Theta)$
on the Kummer surface $\Km (E\sprime\times E)$.
\end{proposition}
\begin{proof}
Let $E_{ij}, F_j$ and $G_i$ $(1\le i, j\le 4)$ be the $(-2)$-curves 
in the double Kummer pencil~(Figure~\ref{fig:Km}) on $\Km (E\sprime\times E)$.
We consider  the divisor
\begin{equation}\label{eq:H}
H\;:=\;E_{12}+2F_2+3E_{32}+4 G_3+5E_{31}+6F_{1}+3E_{21}+4E_{41}+2G_{4},
\end{equation}
and let $C$ be the reduced part of $H-E_{12}$:
\begin{equation}\label{eq:C}
C\;:=\;F_2+E_{32}+ G_3+E_{31}+F_{1}+E_{21}+E_{41}+G_{4},
\end{equation}
which is an $ADE$-configuration of $(-2)$-curves of type $\E_8$.
The complete linear system $|H|$ defines an elliptic pencil
$$
\map{\Phi }{ \Km (E\sprime\times E)}{ \P^1}
$$
with a section $G_1$. 
Since $H E_{13}=0$ and $H E_{14}=0$, 
each of  $E_{13}$ and $E_{14}$ is contained in a fiber of $\Phi$.
We put
$t_0\;:=\;\Phi (H)$,
$t_1:=\Phi (E_{13})$ and  
$t_2:=\Phi(E_{14})$.
Note that $t_0\ne t_1\ne t_2\ne t_0$,
because $H$, $E_{13}$ and $E_{14}$ intersect $G_1$ at distinct points.
By~\cite[Theorem~1]{MR0441982}, 
the fibers of $\Phi$ over $t_1$ and $t_2$ are either
(a) of type $\typeI_{b_1}\sp *$ and $\typeI_{b_2}\sp *$ with $b_1+b_2\le 2$, or
(b) of type $\typeI_0\sp *$ and $\typeIV\sp *$.
Hence there exist exactly eight 
$(-2)$-curves $\Theta_1, \dots, \Theta_8$  in $\Phi\inv (t_1)$ and
$\Phi\inv (t_2)$ that appear in the fiber with odd multiplicity.
We denote by  $\Theta$  the sum of $\Theta_1, \dots, \Theta_8$.
Let $\Delta$ be a projective line, 
and $f: \Delta\to \P^1$  the  double covering
that branches exactly at $t_1$ and $t_2$.
Let $\Bl{Y}$ be the normalization of $\Km (E\sprime\times E)\times_{\P^1}\Delta$.
Then $\Bl{Y}\to \Km (E\sprime\times E)$ is 
a finite double covering that  branches exactly along $\Theta$.
Hence $(C, \Theta)$
is a Shioda-Inose configuration.
\end{proof}
\subsection{The $SIK$ diagram}\label{subsec:SIK}
Let $W$ be either a number field, or a Dedekind domain 
with the quotient field $F$ being a number field
such that $1/2\in W$.
\begin{definition}
Let $\EEE\sprime$ and $\EEE$ be smooth proper families of
elliptic curves over $U:=\Spec W$.
We put
$\AAA:=\EEE\sprime\times_U\EEE$.
A diagram
$$
(\SSIIKK)\;\;:\;\; \YYY\;\;\longleftarrow\;\; \Bl{\YYY}\;\;\longrightarrow\;\; \Km(\AAA)\;\;\longleftarrow\;\;\Bl{\AAA}\;\;\longrightarrow\;\; \AAA
$$
of schemes and morphisms over $U$ is called an \emph{$SIK$ diagram}
of $\EEE\sprime$ and $\EEE$ if 
the left-half 
$\YYY\leftarrow\Bl{\YYY}\rightarrow\Km(\AAA)$
is a Shioda-Inose diagram over $U$, and
the right-half
$\Km(\AAA)\leftarrow\Bl{\AAA}\rightarrow\AAA$
is the Kummer diagram of $\EEE\sprime$ and $\EEE$ over $U$.
\end{definition}
%
%
%
%
\begin{proposition}\label{prop:ext}
Let $E\sprime$ and $E$ be elliptic curves 
defined over a number field $L$.

{\rm (1)}
There exist a finite extension $F$ of $L$,
and an $SIK$ diagram 
\begin{equation*}\label{eq:SIKF}
 (\SSIIKK)_F\;\;:\;\; Y \;\;\longleftarrow\;\; \Bl{Y}\;\;\longrightarrow\;\; \Km(A)\;\;\longleftarrow\;\;\Bl{A}\;\;\longrightarrow\;\; A:=(E\sprime\times E)\otimes F
\end{equation*}
of  $E\sprime\otimes F$ and $E\otimes F$  over $F$.

{\rm (2)}
Moreover,  there exist 
a  non-empty open subset $U$ of $\Spec \Z_F[1/2]$,
smooth proper families $\EEE\sprime$ and $\EEE$ of elliptic curves over $U$
with  the generic fibers 
being isomorphic to $E\sprime\otimes F$ and $E\otimes F$, respectively, 
and an $SIK$ diagram
\begin{equation*}\label{eq:SIKU}
(\SSIIKK)_U\;\;:\;\; \YYY\;\;\longleftarrow\;\; \Bl{\YYY}\;\;\longrightarrow\;\; \Km(\AAA)
\;\;\longleftarrow\;\;\Bl{\AAA}\;\;\longrightarrow\;\; \AAA:=\EEE\sprime\times_U \EEE \phantom{\otimes F a}
\end{equation*}
of $\EEE\sprime$ and $\EEE$ over $U$
such that $(\SSIIKK)_U\otimes F$ is equal  to 
the $SIK$ diagram $(\SSIIKK)_F$ over $F$ in {\rm (1)}  above.
\end{proposition}
\begin{proof}
Our argument for the proof of the assertion (1) is similar to~\cite[\S6]{MR0441982}.
We use the notation in the proof of Proposition~\ref{prop:SIexists}.
Let $F$ be a finite extension of $L$
such that every $2$-torsion point $Q_{ij}:=(u\sprime_i, u_j)$ of $A:=(E\sprime \times E)\otimes F$
is rational over $F$.
Then the blowing-up $\Bl{A}\to A$ 
and the involution $\tiliota_A$ of $\Bl{A}$ are defined over $F$.
Therefore the quotient morphism $\Bl{A}\to \Km(A)$ is defined over $F$,
and every irreducible component of the double Kummer pencil on $\Km(A)$
is rational over $F$.
Since the divisor $H$ is defined over $F$, 
the elliptic pencil $\Phi$ on $\Km (A)$
is defined over $F$. 
Moreover, the points  $t_1, t_2\in \P^1$
are $F$-rational.
Replacing $F$ by a finite extension,
we can assume that $\Bl{Y}$ is defined over $F$, 
and that $\Theta_1, \dots, \Theta_8$
 are rational over $F$.
Then the $(-1)$-curves $T_1, \dots, T_8$ on $\Bl{Y}$ 
are rational over $F$, and 
the contraction $\Bl{Y}\to Y$ is defined over $F$.
Moreover,  the image $R_i\in Y$ of  $T_i\subset \Bl{Y}$  is an $F$-rational point of $Y$.
Thus we have obtained an $SIK$ diagram
$(\SSIIKK)_F$ of $E\sprime$ 
and $E$ over $F$,
and the assertion (1) is proved.
Moreover, $(\SSIIKK)_F$
has the following properties:
\begin{rmenumerate}
\item
Each of the center $Q_{ij}$  of the blowing-up $\Bl{A}\to A$
is rational over $F$, and 
each of the center $R_i$  of the blowing-up $Y\leftarrow\Bl{Y}$
is rational over $F$.
\item
Each irreducible component of the double Kummer pencil on $\Km(A)$
is rational over $F$.
In particular, each
irreducible component $C_i$ of the $\E_8$-configuration $C$
is rational over $F$. (See~\eqref{eq:C}.)
\item
Each irreducible component $\Theta_i$  of the branch curve of the double covering
$\Bl{Y} \to \Km (A)$ is rational over $F$.
\end{rmenumerate}
%
%
%
%
We choose a non-empty open subset $U$ of $\Spec \Z_F[1/2]$,
construct smooth proper families  $\EEE\sprime$ and $\EEE$ of elliptic curves over $U$
with  the generic fibers 
being isomorphic to $E\sprime\otimes F$ and $E\otimes F$, respectively, 
and make a diagram 
$(\SSIIKK)_U$
of schemes and morphisms over $U$
such that
each scheme is smooth and projective over $U$,
and that  $(\SSIIKK)_U\otimes  F$ is equal  to 
the $SIK$ diagram $(\SSIIKK)_F$ over $F$.
We will show that,
after  deleting finitely many  closed points from $U$,
the diagram $(\SSIIKK)_U$ becomes an $SIK$ diagram over $U$.
Note that, since $\EEE\sprime$ and $\EEE$ are families of elliptic curves
(that is, with a section over $U$),
the inversion automorphism $\iota_{\AAA}$ of $\AAA$ is defined over $U$.
We can make $U$ so small that the following hold:
\begin{itemize}
\item
Each  $Q_{ij}\in A$ is the generic fiber
of a closed subscheme $\QQQ_{ij}$ of $\AAA$
that is smooth over $U$,
and these  $\QQQ_{ij}$ are mutually disjoint.
Then $\cup \QQQ_{ij}$  is the fixed locus of $\iota_{\AAA}$, and 
 $\Bl{\AAA}\to \AAA$ is the blowing-up along $\cup \QQQ_{ij}$
by Example~\ref{example:blowup}.
\item
The involution 
$\tiliota_A$ of $\Bl{A}$ extends to an involution $(\tiliota_A)\sp\sim_U$
of $\Bl{\AAA}$ over $U$,
which is 
 a lift of $\iota_{\AAA}$ by Corollary~\ref{cor:compext}.
 By Example~\ref{example:invol},
the morphism $\Km(\AAA)\leftarrow \Bl{\AAA}$ is the quotient morphism
by  $\ang{(\tiliota_A)\sp\sim_U}$.
\item
Each $\Theta_i\subset \Km(A)$  
is the generic fiber
of a closed subscheme $\varTheta_i$ of $\Km(\AAA)$
that is smooth over $U$.
By the 
specialization isometry from $\NS(\Km(A))$ to $\NS(\Km(\AAA)\otimes\kappa_{\smallppp})$
for closed points $\ppp$ of $U$, we see  
that these  $\varTheta_i$ are mutually disjoint.
By Example~\ref{example:doublecov},
the morphism $\YYY\to \Km(\AAA)$ is a double covering
branching exactly  along $\varTheta:=\sum \varTheta_i$.
\item
Each irreducible component $C_i$ of $C$
is the generic fiber
of a closed subscheme $\CCC_i$ of $\Km(\AAA)$
that is smooth over $U$.
We put $\CCC:=\sum \CCC_i$.
Considering  the specialization isometry $\NS(\Km(A))\isomet\NS(\Km(\AAA)\otimes\kappa_{\smallppp})$
for closed points $\ppp$ of $U$,
we see that
$\CCC$ is a flat family of $\E_8$-configurations of $(-2)$-curves over $U$,
and that
$\varTheta$ and $\CCC$ are disjoint.
Hence $(\CCC, \varTheta)\otimes\kappa_P$ is 
a Shioda-Inose configuration on $\Km(\AAA)\otimes\kappa_P$
for every point $P$ of $U$.
\item
Each $R_i\in Y$ is the generic fiber
of a closed subscheme $\RRR_i$ of $\YYY$
that is smooth over $U$,
and these  $\RRR_i$ are mutually disjoint.
The morphism  $\Bl{\YYY}\leftarrow\YYY$ is the blowing-up along $\cup \RRR_i$
by Example~\ref{example:blowup}.
\end{itemize}
Hence $(\SSIIKK)_U$ is an $SIK$ diagram over $U$.
\end{proof}
We consider the $SIK$ diagram $(\SSIIKK)_U$ 
over a  non-empty open subset $U\subset \Spec \Z_F[1/2]$,
and 
the $SIK$ diagram $(\SSIIKK)_F=(\SSIIKK)_U\otimes F$ 
over $F$, as in Proposition~\ref{prop:ext}.
(Remark that 
we have changed  the notation  from~\eqref{eq:EEAEEA}~and~\eqref{eq:YYZYYZ}
to $Y:=\YYY\otimes  F$, $E\sprime:=\EEE\sprime\otimes  F$ and $E:=\EEE\otimes F$.)
By the  isomorphisms of Propositions~\ref{prop:PQ}~and~\ref{prop:RS}, we obtain the following:
\begin{proposition}\label{prop:TSIK}
For each $\sigma\in \Emb (F)$, the diagram 
$(\SSIIKK)\otimes\C$ obtained from $(\SSIIKK)\otimes F$
by $\sigma: F\inj \C$ induces an isomorphism of lattices
$T(Y\sp{\sigma})\cong T(A\sp{\sigma})$
that preserves the Hodge structure.
\end{proposition}
We assume the following,
which are equivalent  by Proposition~\ref{prop:NSKmHom} and 
Corollary~\ref{cor:SIoverC}:
(i) $\rank (\Hom (E\sprime, E))=2$.
(ii) $\Km(A)$ is singular.
(iii) $Y$ is singular.
\begin{proposition}\label{prop:SSPsSIK}
We put
$d(Y):=\disc(\NS(Y))$.
There exists a finite set $N$ of prime integers
containing the prime divisors of $2d(Y)$ 
such that the following holds:
\begin{equation}\label{eq:pnotinN}
p\notin N
\;\;\Rightarrow\;\;
\SSps_p(\YYY)=\begin{cases}
\emptyset & \textrm{if $\chi_p(d(Y))=1$, }\\
\pi_F\inv (p) & \textrm{if $\chi_p(d(Y))=-1$. }
\end{cases}
\end{equation}
\end{proposition}
\begin{proof}
By Proposition~\ref{prop:Hom_char0}, there exists an imaginary quadratic field $K$ such that
$K\cong \End(E\sprime)\otimes\Q\cong \End(E)\otimes \Q$.
We denote by $D$ the discriminant of $K$.
We  choose $N$ in such a way that
$N$ contains all the prime divisors of $2 d(Y) D$,
and that 
if $p\notin N$, then $\pi_F\inv(p)\subset U$ holds.
By Propositions~\ref{prop:TSIK} and~\ref{prop:NSA}, we have
\begin{equation*}
\begin{split}
d(Y):=\disc (\NS(Y))&=-\disc (T(Y\sp{\sigma}))=-\disc (T(A\sp{\sigma}))\\
&=\disc (\NS(A))=-\disc (\Hom(E\sprime, E)).
\end{split}
\end{equation*}
By Proposition~\ref{prop:mnsq},
we have $m^2 d(Y)=n^2 D$ for some non-zero integers $m$ and $n$.
We can assume that $\gcd(m, n)=1$.
Then any $p\notin N$ is prime to $mn$, and hence 
$$
p\notin N\;\;\Rightarrow\;\;\chi_p(d(Y))=\chi_p(D)
$$
holds.
Let $p$ be a prime integer not in $N$.
If $\chi_p(d(Y))=1$, then 
$\SSps_p(\YYY)=\emptyset$ by Proposition~\ref{prop:normalK3}.
Suppose that 
$\chi_p(d(Y))=-1$,
and let $\ppp$ be a point of $\pi_F\inv(p)\subset U$.
Since $\chi_p(D)=-1$,  both of
$E\sprime\sbppp :=\EEE\sprime\otimes \kappa\sbppp$ and 
$E\sbppp :=\EEE\otimes\kappa\sbppp$
are  supersingular 
by Proposition~\ref{prop:redss},
and hence $\Km(\AAA)\otimes \kappa\sbppp=\Km(E\sprime\sbppp\times E\sbppp)$
is supersingular by Propositions~\ref{prop:Hom_charp} and~\ref{prop:NSKmHom}.
By Proposition~\ref{prop:NSZY},
we see that $\YYY \otimes\kappa\sbppp$ is also supersingular.
Hence $\pi_F\inv(p)=\SSps_p(\YYY)$ holds.
\end{proof}
From the equality~\eqref{eq:LHom} and Propositions~\ref{prop:LKm},~\ref{prop:LSI}, we obtain 
the following:
\begin{proposition}\label{prop:LSIK}
Let $N$ be a finite set of prime integers with the properties 
given  in Proposition~\ref{prop:SSPsSIK}.
Suppose that $p\notin N$ satisfies $\chi_p(d(Y))=-1$, and 
let $\ppp$ be a point of $\pi_F\inv (p)=\SSps_p(\YYY)$.
Then the diagram 
$(\SSIIKK)_U$ induces an isomorphism of lattices
$$
L(\YYY,\ppp)\wcong (\Hom(E\sprime, E)\isomet \Hom(E\sprime\sbppp, E\sbppp))\sperp[-1],
$$
where
$E\sprime\sbppp :=\EEE\sprime\otimes \kappa\sbppp$ and 
$E\sbppp :=\EEE\otimes\kappa\sbppp$.
\end{proposition}
\subsection{Shioda-Mitani theory}\label{subsec:ShiodaMitani}
In this subsection, we work over $\C$,
and review the Shioda-Mitani theory~\cite{MR0382289} on 
product abelian surfaces.
%
%
%
%
Let $M[a,b,c]$ be a matrix in the set $\QQQ_D$
defined in~\eqref{eq:QQQ},
where $D=b^2-4ac$ is a negative integer.
Let  $\sqrt{D}\in \C$ be in the upper half-plane, and we put
\begin{equation}\label{eq:taus}
\tau\sprime:={(-b+\sqrt{D})}/{(2a)},
\quad
\tau:={(b+\sqrt{D})}/{2}.
\end{equation}
We consider the complex elliptic curves
\begin{equation}\label{eq:anaEs}
\ana E\sprime:=\C/(\Z+\Z\tau\sprime),\quad
\ana E:=\C/(\Z+\Z\tau).
\end{equation}
%
%
%
%
\begin{proposition}[\S3 in \cite{MR0382289}]\label{prop:SM1}
The oriented transcendental lattice $\ori{T}(\ana E\sprime\times \ana E)$ of 
the product abelian surface $\ana E\sprime\times \ana E$
is represented by $M[a,b,c]\in \QQQ_D$.
\end{proposition}
%
%
%
%
Suppose that 
$D$ is a negative fundamental discriminant,
and that  $M[a,b,c]$ is in the set $\QQQ\primitive_D$
defined in~\eqref{eq:QQQprim}.
We put  $K:=\Q(\sqrt{D})\subset \C$.
Then
$$
I_0:=\Z+\Z\tau\sprime
$$
is a fractional ideal of $K$,
and 
$\Z+\Z\tau$ is equal to $\Z_K$.
\begin{proposition}[(4.14) in \cite{MR0382289}]\label{prop:SM2}
Let $J_1$ and $J_2$ be  fractional ideals of $K$.
Then the product abelian surface $\C/J_1\times \C/J_2$
is isomorphic to $\ana E\sprime\times \ana E=\C/I_0\times \C/\ZK$ if and only if 
$[J_1][J_2] =[I_0]$ holds in the ideal class group $\Cl_D$.
\end{proposition}
Recall the definition of  $\Psi: \Cl_D\isom \;\ori{\lats}\primitive_D$
in Proposition~\ref{prop:ClD}.
The image of  $[I_0]\in \Cl_D$ 
  by $\Psi$ is represented by $M[a,b,c]$.
Hence  we obtain the following:
\begin{corollary}\label{cor:SM}
For  fractional ideals $J_1$ and $J_2$ of $K$, we have 
$$
[\,\ori{T}(\C/J_1\times \C/J_2)\,]\;=\;\Psi([J_1][J_2]).
$$
\end{corollary}
\subsection{Proof of Theorem~\ref{thm:almostallp}}\label{subsec:thm:almostallp}
Let $X\to\Spec F$ and $\XXX\to U$ be as in \S\ref{sec:Introduction}.
We choose $\sigma\in \Emb(F)$,
and let $M[a,b,c]\in \QQQ_{d(X)}$ be a matrix representing 
$[\,\ori{T}(X\sp{\sigma})\,]\in \ori{\lats}_{d(X)}$,
where $d(X):=\disc(\NS(X))$.
We define complex elliptic curves
$\ana E\sprime$ and $\ana E$ by~\eqref{eq:taus} and~\eqref{eq:anaEs}.
Then there exist elliptic curves
$E\sprime$ and $E$
defined over a number field $L\subset \C$
such that $E\sprime\otimes\C$ and $E\otimes \C$
are isomorphic to $\ana E\sprime$ and $\ana E$,
respectively.
By replacing $L$ with
a finite extension if necessary,
we have an $SIK$ diagram
$$
\YYY\;\;\longleftarrow\;\; \Bl{\YYY}\;\;\longrightarrow\;\; \Km(\AAA)\;\;\longleftarrow\;\;\Bl{\AAA}
\;\;\longrightarrow\;\; \AAA:=\EEE\sprime \times_{U_L}\EEE
$$
over a  non-empty open subset $U_L$ of $\Spec \Z_L[1/2]$
such that the generic fibers of $\EEE\sprime $ and $\EEE$ are isomorphic to
$E\sprime$ and $E$, respectively.
We put 
$$
A\;:=\;\AAA\otimes L \;=\; E\sprime\times E \quand Y\;:=\;\YYY\otimes L.
$$
Then we see from
Proposition~\ref{prop:SM1}
that $[\,\ori{T} (A\otimes \C)\,]$ is represented by 
the matrix $M[a,b,c]$.
Therefore we have
$[\,\ori{T} (A\otimes \C)\,]=[\,\ori{T}(X\sp{\sigma})\,]$.
On the other hand,
we have $[\,\ori{T}(Y\otimes \C)\,]=[\,\ori{T} (A\otimes \C)\,]$
by Proposition~\ref{prop:TSIK}.
Hence  we obtain
$$
[\,\ori{T}(Y\otimes \C)\,]\,=\,[\,\ori{T}(X\sp{\sigma})\,].
$$
By  the Torelli theorem for  $K3$ surfaces~\cite{MR0284440} or 
the Shioda-Inose theorem (Theorem~\ref{thm:SI}),
the complex $K3$ surfaces  $Y\otimes\C$ and $X\sp{\sigma}$ are isomorphic.
Hence 
$d(Y):=\disc(\NS(Y))$ is equal to $d(X)$.
Moreover,  
there exists a number field $M\subset \C$ containing 
both of $\sigma(F)\subset \C$ and $L\subset \C$
such that $X\otimes M$ and $Y\otimes M$ are isomorphic over $M$.
Then $\XXX\times \Spec\Z_M$ and $\YYY\times \Spec\Z_M$
are isomorphic over the generic point of $\Spec \Z_M$,
and hence 
there exists  a  non-empty open subset $V$ of $\Spec \Z_M$
such that 
$$
\XXX_V:=\XXX\otimes V\quand
\YYY_V:=\YYY\otimes V
$$
are isomorphic over $V$ by Corollary~\ref{cor:isomV}.
Let 
$$
\pi_{M, F}\;:\;\Spec \Z_M\to \Spec \Z_F
\quand
\pi_{M, L}\;:\;\Spec \Z_M\to \Spec \Z_L
$$
be the natural projections.
By deleting finitely many closed points from $V$, we can assume that 
$\pi_{M, F}(V)\subset U$ and $\pi_{M, L}(V)\subset U_L$.
Then we have
$$
\pi_{M, F}\inv (\SSps_p(\XXX))\cap V =\SSps_p(\XXX_V)=\SSps_p(\YYY_V)=\pi_{M, L}\inv (\SSps_p(\YYY))\cap V
$$
for any $p\in \pi_M(V)$.
We choose a finite set $N$ of prime integers
in such a way that the following hold:
\begin{rmenumerate}
\item $N$ contains all the prime divisors of $2 d(X)=2 d(Y)$,
\item  if  $p\notin N$, then $\pi_M\inv (p)\subset V$, and 
hence $\pi_F\inv (p)\subset U$ and $\pi_L\inv (p)\subset U_L$  hold, 
\item  $N$ satisfies the condition~\eqref{eq:pnotinN} for $\YYY$.
\end{rmenumerate}
Then   $N$ satisfies the condition~\eqref{eq:pnotinN1} for $\XXX$.
Hence Theorem~\ref{thm:almostallp} is proved.
\subsection{Proof of Theorem~\ref{thm:toF}(T)}\label{subsec:thm:toFT}
Let $S$ be as in the statement of Theorem~\ref{thm:toF}.
Since $D=\disc(\NS(S))$ is assumed to be 
a fundamental discriminant,
there exists an imaginary quadratic field $K$
with discriminant $D$.
We fix an embedding $K\inj \C$ once and for all.
For a finite extension $L$ of $K$,
we denote by $\Emb(L/K)$ the set of embeddings  of $L$ into $\C$
whose  restrictions to $K$ are  the fixed one.

\par
\smallskip
We recall the theory of complex multiplications.
See~\cite[Chap.~II]{MR1312368}, for example, for detail.
Let $\cloQ\subset\C$
be the algebraic closure of $\Q$ in $\C$, and
let $\ELL (\ZK)$ be the set of $\cloQ $-isomorphism classes $[E]$
of elliptic curves $E$ defined over $\cloQ $
such that $\End(E)\cong\ZK$.
Then $\ELL(\ZK)$ consists of $h$ elements,
where $h$ is the class number $|\Cl_D|$ of $\ZK$.
We denote by
$$
\alpha_1, \dots, \alpha_h\;\;\in\;\; \cloQ \;\;\subset\;\;\C
$$
the $j$-invariants $j(E)$ of the 
isomorphism classes $[E]\in \ELL(\ZK)$, and put
$$
\Phi_D (t) \;:=\; (t-\alpha_1)\cdots (t-\alpha_h).
$$
Then $\Phi_D(t)$ is a polynomial in $\Z[t]$,
which is called the \emph{Hilbert class polynomial}
of $\ZK$.
It is known that  $\Phi_D (t) $ is irreducible in $K[t]$.
%
%
The  field $\HHH:=K(\alpha_1)\subset \C$ is the maximal unramified  abelian  extension of $K$,
which is called the \emph{Hilbert class field} of $K$.
%
%
We define an action of $\Cl_D$ on $\ELL(\ZK)$ by
$$
[I]*[E]\;:=\;[\C/I\inv I_E]
\qquad\textrm{for} \;\; [I]\in \Cl_D
\;\;\textrm{and}\;\; [E]\in \ELL(\ZK),
$$
where $I_E\subset K$ is a fractional ideal  such that 
$E\cong \C/I_E$.
On the other hand,
for  an elliptic curve $E$ defined over $\cloQ$ and  $\gamma\in \Gal (\cloQ /K)$,
we denote by $E\sp{\gamma}$ the elliptic curve
obtained from $E$ by letting $\gamma$ act
on the defining equation for $E$.
Then  $\Gal(\cloQ /K)$ 
acts on $\ELL(\ZK)$ 
by 
$[E]\sp{\gamma}:=[E\sp{\gamma}]$.
The following is the central result in the theory of complex multiplications.
%
%
\begin{theorem}\label{thm:CM}
There exists a homomorphism
$F:\Gal(\cloQ /K)\to\Cl_D$
such that
$[E]\sp{\gamma}=F(\gamma)*[E]$
holds for any $[E]\in \ELL(\ZK)$ and $\gamma\in \Gal(\cloQ /K)$.
Moreover, this homomorphism 
$F$ induces an isomorphism
$\Gal(\HHH/K)\cong \Cl_D$.
\end{theorem}
%
%
We put
$H:=K[t]/(\Phi_D)$,
and denote by $\alpha\in H$ the class of $t\in  K[t]$ modulo 
the ideal $(\Phi_D)$.
Then we have $\Emb (H/K)=\{\sigma_1, \dots, \sigma_h\}$,
where $\sigma_i$ is given by $\sigma_i(\alpha)=\alpha_i$.
Moreover, we have
$\HHH=\sigma_1(H)=\dots=\sigma_h(H)$ in $\C$.
Let $\Ea$ be an elliptic curve defined over $H$
such that
$$
j(\Ea)=\alpha\;\in\; H.
$$
A construction of 
such an elliptic curve is given, 
for example,
in~\cite[\S1 in Chap.~III]{MR817210}.
For each $\sigma_i\in \Emb (H/K)$,
we denote by $\Ea^{\sigma_i}$ the elliptic curve 
defined over $\HHH=\sigma_i(H)\subset \overline{\Q}$
obtained from $\Ea$ by applying $\sigma_i$ to the coefficients of the defining equation.
Then 
we have $j(\Ea^{\sigma_i})=\alpha_i\in\HHH$,
and there exists a unique ideal class $[I_i]\in \Cl_D$  of $K$ such that
$\Ea^{\sigma_i}\cong \C/I_i$.
Moreover,  we have
\begin{equation}\label{eq:ClDIi}
\ELL(\ZK)=\{[\Ea\sp{\sigma_1}], \dots, [\Ea\sp{\sigma_h}]\}
\quand  
\Cl_D=\{[I_1], \dots, [I_h]\}.
\end{equation}
Since $\Gal (\HHH/K)$ is abelian,
there exists a canonical isomorphism $\Gal(H/K)\isom\Gal(\HHH/K)$,
which we will denote by $\gamma\mapsto\wt{\gamma} $.
By Theorem~\ref{thm:CM},
we have an isomorphism $\Gal (H/K)\cong \Cl_D$ 
denoted by $\gamma\mapsto[I_\gamma]$ such that 
\begin{equation}\label{eq:Egamma}
\Ea^{\sigma_i \circ \gamma}= (\Ea^\gamma)\sp{\sigma_i}=(\Ea\sp{\sigma_i})\sp{\wt{\gamma} } \wcong \C/{I_\gamma}\inv I_i
\end{equation}
holds for any $i=1, \dots, h$ and any $\gamma\in \Gal(H/K)$.
\par
There exist a finite extension
$F$ of $H$
and a non-empty open subset $U$ of $\Spec \Z_F[1/2]$
such that,
for each $\gamma\in \Gal(H/K)$, there exist
 smooth proper families of elliptic curves $\EEE\sb{\alpha}^\gamma$ and $\EEE\sb{\alpha}$
 over $U$ whose  generic fibers 
are isomorphic to $\Ea^\gamma\otimes F$ and $\Ea\otimes F$, respectively, and 
 an $SIK$ diagram
$$
(\SSIIKK)^{\gamma}\;\;\;:\;\;\;
\YYYgamma\;\;\longleftarrow\;\; 
\Bl{\YYY}\sp\gamma\;\;\longrightarrow\;\; 
\Km(\AAA^\gamma)\;\;\longleftarrow\;\;
\Bl{\AAA}^\gamma\;\;\longrightarrow\;\; 
\AAA^\gamma:=\EEE\sb{\alpha}^\gamma\times_{U}\EEE\sb{\alpha}
$$
of $\EEE\sb{\alpha}^\gamma$ and $\EEE\sb{\alpha}$ over $U$.
We then put
$Y^\gamma:=\YYYgamma\otimes F$ and $A^\gamma:=\AAA^\gamma\otimes F$.
Let $\sigma$ be an element of $\Emb (F/K)$.
If the restriction of $\sigma$ to $H$ is equal to $\sigma_i$, then
we have the following equalities in $\ori{\lats}\primitive_D$:
$$
\renewcommand{\arraystretch}{1.4}
\begin{array}{lcll}
[\,\ori{T} ((Y\sp\gamma)\sp{\sigma})\,] &=&  [\,\ori{T} ((A\sp\gamma)\sp{\sigma})\,] &\textrm{by Proposition~\ref{prop:TSIK}}\\
&=& [\,\ori{T} (\Ea\sp{\sigma_i\circ\gamma}\times \Ea\sp{\sigma_i})\,] &\\
&=&[\,\ori{T} (\C/({I_\gamma\inv} I_i) \times \C/I_i)\,] & \textrm{by~\eqref{eq:Egamma}}\\
& =&\Psi([I_\gamma]\inv [I_i]^2) &\textrm{by Corollary~\ref{cor:SM}}.
\end{array}
$$
Note that the restriction map
$\Emb(F/K)\to \Emb(H/K)$
is surjective.
Therefore, by Proposition~\ref{prop:ClD} and the equalities~\eqref{eq:ClDIi},
we see that the subset 
$$
\set{ [\,\ori{T} ((Y\sp\gamma)\sp{\sigma})\,] }{\sigma\in \Emb(F/K)}
\;\;\;=\;\;\;
\set{\Psi([I_\gamma]\inv [I_i]^2)}{i=1, \dots, h}
$$
of $\ori{\lats}\primitive_D$
coincides with the lifted genus
that contains $\Psi([I_\gamma]\inv )$.
Since the homomorphism  $\gamma\mapsto [I_\gamma]$ from $\Gal(H/K)$ to $\Cl_D$
is an isomorphism,
we have a unique element $\gamma(S)\in \Gal(H/K)$ such that
$\Psi([I_{\gamma(S)}]\inv )$ is equal to $[\,\ori{T}(S)\,]$.
We put  
$$
\XXX\;:=\; \YYY^{\,\gamma(S)}\quand X\;:=\; Y^{\gamma(S)}.
$$
Then $X$ has the property required in Theorem~\ref{thm:toF}(T).
%
%
\subsection{Proof of Theorem~\ref{thm:toF}(L)}\label{subsec:thm:toFL}
We continue to use the notation fixed in the previous subsection.
We consider $\Ea$ as  being defined over $F$.
Replacing $F$ by a finite extension if necessary,
we can assume that $F$ is Galois over $\Q$, and that 
\begin{equation}\label{eq:EndF}
\End _F(\Ea)=\End (\Ea)
\end{equation}
holds so that $\Lie: \End(\Ea)\to F$ is defined. 
Since $j(\Ea)=\alpha$ is a root of $\Phi_D$ and 
$F$ contains $K$, we have the $\Lie$-normalized isomorphism
\begin{equation}\label{eq:normalized}
\End (\Ea)\wcong \Z_K.
\end{equation}
Making the base space $U$ 
of the $SIK$ diagram $(\SSIIKK)^{\gamma(S)}$  
smaller if necessary,
we can assume the following:
\begin{enumerate}
\renewcommand{\theenumi}{\roman{enumi}}
\item $U=\pi_F\inv(\pi_F(U))$, and $p\notdiv 2D$ for any $p\in \pi_F(U)$,
\item if $p\in \pi_F(U)$, then  $\Phi_D(t)\bmod p$ has no multiple roots in $\cloF_p$, and
\item for $p\in \pi_F(U)$,  we have the following equivalence:
$$
\chi_p(D)=-1\;\;\;\Leftrightarrow\;\;\;
\SSps_p(\XXX)\ne\emptyset \;\;\;\Leftrightarrow\;\;\;
\SSps_p(\XXX)=\pi_F\inv (p).
$$
\end{enumerate}
Let $p$ be a prime integer in $\pi_F(U)$
such that
$\chi_p(D)=-1$,
so that $\SSps_p(\XXX)=\pi_F\inv (p)$.
We show that,
under the assumption that $D$ is odd, 
the set of isomorphism classes of supersingular reduction lattices 
$\shortset{[L(\XXX,\ppp)]}{\ppp\in \pi_F\inv (p)}$
coincides with a genus.
\par
\smallskip
Let $B$ denote the quaternion algebra over $\Q$
that ramifies exactly at $p$ and $\infty$.
%
We consider pairs $(R,Z)$ of
a $\Z$-algebra $R$ and a subalgebra $Z\subset R$
such that $R$ is isomorphic to a maximal order of $B$,
and that $Z$ is isomorphic to $\ZK$.
We say that two such pairs $(R,Z)$ and $(R\sprime,Z\sprime)$
are isomorphic if 
there exists an isomorphism $\varphi:R\isom R\sprime$
satisfying  $\varphi(Z)=Z\sprime$.
We denote by $\RRR$
the set of isomorphism classes $[R,Z]$ of these pairs.
Next we consider pairs $(R,\rho)$ of
a $\Z$-algebra $R$ isomorphic to a maximal order of $B$
and an embedding $\rho: \ZK\inj R$ as a $\Z$-subalgebra.
We say that two such pairs $(R,\rho)$ and $(R\sprime,\rho\sprime)$
are isomorphic if 
there exists an isomorphism $\varphi:R\isom R\sprime$
satisfying  $\varphi\circ\rho=\rho\sprime$.
We denote by $\ori{\RRR}$
the set of isomorphism classes $[R,\rho]$ of these pairs.
For an embedding $\rho:\ZK\inj R$,
we denote by $\bar\rho$ the composite of 
the non-trivial automorphism of $\ZK$ and $\rho$.
The natural map 
$$
\map{\Pi_{\RRR}}{\ori{\RRR}}{ \RRR}
$$
given by $[R, \rho]\mapsto [R, \rho(\ZK)]$
is surjective,
and its fiber consists either 
of two elements  $[R, \rho]$ and $[R,\bar\rho]$,
or of a single element $[R, \rho]=[R,\bar\rho]$.
\par
\smallskip
Let $\ppp$ be a point of $\pi_F\inv (p)$.
We denote by $F\sbsqppp$ the completion of $F$ at $\ppp$,
and put 
$$
E\sbsqppp:=\EEE\sb{\alpha}\otimes F\sbsqppp
\quand
E\sbppp:=\EEE\sb{\alpha}\otimes\kappa\sbppp.
$$
Then we have  canonical  isomorphisms
\begin{equation}\label{eq:Endsbsqppp}
\End_{F\sbsqppp}(E\sbsqppp)\wcong  \End(E\sbsqppp)\wcong \End(\Ea)
\end{equation}
by the assumption~\eqref{eq:EndF}, and hence $\Lie: \End(E\sbsqppp)\to F\sbsqppp$ is defined.
We put
$$
R\sbppp\;:=\;\End(E\sbppp),
$$
which is isomorphic to a maximal order of $B$ by Proposition~\ref{prop:Endss1}, 
and denote by
$$
\mapisomet{\rho\sbppp}{ \End(E\sbsqppp)}{R\sbppp}
$$
the specialization isometry.
Using the  isomorphisms~\eqref{eq:Endsbsqppp}
and the $\Lie$-normalized isomorphism~\eqref{eq:normalized},
we obtain an element
$[R\sbppp,\rho\sbppp]$ of $\ori{\RRR}$.
We denote by 
\begin{equation*}\label{eq:SStoRRR}
\map{\ori{r} }{\SSps_p(\XXX)}{ \ori{\RRR}}
\end{equation*}
the map given by $\ppp\mapsto [R\sbppp,\rho\sbppp]$.
\begin{lemma}\label{lemma:surj}
The map $\ori{r}$ is surjective.
\end{lemma}
\begin{proof}
First we show that
the map $r:=\Pi_{\RRR}\circ \ori{r}$ from $\SSps_p(\XXX)$ to $\RRR$
is surjective.
Let $[R, Z]$ be an element of $\RRR$.
By Proposition~\ref{prop:Endss2},
there exists a supersingular elliptic curve $C_0$ 
in characteristic $p$ with an isomorphism 
$\psi: \End(C_0)\isom R$.
Let 
$\alpha_0\in \End(C_0)$
be an element 
such that the subalgebra $\Z+\Z\,\alpha_0$ corresponds to $Z\subset R$ by $\psi$.
By Proposition~\ref{prop:lift},
there exists a  lift $(C, \alpha)$ of $(C_0, \alpha_0)$,
where $C$ is an elliptic curve defined over a finite extension of $\Q_p$.
Since $\Z+\Z\,\alpha\subseteq \End(C)$ is isomorphic to $\ZK$, we have
$\End(C)\cong\ZK$, and hence the $j$-invariant of $C$ is a root of 
the Hilbert class polynomial $\Phi_D$ in $\cloQ _p$.
Since  the set of roots of 
$\Phi_D$ in $\cloQ _p$ is in one-to-one correspondence with
$\pi_H\inv(p)$ by  the assumption (ii) on $U$,
and $U$ contains  $\pi_F\inv (p)$ by  the assumption (i) on $U$, 
there exists $\ppp\in \pi_F\inv (p)\subset U$
such that 
$$
j(E\sbsqppp)=j(C).
$$
 By applying Proposition~\ref{prop:indep2} with $g=\id$, we have
$r(\ppp)=[R, Z]$.
To prove that $\ori{r}$ is surjective,
therefore, it is enough to show that,
for each $\ppp\in \pi_F\inv(p)$,
there exists $\ppp\sprime \in \pi_F\inv(p)$ such that
$[R\sb{\smallppp\sprime}, \rho\sb{\smallppp\sprime}]=[R\sbppp, \overline{\rho\sbppp}]$
holds in $\ori{\RRR}$.
We choose an element $g\in \Gal (F/\Q)$ such that 
the restriction of $g$ to $K$  is the non-trivial element of $\Gal(K/\Q)$,
and let $\ppp\sprime$ be the image of $\ppp$ by the action of $g$ on $\pi_F\inv (p)$.
Consider the diagram 
$$
\renewcommand{\arraystretch}{1.4}
\begin{array}{ccccccc}
F\sbsqppp &&\hskip -.3cm\mapleftsp{\Lie}\hskip -.3cm&&\End(E\sbsqppp) & \maprightsp{\rho\sbppp} & \End(E\sbppp) \\
			          & \NWinj \hskip .5cm &                 &\hskip .8cm\NEisom \;\; \smallmath{\lambda}\hskip -.5cm&                            && \\
\smallmath{f_g} \Big\downarrow\wr   & \hskip .5cm F&\hskip -.3cm\mapleftsp{\Lie}\hskip -.3cm & \End(\Ea) &   \smallmath{e_g}\Big\downarrow\wr& &\smallmath{E_g} \Big\downarrow\wr \\
			        & \SWinj  \hskip .7cm&                 &\hskip .8cm\SEisom \;\; \smallmath{\lambda\sprime} \hskip -.5cm&                            && \\
F\sb{[\smallppp\sprime]}&&\hskip -.3cm\mapleftsp{\Lie}\hskip -.3cm &&\End(E\sb{[\smallppp\sprime]}) &
\maprightsp{\rho\sb{\smallppp\sprime}} & \End(E\sb{\smallppp\sprime}),
\end{array}
$$
where $\lambda$ and $\lambda\sprime$ are the canonical isomorphisms~\eqref{eq:Endsbsqppp},
and the  vertical isomorphisms $f_g$, $e_g$ and $E_g$ are given by the action of $g$.
Then we have
$e_g\circ \lambda=\conj{\lambda\sprime}$,
where $\conj{\lambda\sprime}$ is the composite of the nontrivial automorphism of $\End(\Ea)\cong\ZK$ and $\lambda\sprime$.
By Proposition~\ref{prop:indep2},
we have $E_g\circ \rho\sbppp=\rho\sb{\smallppp\sprime}\circ e_g$,  and hence 
$[R\sb{\smallppp\sprime}, \rho\sb{\smallppp\sprime}]=[R\sbppp, \overline{\rho\sbppp}]$. 
\end{proof}
Suppose that the ideal class $[I_{\gamma(S)}]\in \Cl_D$ is represented by an ideal $J\subset \ZK$.
We can regard $J$ as an ideal of $\End (\Ea)$ by the $\Lie$-normalized isomorphism~\eqref{eq:normalized}.
By~\cite[Corollary 7.17]{MR1028322},
we can choose  $J\subset \ZK$  in such a way that
$$
d_J:=\deg\phi^J=[\End(\Ea):J]
$$
is prime to $D$. (See Remark~\ref{rem:assumps}~(2).)
For any $\sigma_i\in \Emb (H/K)$,
we have the following isomorphisms 
of complex elliptic curves:
$$
\renewcommand{\arraystretch}{1.2}
\begin{array}{lcll}
(\Ea^{\gamma(S)})\sp{\sigma_i}&\wcong & \C/{I_{\gamma(S)} }\inv I_i &\textrm{by~\eqref{eq:Egamma}}\\
&\wcong  & \C/J\inv I_i &\textrm{by~$[J]=[I_{\gamma(S)}]$}\\
&\wcong& (\Ea\sp{\sigma_i})^J&\textrm{by $\Ea\sp{\sigma_i}\cong \C/I_i$ and~Proposition~\ref{prop:anaEJ}}\\
&\wcong& (\Ea\sp J)^{\sigma_i} &\textrm{since the construction of $E\to E^J$ is algebraic.}
 \end{array}
$$
Hence we have
\begin{equation}\label{eq:EagammaSFbar}
\Ea^{\gamma(S)}\otimes \overline{F} \wcong \Ea^J \otimes \overline{F}.
\end{equation}
We then consider $J$ as 
an ideal of 
$\End(E\sbsqppp)$ by the canonical isomorphisms~\eqref{eq:Endsbsqppp}, and consider 
 the left-ideal
$R\sbppp\rho\sbppp(J)$
of $R\sbppp$ generated by $\rho\sbppp(J)$.
From the isomorphism~\eqref{eq:EagammaSFbar},
we obtain  an isomorphism
$$
\EEE_\alpha\sp{\gamma(S)} \otimes\overline{F}\sbsqppp \wcong {E\sbsqppp}^J \otimes\overline{F}\sbsqppp.
$$
Then, by~Proposition~\ref{prop:isogs},
we have
$$
\EEE_\alpha\sp{\gamma(S)}\otimes\bar{\kappa}\sbppp \wcong  {E\sbppp}^{R\sbppp\rho\sbppp(J)}\otimes\bar{\kappa}\sbppp.
$$
Therefore we have the following equalities in the set $\lats_{p^2d_J^2D}$:
$$
\renewcommand{\arraystretch}{1.2}
\begin{array}{lcll}
&& [\;L(\XXX,\ppp) [-d_J]\;]&\\
&= & 
[\;(\Hom({E\sbsqppp}^J, {E\sbsqppp}) \isomet\Hom({E\sbppp}^{R\sbppp\rho\sbppp(J)},{E\sbppp}))\sperp\;[d_J]\;] 
&\textrm{by Proposition~\ref{prop:LSIK}}\\
&= &[\;(J\isomet R\sbppp\rho\sbppp(J))\sperp\;] &\textrm{by Proposition~\ref{prop:JRJsperp}}.
 \end{array}
$$
By the surjectivity of the map $\ori{r}$, 
we complete the proof   of  Theorem~\ref{thm:toF}(L) by the following proposition,
which will be proved in the next section.
\begin{proposition}\label{prop:Jgenus}
Let $J$ be an ideal of $\ZK$.
Suppose that $D$ is odd, and that $d_J=N(J)=[\ZK:J]$ is prime to $D$.
Then the set
$$
\set{[(J\isomet R \rho (J))\sperp]}{[R, \rho]\in\ori{\RRR}}
$$
coincides with  a genus in $\lats_{p^2d_J^2D}$.
\end{proposition}
\begin{remark}\label{rem:assumps}
(1) 
We make use of the assumption  that $D$ is odd in  Theorem~\ref{thm:toF}(L)
only in the proof of Proposition~\ref{prop:Jgenus}.
(2) The condition  $\gcd(N(J), D)=1$ is assumed  only 
in order to simplify the proof of Proposition~\ref{prop:Jgenus}.
\end{remark}
\section{The maximal orders of a quaternion algebra}\label{sec:Endss}
Let $K$, $D$, $p$, $B$ and $\ori{\RRR}$ be as in the previous section.
We assume that $D$ is odd.
We describe the set $\ori{\RRR}$ following  Dorman~\cite{MR1024555},
and prove Proposition~\ref{prop:Jgenus}.
\subsection{Dorman's description of $\ori{\RRR}$}
Note that  $D$ is a square-free negative integer 
satisfying $D\equiv 1\,\bmod 4$.
We choose  a prime integer $q$ that satisfies
\begin{equation}\label{eq:chi}
\chi_l(-pq)=1\quad \textrm{for all prime divisors  $l$ of $D$}.
\end{equation}
Then  the $\Q$-algebra
\begin{equation*}\label{eq:matB}
B\;:=\;\set{[\alpha, \beta]}{\alpha, \beta\in K},
\quad 
\textrm{where}\quad 
[\alpha, \beta]:=\left(\hskip -2pt 
\begin{array}{cc}
\alpha & \beta\\
-pq\conj{\beta} & \conj{\alpha}
\end{array} \hskip -2pt 
\right), 
\end{equation*}
is a quaternion algebra 
that ramifies exactly 
at $p$ and $\infty$.
The canonical involution of $B$ is  given by
$[\alpha, \beta]\sp *=[\conj{\alpha}, -\beta]$.
Hence  the bilinear  form~\eqref{eq:bilB} on $B$ is given by
\begin{equation}\label{eq:bilmatB}
([\alpha, \beta], [\alpha\sprime, \beta\sprime])\;=\;
\Tr_{K/\Q} (\alpha\conj{\alpha\sprime})+pq \Tr_{K/\Q} (\beta\conj{\beta\sprime}).
\end{equation}
Note that we have
\begin{equation}\label{eq:abc}
[\gamma, 0][\alpha, \beta]=[\gamma\alpha,\gamma\beta]
\quand
[\alpha, \beta][\gamma, 0]=[\gamma\alpha,\conj{\gamma}\beta].
\end{equation}
For simplicity, we use the following notation:
$$
[S, T]\;:=\; \set{[\alpha, \beta]\in B}{\alpha\in S, \, \beta\in T}
\qquad
\textrm{for subsets $S$ and  $T$ of $K$.}
$$
For $u\in B\sptimes$, we denote by $\theta_u: B\isom B$
the inner automorphism $\theta_u(x):=uxu\inv$.
We have a natural embedding $\iota: K\inj B$
given by $\iota(\alpha):=[\alpha, 0]$.
By the Skolem-Noether theorem~(\cite[\S10 in Chap.~8]{MR0098114}), 
we see that, if
$\iota\sprime : K\inj B$
is another embedding as a $\Q$-algebra, 
then there exists $u\in B\sptimes$ such that
$\theta_u\circ \iota=\iota\sprime$ holds.
On the other hand, 
we have $\theta_u\circ \iota=\iota$ if and only if $u\in [K\sptimes, 0]$.
Hence 
we have a canonical identification
$$
\ori{\RRR}\wcong [K\sptimes, 0]\backslash \setR,
$$
where $\setR$ is the set of maximal orders $R$ of $B$ such that $R\cap[K,0]=[\ZK,0]$ holds.
We will examine the set $\setR$.
\par
\smallskip
For a $\Z$-submodule $\lat\subset B$ of rank $4$, we put
$$
\NB(\lat):=[[\ZK,\ZK]: n\lat]/n^4,
$$
where $n$ is a non-zero integer such that $n\lat\subset[\ZK,\ZK]$.
An order  $R$ of $B$ is maximal if and only if
$R$ is  of discriminant $p^2$ as a lattice.
Since the discriminant of  $[\ZK,\ZK]$ is   $p^2q^2 |D|^2$,
we obtain the following:
\begin{lemma}\label{lem:NZ}
An order  $R$ of $B$ is maximal if and only if
$\NB(R)=1/q|D|$.
\end{lemma}
We denote by
$\shortmap{\pr_2}{B}{K}$
the projection given by  $\pr_2([\alpha, \beta]):=\beta$.
\begin{lemma}\label{lem:NM}
Let $R$ be an element of $\setR$.
Then $M_R:=\pr_2(R)$ is a fractional ideal of $K$ with 
$N (M_R)=1/q|D|$.
\end{lemma}
\begin{proof}
It is obvious that $M_R\subset K$ is a finitely generated $\ZK$-module
by the formula~\eqref{eq:abc}.
Since $[K,0]\cap R=[\ZK, 0]$, we have
$N(M_R)=\NB(R)=1/q|D|$ by Lemma~\ref{lem:NZ}.
\end{proof}
From the condition~\eqref{eq:chi} on $q$, $\chi_p(D)=-1$ and $D\equiv 1\bmod 4$, 
we deduce that  $q$ splits completely in $K$.
We choose 
an ideal $Q\subset \ZK$ such that 
$(q)=Q\conj{Q}$.
We also denote by $\DDD$ the principal ideal 
$(\sqD)\subset \ZK$.
Let $R$ be an element of $ \setR$.
By Lemma~\ref{lem:NM},
the fractional ideal
$$
I_R:=\DDD Q M_R\qquad (M_R:=\pr_2(R))
$$
satisfies
$N(I_R)=1$.
Since $[K,0]\cap R=[\ZK, 0]$, we can define a map
$$
\map{f_R}{M_R}{K/\ZK}
$$
by $f_R(\beta):=\alpha+\ZK$ for $[\alpha, \beta]\in R$.
By the formula~\eqref{eq:abc}, we see that
$f_R$ is a homomorphism of $\ZK$-modules,
and we have
$f_R(\gamma\beta)=f_R(\conj{\gamma}\beta)$
for any $\gamma\in \ZK$ and $\beta\in M_R$.
Therefore $f_R(\sqD\beta)=\sqD f_R(\beta)=0$ holds for any $\beta\in M_R$.
Thus $f_R$ induces a homomorphism
$$
\map{\wt{f}_R}{M_R/\DDD M_R}{\DDD\inv/\ZK}
$$
of torsion $\ZK$-modules.
\begin{lemma}\label{lem:wtfisom}
The homomorphism $\wt{f}_R$ is an isomorphism.
\end{lemma}
\begin{proof}
Since $|M_R/\DDD M_R|=|\DDD\inv/\ZK|=|D|$,
it is enough to show that $\wt{f}_R$ is injective.
Let $\FFF$ be the fractional ideal such that
$\Ker (f_R)=\FFF\DDD M_R=\FFF Q\inv I_R$.
Suppose that
$\beta, \beta\sprime\in \Ker (f_R)$.
Then $[0, \beta]\in R$ and $[0, \beta\sprime]\in R$ hold,
and hence
$[0, \beta]\cdot [0, \beta\sprime]=
[-pq\beta\conj{\beta\sprime}, 0]$ is also in $R$.
From $[K,0]\cap R=[\ZK, 0]$, we have
$-pq\beta\conj{\beta\sprime}\in \ZK$.
Since $N(I_R)=1$, we have   
$$
pq (\FFF Q\inv I_R)(\overline{\FFF Q\inv I_R})=p \FFF\overline{\FFF}\;\subseteq \;\ZK.
$$
Since $\gcd(p,D)=1$ and $\ZK\subseteq \FFF\subseteq \DDD\inv$,
we have  $\FFF=\ZK$.
\end{proof}
Since $\wt{f}_R$ is an isomorphism, there exists
a unique element
$$
\mu_R+\DDD M_R=\mu_R+Q\inv I_R \;\in \;M_R/\DDD M_R
$$
such that
$\wt{f}_R(\mu_R+Q\inv I_R)=(1/\sqD )+\ZK$.
\begin{lemma}\label{lem:ft}
For any $\beta\in M_R$, 
we have
$$
f_R(\beta)= pq \sqD \conj{\mu_R}\beta +\ZK=pq \sqD \mu_R\conj{\beta}+\ZK.
$$
\end{lemma}
\begin{proof}
For  $\beta\in M_R$,
we have $[0,\sqD\beta]\in R$.
Since
$[1/\sqD, \mu_R]\in R$, we have
\begin{eqnarray*}
[0,\sqD\beta]\cdot[1/\sqD, \mu_R]&=\;\;[\,-pq\sqD \conj{\mu_R}\beta, \,-\beta\,] \;\;&\in\;\; R\quand \\[0pt]
[1/\sqD, \mu_R]\cdot[0,\sqD\beta]&=\;\;[\,pq\sqD \mu_R\conj \beta,\, \beta\,] \;\;\;\;\;\;\;&\in\;\; R,
\end{eqnarray*}
from which the desired description of $f_R$ follows.
\end{proof}
By the definition of $\mu_R$, we have
$f_R(\mu_R)=pq \sqD \conj{\mu_R}\mu_R+\ZK=(1/\sqD)+\ZK$.
Therefore we have
$$
pqD|\mu_R|^2-1\;\in\; \DDD\cap \Q=D\Z,
$$
where the second equality follows from the assumption that $D$ is odd.
\begin{lemma}\label{lem:normZ}
Let $I$ be a fractional ideal with  $N(I)=1$,
and let $x, x\sprime\in \DDD\inv Q\inv I$
satisfy  $x\sprime -x \in Q\inv I$.
Then we have $qD|x|^2\in \Z$ and 
$qD |x\sprime|^2 \equiv qD |x|^2\,\bmod D$.
\end{lemma}
\begin{proof}
Since $I\conj{I}=\ZK$, we have 
$$
qD|x|^2\;\in\; q D (\DDD^{-2} Q\inv  I\conj{Q\inv I})\cap \Q=\ZK\cap \Q=\Z.
$$
We put $x\sprime=x+y$ with $y\in Q\inv I$.
Then we have
$q|y|^2\in \ZK\cap \Q=\Z$.
Since $D$ is odd, we have
$\DDD\inv\cap \Q=\Z$, and hence
$q (x \conj{y}+y \conj{x})\in  \DDD\inv \cap \Q=\Z$ holds.
\end{proof}
We define $\TTT$ to be the set of all pairs
$(I, \mu+Q\inv I)$,
where 
$I$ is a fractional ideal of $K$ such that $N(I)=1$, and 
$\mu+Q\inv I$ is an element of $\DDD\inv Q\inv I/Q\inv I$ such that
$pqD|\mu|^2\equiv 1\,\bmod D$.
Then  we have   a map $\shortmap{\tau}{\setR}{\TTT}$ given  by 
$$
\tau(R)\;:=\; (I_R, \mu_R+Q\inv I_R )\;\in\; \TTT.
$$
\begin{proposition}\label{prop:tau}
The  map $\tau$ is a bijection.
\end{proposition}
\begin{proof}
The maximal order $R$ is uniquely recovered  from  $(I_R, \mu_R+Q\inv I_R )$ by
$$
R\;=\;\set{[\alpha, \beta]}{\beta\in \DDD\inv Q\inv I_R, \; \alpha\equiv  pq \sqD \conj{\mu_R}\beta \,\bmod \ZK}.
$$
Hence $\tau$ is injective.
Let an element 
$t:=(I, \mu+Q\inv I)$
of $\TTT$ be given.
We put
$M_t:=\DDD\inv Q\inv I$,
and define $f_t: M_t\to \DDD\inv /\ZK$ by
$$
f_t(\beta) \;:=\; pq\sqD\conj{\mu}\beta +\ZK.
$$
Note that the definition of $f_t$
does not depend on the choice of 
the representative $\mu$ of $\mu+Q\inv I$.
Since $M_t\conj{M_t}=(1/Dq)$,
we see that 
${\mu}\conj{\beta}-\conj{\mu}\beta$ is contained in 
$\DDD(1/Dq)=\DDD\inv(1/q)$ for any $\beta\in M_t$.
(Note that $\gamma-\conj{\gamma}\in \DDD$ for any $\gamma\in \ZK$.) 
Therefore we have
\begin{equation}\label{eq:secondeq}
f_t(\beta) \;=\; pq\sqD\mu\conj{\beta} +\ZK
\end{equation}
for any $\beta\in M_t$.
We put
$$
R_t\;:=\;\set{[\alpha, \beta]}{\beta\in M_t, \;\;\alpha\equiv f_t(\beta)\,\bmod \ZK}.
$$
We prove that $\tau$ is surjective
by showing that 
$R_t\in \setR$.
It is obvious that $R_t$ is a  $\Z$-module of rank $4$
satisfying  $R_t\cap[K, 0]=[\ZK, 0]$.
We  show that $R_t$ is closed under the product.
Since $f_t$ is a homomorphism of $\ZK$-modules,
we have
$[\ZK, 0]R_t=R_t$.
By the formula~\eqref{eq:secondeq}, we have
$R_t [\ZK, 0]=R_t$.
Hence it is enough to prove that
\begin{equation*}\label{eq:prod}
[\,pq\sqD\mu\conj{\beta},\;\beta\,]\cdot [\,pq\sqD\conj{\mu}\beta\sprime,\;\beta\sprime\,]
=[\,p^2q^2 D |\mu|^2\conj{\beta}\beta\sprime-pq\beta\conj{\beta\sprime}, \;\; pq \sqD \mu (\conj{\beta}\beta\sprime-\beta\conj{\beta\sprime})\,]
\end{equation*}
is in $R_t$ for any $\beta, \beta\sprime\in M_t$.
Because
$$
pq \beta\conj{\beta\sprime}\in pqM_t\conj{M}_t=(p/D)\quand
1-pqD|\mu|^2\equiv 0\,\bmod D,
$$
we have a congruence
$pq \beta\conj{\beta\sprime}\equiv p^2q^2 D |\mu|^2 \beta\conj{\beta\sprime}\,\bmod\ZK$.
Hence 
$$
f_t(pq \sqD \mu (\conj\beta\beta\sprime-\beta\conj{\beta\sprime}))
\equiv 
p^2q^2 D |\mu|^2\conj{\beta}\beta\sprime-pq\beta\conj{\beta\sprime}
\,\bmod\ZK.
$$
Therefore $R_tR_t=R_t$ is proved,
and hence $R_t$ is an order.
Because
$N(M_t)=1/q|D|$, we see that  $R_t$ is maximal by Lemma~\ref{lem:NZ}.
Hence $R_t\in \setR$.
\end{proof}
We make $[K\sptimes, 0]$ act on the set $\TTT$ by
$$
u\cdot(I, \mu+Q\inv I):=(I u\conj{u}\inv, \mu u\conj{u}\inv+Q\inv I u\conj{u}\inv).
$$
Then $\ori{\RRR}\cong  [K\sptimes, 0]\backslash \setR $
is canonically identified with $[K\sptimes, 0]\backslash \TTT$.
\par
\smallskip
Let $\III_1\subset \III_D$ denote the group of
fractional ideals $I$ with $N(I)=1$.
We put
$$
\PPP_1:=\III_1\cap\PPP_D
\quand
\CCC_1:=\III_1/\PPP_1.
$$
Then $\CCC_1$ is a subgroup of the ideal class group $\Cl_D=\III_D/\PPP_D$.
Since the homomorphism $\III_D\to\III_1$ given by $I\mapsto I\conj{I}\inv$ is surjective
and  $[I][\conj{I}]\inv$ is equal to $[I]^2$ in $\Cl_D$, we see that 
the subgroup $\CCC_1$ of $\Cl_D$ is equal to $\Cl_D^2$.
\begin{lemma}\label{lem:surjR}
The map $R\mapsto [I_R]$
from $ \setR$ to $\CCC_1=\Cl_D^2$
 is surjective.
\end{lemma}
\begin{proof}
We will show that,
for each  $I\in \III_1$,
there exists $\mu\in \DDD\inv Q\inv I$ such that
$pqD|\mu|^2\equiv 1\,\bmod D$ holds.
Since $I\conj{I}=\ZK$, there exists an ideal
$A\subset\ZK$ such that $A+\conj{A}=\ZK$ and $I=A\conj{A}\inv$.
Since $A$ and $\conj{A}$ have no common prime divisors,
the norm $n:=N(A)$ is prime to $D$.
Hence  there exists an integer $m$
such that $nm\equiv  1\,\bmod D$ holds.
By the condition~\eqref{eq:chi} on $q$,
there exists an integer $z$ such that
$-pqz^2\equiv 1\,\bmod D$ holds.
Therefore we have an element
$$
znm\; \in \; A\conj{A}\; \subset \;A \;\subset\;  I\; \subset \; Q\inv I
$$
such that
$-pq(znm)^2\equiv 1\,\bmod D$.
Then the element $\mu:=znm/\sqD$ of 
$\DDD\inv Q\inv I$
satisfies $pqD|\mu|^2=-pq (znm)^2\equiv 1\,\bmod D$.
\end{proof}
%
%
%
\subsection{Proof of Proposition~\ref{prop:Jgenus}}
Let $J$ be a non-zero ideal in $\ZK$ such that
$d_J:=N(J)$ is prime to $D$.
Then we have
\begin{equation}\label{eq:JD}
\conj{J}\cap \DDD=\conj{J}\DDD.
\end{equation}
\begin{lemma}
Let $R$ be an element of $\setR$,
and let $RJ$ denote the left-ideal of $R$ generated  by $[J, 0]\subset R$.
Then we have
$$
(J\isomet RJ)\sperp\;=\; [0, Q\inv I_R\conj{J}],
$$
where $[0, Q\inv I_R\conj{J}]$ is the lattice 
such that the  underlying $\Z$-module is $Q\inv I_R\conj{J}\subset K$,
and that the bilinear form is given by $(x, y):=pq\Tr_{K/\Q}(x\conj{y})$.
\end{lemma}
\begin{proof}
For simplicity, we put
$M:=M_R$, $I:=I_R$, $f:=f_R$ and $\mu:=\mu_R$.
Since $J\otimes \Q=K$
and $B=[K,0]\perp [0, K]$ by~\eqref{eq:bilmatB},
we see that 
$(J\isomet RJ)\sperp$ is equal to $[0, K]\cap RJ$.
Let $\gamma, \gamma\sprime$ be a basis of $J$ as a $\Z$-module.
For an element $x\in K$, we have the following equivalence:
\begin{itemize}
\setlength{\itemsep}{3pt}
\item[] $[0, x]\in RJ$,
\item[$\Leftrightarrow$] 
there exist $[\alpha, \beta], [\alpha\sprime, \beta\sprime]\in R$ such that
$\alpha\gamma +\alpha\sprime\gamma\sprime=0$ and $\beta\,\conj{\gamma}+\beta\sprime\conj{\gamma\sprime}=x$,
\item[$\Leftrightarrow$]
there exist $\beta, \beta\sprime\in M$ and $a, a\sprime\in \ZK$ such that
$\beta\,\conj{\gamma}+\beta\sprime\conj{\gamma\sprime}=x$ and
$pq\sqD \mu (\conj\beta\gamma+\conj{\beta\sprime}\gamma\sprime)+a\gamma+a\sprime\gamma\sprime=0$,
\item[$\Leftrightarrow$]
$x\in \conj{J}M$ and $pq\sqD\mu\conj{x}\in J$.
\end{itemize}
Suppose that $x\in \conj{J}M$ and $pq\sqD\mu\conj{x}\in J$.
Then we have $f(x)=0$ by Lemma~\ref{lem:ft}, and hence 
$$
x\;\in\; \conj{J}M\cap \Ker(f)=\conj{J}\DDD M =Q\inv I\conj{J}
$$
by Lemma~\ref{lem:wtfisom}~and the equality~\eqref{eq:JD}.
Conversely, suppose that $x\in Q\inv I\conj{J}$.
Then we have $x\in \conj{J}M$.
On the other hand,  there exist $\xi, \xi\sprime\in Q\inv I$ such that
$x=\xi\,\conj{\gamma}+\xi\sprime\conj{\gamma\sprime}$.
Since $\xi, \xi\sprime \in \DDD M$, we have
$f(\xi)=f(\xi\sprime)=0$,
and hence
both of $pq\sqD\mu\conj{\xi}$ and $pq\sqD\mu\conj{\xi\sprime}$ are in $\ZK$.
Therefore we have
$$
pq\sqD\mu\conj{x}=(pq\sqD\mu\conj{\xi})\gamma+(pq\sqD\mu\conj{\xi\sprime})\gamma\sprime\in J,
$$
and thus $[0, x]\in RJ$ holds.
\end{proof}
We define an orientation of the $\Z$-module $Q\inv I_R\conj{J}\subset K$
by~\eqref{eq:oriinK}.
Then, for each $R\in\setR$,  we obtain  an oriented lattice
$[0, Q\inv I_R\conj{J}]$ of discriminant
$$
(pq)^2\cdot N(Q\inv I_R\conj{J})^2\cdot \disc(\ZK)=-p^2 d_J^2 D.
$$
On the other hand, recall that 
$\Psi ([Q\inv I_R\conj{J}])\in \ori{\lats}\primitive_D$ is represented by an oriented lattice
such that the  underlying $\Z$-module is $Q\inv I_R\conj{J}\subset K$,
and that the bilinear form is given by 
$$
(x, y)=\frac{1}{N(Q\inv I_R\conj{J})}\Tr_{K/\Q}(x\conj{y})=\frac{q}{d_J} \Tr_{K/\Q}(x\conj{y}).
$$
Therefore the isomorphism class of the oriented lattice
$(J\isomet RJ)\sperp= [0, Q\inv I_R\conj{J}]$
is equal to
$$
\Psi ([Q\inv I_R\conj{J}])[pd_J]\;\in\; \ori{\lats}_{p^2 d_J^2 D}.
$$
By Lemma~\ref{lem:surjR}, we have
$\shortset{[I_R]}{R\in \setR}=\Cl_D^2$.
Hence  Proposition~\ref{prop:ClD} implies that  the subset 
$\shortset{\Psi ([Q\inv I_R\conj{J}])}{R\in \setR}$ of $\ori{\lats}\primitive_D$
is a lifted genus $\ori{\GGG}$.
Consequently, the set 
$$
\set{[(J\isomet RJ)\sperp]}{R\in \setR}\;=\;\ori{\GGG}\,[pd_J]
$$
is also a lifted genus in $\ori{\lats}_{p^2 d_J^2 D}$.
Thus Proposition~\ref{prop:Jgenus} is proved.
%
%
%
%
\bibliographystyle{plain}
\def\cprime{$'$} \def\cprime{$'$} \def\cprime{$'$}

\end{document}